\documentclass[12pt]{amsart}
\usepackage{amsfonts, amssymb, amsmath, amsthm, eucal,
latexsym, array, pictex}
\usepackage{fullpage}
\usepackage{verbatim}

\input {cyracc.def}

\font\tencyr=wncyr10 
\font\tencyi=wncyi10 
\font\tencysc=wncysc10 

\def\rus{\tencyr\cyracc}
\def\rusi{\tencyi\cyracc}
\def\rusc{\tencysc\cyracc}

\newtheorem{thm}{Theorem}[section]
\newtheorem{conj}{Conjecture}[section]
\newtheorem{lm}[thm]{Lemma} 
\newtheorem{cl}[thm]{Corollary}
\newtheorem{prop}[thm]{Proposition}

\theoremstyle{remark}
\newtheorem{ex}{Example}[section]
\newtheorem{rmk}{Remark}[section]
\newtheorem{qn}{Question}[section]

\theoremstyle{definition}
\newtheorem{df}{Definition}[section]

\renewcommand{\iff}{if and only if }

\newcommand{\gt}{\mathfrak}

\newcommand{\GL}{{\rm GL}}

\newcommand{\id}{{\rm id}}
\newcommand{\ind}{{\rm ind\,}}

\newcommand{\rk}{\mathrm{rk\,}}

\newcommand{\Lie}{\mathrm{Lie\,}}
\newcommand{\Ker}{{\rm Ker\,}}

\newcommand{\gr}{\mathrm{gr\,}}
\renewcommand{\Im}{{\rm Im\,}}

\newcommand{\ad}{\mathrm{ad\,}}
\newcommand{\Ad}{{\rm Ad}}

\newcommand{\g}{\mathfrak g}
\newcommand {\ka}{{\mathfrak k}}
\newcommand {\el}{{\mathfrak l}}

\newcommand {\n}{{\mathfrak n}}

\newcommand {\p}{{\mathfrak p}}
\newcommand {\q}{{\mathfrak q}}
\newcommand {\es}{{\mathfrak s}}

\newcommand {\z}{{\mathfrak z}}

\newcommand {\gln}{{\mathfrak {gl}}_n}
\newcommand {\sln}{{\mathfrak {sl}}_n}

\newcommand {\spn}{{\mathfrak {sp}}_{2n}}

\newcommand {\son}{{\mathfrak {so}}_{n}}
\newcommand {\tri}{\mathfrak {sl}_2}

\newcommand {\cH}{{\mathcal H}}

\newcommand {\cN}{{\mathcal N}}
\newcommand {\co}{{\mathcal O}}
\newcommand {\cS}{{\mathcal S}}

\newcommand {\VV}{{\mathbb V}}
\newcommand {\BZ}{{\mathbb Z}}
\newcommand {\md}{/\!\!/}
\newcommand {\GR}[2]{{\textrm{{\bf #1}}}_{#2}}

\renewcommand{\le}{\leqslant}
\renewcommand{\ge}{\geqslant}

\newenvironment{E6}[6]{%
\begin{tabular}{@{}c@{}}
{#1}--{#2}--\lower3.3ex\vbox{\hbox{{#3}\rule{0ex}{2.5ex}}
\hbox{\hspace{0.4ex}\rule{.1ex}{1ex}\rule{0ex}{1.4ex}}\hbox{{#6}\strut}}--{#4}--{#5}
\end{tabular}}

\newenvironment{E7}[7]{%
\begin{tabular}{@{}c@{}}
{#1}--{#2}--{#3}--\lower3.3ex\vbox{\hbox{{#4}\rule{0ex}{2.5ex}}
\hbox{\hspace{0.4ex}\rule{.1ex}{1ex}\rule{0ex}{1.4ex}}\hbox{{#7}\strut}}--{#5}--{#6}
\end{tabular}}

\newfam\eusfam%
\font\euszw=eusm10 scaled 1200%
\font\eusac=eusm7 scaled 1200%
\font\eusacc=eusm7 scaled 1000%
\textfont\eusfam=\euszw\scriptfont\eusfam=\eusac%
\scriptscriptfont\eusfam=\eusacc%
\begin{document}

\vskip1ex

\title[Invariants of centralisers]{On symmetric invariants of centralisers \\
in reductive Lie algebras}
\author{D.\,Panyushev}
\address[D.P.] {Independent University of Moscow, Bol'shoi Vlasevskii per. 11,
119002, Moscow  Russia}
\email{panyush@mccme.ru}
\author{A.\,Premet}
\address[A.P.] {School of Mathematics, The University of Manchester, Oxford Rd, M13 9PL, UK.}
\email{sashap@maths.man.ac.uk}
\author{O.\,Yakimova}
\address[O.Y.] {Mathematisches Institut, Universit\"at zu K\"oln, Weyertal 86-90,
50931 K\"oln Germany}
\email{yakimova@mpim-bonn.mpg.de}
\thanks{D.P. and O.Y were supported in part by RFBI Grant 05-01-00988}
\begin{abstract}
Let $\g$ be a finite dimensional simple Lie algebra of rank $l$ over
an algebraically closed field of characteristic $0$. Let $e$ be a
nilpotent element of $\g$ and let $\g_e$ be the centraliser of $e$
in $\g$. In this paper we study the algebra $\cS(\g_e)^{\g_e}$ of
symmetric invariants of $\g_e$. We prove that if $\g$ is of type
$\bf A$ or $\bf C$, then $\,\cS(\g_e)^{\g_e}$ is always a graded
polynomial algebra in $l$ variables, and we show that this continues
to hold for {\it some} nilpotent elements in the Lie algebras of
other types. In type $\bf A$ we prove that the invariant algebra
$\cS(\g_e)^{\g_e}$ is freely generated by a regular sequence in
$\,\cS(\g_e)$ and describe the tangent cone at $e$ to the nilpotent
variety of $\g$.
\end{abstract}
\maketitle \maketitle


\vskip1ex

\vskip1ex

\vskip1ex

\vskip1ex

\vskip1ex

\vskip1ex

\vskip1ex

\vskip1ex

\tableofcontents

\section*{ Introduction}

\noindent 0.1. Let $\g$ be a finite-dimensional reductive Lie
algebra of rank $l$ over an algebraically closed field $\mathbb K$
of characteristic zero, and let $G$ be the adjoint group of $\g$.
Let ${\cN}(\g)$ denote the nilpotent cone of $\g$, i.e., the set of
all nilpotent elements of $\g$. Given $x\in\gt g$ we denote by
$\g_x$ and $G_x$ the centraliser of $x$ in $\gt g$ and $G$,
respectively. It is well-known that $\gt
g_x=\Lie\,G_x=\Lie\,G_x^\circ$ (here and in what follows $H^\circ$
stands for the identity component of an algebraic group $H$).

Inspired by a conversation with J.~Brundan at the Oberwolfach
meeting on enveloping algebras in March 2005, the second author
put forward the following conjecture:

\begin{conj}  \label{conj:main}
 For any $x\in\g$
the invariant algebra $\cS(\g_x)^{\g_x}$ is a graded polynomial
algebra in $l$ variables.
\end{conj}


In order to prove (or disprove) Conjecture~\ref{conj:main} it
suffices to consider the case where $\g$ is simple and
$x\in\cN(\g)$. The conjecture is known to hold for some
$x\in\cN(\g)$. For example, when $x=0$, it is an immediate
consequence of the Chevalley Restriction Theorem. At the other
extreme, when $x\in\cN(\g)$ is regular, the centraliser $\g_x$ is
abelian of dimension $l$ and we have $\cS(\g_x)^{\gt
g_x}\,=\,\cS(\g_x)\,\cong\,{\mathbb K}[X_1,\ldots,X_l]$ with $\deg
X_i=1$ for all $i$.

Conjecture~\ref{conj:main} is closely related to an earlier
conjecture of A.~Elashvili (initiated by a question of A.~Bolsinov).
Recall that the {\it index} of a finite-dimensional Lie algebra $\gt
s$ over $\mathbb K$, denoted $\ind \gt s$, is defined as the minimal
dimension of the stabilisers of linear functions on $\gt s$. In
other words, $\ind \gt s\,=\,\min\, \{\dim\,\gt s^f\,|\,\,f\in \gt
s^*\}$ where $\gt s^f=\{x\in\gt s\,|\,\,f([x,\gt s])=0\}$.
Elashvili's conjecture states that
\[
\ind \g_x=\,l=\,\rk\g  \qquad (\forall\, x\in\g).
\]
According to Vinberg's inequality, $\ind \g_x\ge l$ for all
$x\in\g$ (see \cite[1.6 \& 1.7]{dima2}, but the equality is {\it
much} harder to establish.

During the last decade Elashvili's conjecture drew attention of
several Lie theorists. Similar to Conjecture~\ref{conj:main} it
reduces to the case in which $\g$ is simple and $x\in\cN(\g)$. For
the spherical nilpotent orbits, Elashvili's conjecture was proved in
\cite{dima2} and \cite{dima3} by the first author. For $\g$
classical, Elashvili's conjecture was recently proved in \cite{fan}
by the third author. In 2004, J.-Y.~Charbonnel published a case-free
proof of Elashvili's conjecture applicable to all simple Lie
algebras; see \cite{char}. Unfortunately, the argument in
\cite{char} has a gap in the final part of the proof, which was
pointed out by L.~Rybnikov. At present we are unable to close this
gap. Answering a question of Elashvili, W.~de Graaf used a computer
programme to verify the conjecture for all nilpotent elements in the
Lie algebra of type ${\bf E}_8$ (unpublished). Since there are many
nilpotent orbits in the Lie algebras of exceptional types, it is
difficult to present the results of such computations in a concise
way.

To summarise, Elashvili's conjecture holds for the Lie algebras of
type $\mathbf A$, $\mathbf B$, $\mathbf C$, $\mathbf D$ and
${\mathbf G}_2$ and remains a challenging open problem for the Lie
algebras of type $\mathbf E$ and ${\mathbf F}_4$. We feel that it
would be very important to find a conceptual proof of Elashvili's
conjecture applicable to all finite-dimensional simple Lie algebras.

\smallskip

\noindent 0.2. The main goal of this paper is to prove
Conjecture~\ref{conj:main} for all nilpotent elements in the Lie
algebras of type $\mathbf A$ and $\mathbf C$. Our methods also work
for some nilpotent elements in the Lie algebras of type $\mathbf B$
and $\mathbf D$ and for a few nilpotent orbits in the exceptional
Lie algebras.

From now on, we fix a nonregular element $e\in\cN(\g)\setminus\{0\}$
and include it into an $\gt{sl}_2$-triple $(e,h,f)$ of $\g$. Let
$(\,\cdot\,,\,\cdot\,)$ denote the the scalar multiple of the
Killing form of $\g$ such that $(e,f)=1$, and put $\chi
=(e,\,\cdot\,)$. The map $\kappa$ from $\gt g$ to $\g^*$ which takes
$x$ to $(x,\,\cdot\,)$ extends uniquely to a $G$-equivariant
isomorphism between the symmetric algebra $\cS(\g)$ and the
coordinate algebra ${\mathbb K}[\g]$ of $\g$. This isomorphism of
graded algebras will be denoted by the same letter $\kappa$ and
referred to as a {\it Killing isomorphism}. The $G$-equivariance of
$(\,\cdot\,,\,\cdot\,)$ implies that $\g_e=[e,\g]^\perp$. On the
other hand, $\g=[e,\g]\oplus\g_f$ by the $\gt{sl}_2$-theory. It
follows that the Killing isomorphism $\kappa$ induces an algebra
isomorphism
\[
\kappa_e: \cS(\g_e)\stackrel{\sim}{\longrightarrow}\,{\mathbb
K}[\g_f],\quad x\mapsto (x,\,\cdot\,)_{\vert{\g_f}}\quad\ \
\,(\forall\,x\in\g_e).
\]
The coordinate algebra ${\mathbb K}[\g_f]$ carries a natural ${\mathbb
Z}$-grading in which the linear forms on $\g_f$ have degree $1$.
Each nonzero $\varphi\in {\mathbb K}[\g_f]$ is expressed uniquely
as
\[
\varphi=\varphi_k+\,\mbox{terms of higher degree},
\]
where $\varphi_k$ is a nonzero homogeneous element of degree
$k=k(\varphi)$. We say that $\varphi_k$ is the {\it initial term}
of $\varphi$, written $\varphi_k={\rm in}(\varphi)$. For
$\varphi=0$ we set ${\rm in}(\varphi)=0$.

Let ${\mathcal S}_e$ denote the {\it Slodowy slice} $e+\g_f$ at $e$
through the adjoint orbit $G\cdot e$. The translation map $x\mapsto
e+x$ induces an isomorphism of affine varieties
$\tau:\g_f\stackrel{\sim}{\rightarrow}{\mathcal S}_e$. The
comorphism $\tau^*$ maps the coordinate algebra ${\mathbb
K}[{\mathcal S}_e]$ isomorphically onto ${\mathbb K}[\g_f]$.

Let $F$ be a homogeneous element in $\cS(\g)$. Then
$\kappa(F)\in {\mathbb K}[\g]$ and
$\kappa(F)_{\vert\,{\mathcal S}_e}\in{\mathbb K}[{\mathcal S}_e]$.
The above discussion shows
that $\tau^*(\kappa(F)_{\vert\,{\mathcal S}_e})\in{\mathbb K}[\g_f]$
and $\kappa_e^{-1}\big({\rm
in}(\tau^*(\kappa(F)_{\vert\,{\mathcal S}_e}))\big)\in
\cS(\g_e)$. We now put
\[
^{e\!}F := \kappa_e^{-1}\big({\rm in}(\tau^*(\kappa(F)_{\vert\,{\mathcal S}_e}))\big).
\]
Thus, to each homogeneous $F\in \cS(\g)$ we assign a homogeneous
element $^{e\!}F\in\cS(\g_e)$. Roughly speaking, $^{e\!}F$ is the
initial component of $F_{|\kappa(\cS_e)}$.

\begin{prop}           \label{invariants}
If $F$ is a homogeneous element of $\cS(\g)^G$, then
$^{e\!}F\in\cS(\g_e)^{G_e}$.
\end{prop}
We give two proofs of Proposition~\ref{invariants}. The first
proof relies in a crucial way on some properties of the
quantisation  of the coordinate algebra ${\mathbb K}[{\mathcal
S}_e]$ introduced in \cite{Sasha1} (see also \cite{gg}). The
second (elementary) proof is given in the Appendix.

\smallskip

\noindent 0.3. Of particular interest are those homogeneous
generating sets $\{F_1,\ldots, F_l\}\subset\cS(\g)^{\g}$ for which
the resulting systems ${^{e\!}F}_1,\ldots,{^{e\!}F}_l$ are
algebraically independent. In Section~2 we show that if Elashvili's
conjecture holds for $\g_e$, then for any homogeneous system of
basic invariants $F_1,\ldots, F_l$ in $\cS(\g)^\g$ we have the
inequality
\begin{eqnarray}\label{sum-deg}
\textstyle{\sum}_{i=1}^l\deg{^e\!}F_i \le
(\dim\g_e+\rk\g)/2.\end{eqnarray} Furthermore, $^{e\!}F_1,\ldots,
{^{e\!}F}_l$ are algebraically independent in $\cS(\g_e)$ if and
only if the equality holds in (\ref{sum-deg}), that is
$\textstyle{\sum}_{i=1}^l\deg{^e\!}F_i = (\dim\g_e+\rk\g)/2$. If
this happens, we say that the system $F_1,\ldots,F_l$ is {\it
good} for $e$.

Given a linear function $\gamma$ on $\g_e$ we denote by
$\g_e^\gamma$ the stabiliser of $\gamma$ in $\g_e$ and set
$$(\g_e^*)_{\rm sing}:=\,\{\gamma\in\g_e^*\,|\,\, \dim\g_e^\gamma
>\ind\g_e\}.$$ The complement $\g_e^*\setminus(\g_e^*)_{\rm sing}$
consists of all {\it regular} linear functions of $\g_e$. We prove
in Section~2 that if Elashvili's conjecture holds for $\g_e$, then
for any good generating set $\{F_1,\ldots, F_l\}\subset\cS(\g)^\g$
the differentials $d_\gamma ({^{e\!}}F_1),\ldots,
d_\gamma({^{e\!}}F_l)$ are linearly independent at $\gamma\in\gt
g_e^*$ if and only if $\gamma$ is regular in $\gt g_e^*$. When
$e=0$, this is a classical result of Lie Theory often referred to as
{\it Kostant's differential criterion for regularity} (note that any
homogeneous generating system in $\cS(\g)^\g$ is good for $e=0$ and
Elashvili's conjecture is true in this case). When $e$ is regular
nilpotent, the statement follows from another theorem of Kostant
saying that the restriction of the adjoint quotient map to the
Slodowy slice $\cS_e$ is an isomorphism of algebraic varieties.
Beyond these two extreme cases our result seems to be new. It should
be stressed, however, that if $\g$ is not of type $\mathbf A$ or
$\mathbf C$, then there may exist nilpotent elements in $\g$ which
do not admit good generating systems in $\cS(\g)^\g$. One such
element in $\g=\mathfrak{so}_{12}$ is exhibited in
Example~\ref{ex-so}. Quite surprisingly, the root vectors in Lie
algebras of type ${\mathbf E}_8$ provide yet another example of this
kind.

\smallskip

\noindent 0.4. Our proof of the above results relies on some
geometric properties of Poisson algebras of Slodowy slices
(established in \cite{Sasha1} and \cite{gg}) and  a theorem of
Odesskii-Rubtsov \cite{or} on polynomial Poisson algebras with a
regular structure of symplectic leaves. All necessary background on
polynomial Poisson algebras is assembled in Section~1.

Let ${\mathcal A}={\mathbb K}[x_1,\ldots, x_n]$ be a polynomial
algebra in $n$ variables. For $g_1,\ldots, g_m\in\mathcal A$, we
denote by ${\mathcal J}(g_1,\ldots,g_m)$ the Jacobian locus of
$g_1,\ldots, g_m$, i.e., the set of all $\xi\in{\rm
Specm}\,\mathcal A$ for which the differentials $d_\xi
g_1,\ldots,d_\xi g_m$ are linearly dependent. Suppose  ${\mathcal
A}$ is a Poisson algebra and let $\pi\in\mathrm{Hom}_{\mathcal
A}\,(\Omega^2({\mathcal A}),{\mathcal A})$ be the corresponding
Poisson bivector. Let $Z({\mathcal A})$ denote the Poisson centre
of ${\mathcal A}$. The defect of the skew-symmetric matrix
$\big(\{x_i,x_j\}\big)_{1\le i,j\le n}$ with entries in $\mathcal
A$ is called the {\it index} of ${\mathcal A}$ and denoted
$\ind\mathcal A$. It is well-known (and easily seen) that
$\mathrm{tr.\, deg}_{\mathbb K}\, Z(\mathcal A)\le\ind\mathcal A.$
We denote by $\mathrm{Sing}\,\pi$ the set of all $\xi\in{\rm
Specm}\,{\mathcal A}$ for which $\rk\pi(\xi)<n-\ind\mathcal A$. A
subset $\{Q_1,\ldots\, Q_l\}\subset Z({\mathcal A})$ is said to be
{\it admissible\/} if $l=\ind\mathcal A$ and the Jacobian locus
${\mathcal J}(Q_1,\ldots,Q_l)$ has codimension $\ge 2$ in $\mathbb
A^n$. We say that $({\mathcal A}, \pi)$ is a {\it quasi-regular}
Poisson algebra if $Z(\mathcal A)$ contains an admissible subset
and $\mathrm{Sing}\, \pi$ has codimension $\ge 2$ in ${\rm
Specm}\,\mathcal A$.

Assume now that ${\mathcal A}=\mathbb K[x_1,\ldots, x_n]$ is graded
and each $x_i$ is homogeneous of positive degree. Let $f_1,\ldots,
f_s$ be a collection of homogeneous elements in $\mathcal A$ such
that the Jacobian locus ${\mathcal J}(f_1,\ldots, f_s)$ has
codimension $\ge 2$ in ${\rm Specm}\, \mathcal A$, and denote by $R$
the subalgebra of $\mathcal A$ generated by $f_1,\ldots, f_s$.
Inspired by Skryabin's result \cite[Theorem~5.4]{sk} on modular
invariants of finite group schemes we prove that if an element
$\tilde{f}\in\mathcal A$ is algebraic over $R$, then necessarily
$\tilde{f}\in R$. This has the following consequence:
\begin{thm}\label{homo} Let ${\mathcal
A}=\mathbb K[x_1,\ldots,x_n]$ be a quasi-regular Poisson algebra of
index $l$ and suppose that $\mathcal A=\bigoplus_{k\ge 0}{\mathcal
A}(k)$ is graded in such a way that $x_i\in {\mathcal A}(r_i)$ for
some $r_i>0$, where $1\le i\le n$. Suppose further that $Z({\mathcal
A})$ contains an admissible set $\{Q_1,\ldots,Q_l\}$ consisting of
homogeneous elements of $\mathcal A$. Then $Z(\mathcal A)=\,{\mathbb
K}[Q_1,\ldots, Q_l]$.
\end{thm}

\noindent 0.5. In this paper, we mostly apply Theorem~\ref{homo} to
the pair $({\mathcal A},\pi)=\big(\cS(\g_e),\pi_e^{PL}\big)$ where
$\pi_e^{PL}$ is the Poisson bivector of $\cS(\g_e)$ induced by the
Lie bracket of $\g_e$. In this situation $Z(\mathcal
A)=\cS(\g_e)^{\g_e}$. (One noteworthy application of
Theorem~\ref{homo} to quantisations of Slodowy slices can be found
in Remark~\ref{slodowy}.) Of course, before applying
Theorem~\ref{homo} to the pair $\big(\cS(\g_e),\pi_e^{PL}\big)$ we
have to make sure that our nilpotent element qualifies. That is to
say, we must check that $e$ admits a good generating system
$F_1,\ldots, F_l$, that Elashvili's conjecture holds for $\g_e$, and
that ${\mathcal J}({^e\!}F_1,\ldots, {^e\!}F_l)=(\g_e^*)_{\rm sing}$
has codimension $\ge 2$ in $\g_e^*$. Our main result is the
following:

\begin{thm}\label{inv}
Suppose $e$ admits a  good generating system $F_1,\ldots, F_l$ in
$\,\cS(\g)^\g\,$ and assume further that Elashvili's conjecture
holds for $\g_e$ and $(\g_e^*)_{\rm sing}$ has codimension $\ge 2$
in $\g_e^*$. Then $\cS(\g_e)^{\g_e}=\,\cS(\g_e)^{G_e}$ is a
polynomial algebra in $\,{^e\!}F_1,\ldots, {^e\!}F_l$.
\end{thm}

Suppose $\g$ is of type ${\mathbf A}_n$ or ${\mathbf C}_n$, where
$n\ge 2$, and let $e\in\cN(\g)$. By \cite{fan},  Elashvili's
conjecture holds for $\g_e$. In Section~3 we show that the singular
locus $(\g_e^*)_{\rm sing}$ has codimension $\ge 2$ in $\g_e^*$,
whilst our results in Section~4 imply that in types $\mathbf A$ and
$\mathbf C$ the invariant algebra $\cS(\g)^\g$ contains a
homogeneous generating set which is good for {\it all} nilpotent
elements in $\g$ (this is no longer true in types $\mathbf B$ and
$\mathbf D$). Applying Theorem~\ref{inv} we are able to conclude
that Conjecture~\ref{conj:main} holds for all nilpotent elements in
$\g$.

Apart from the the above-mentioned results, we show in Sections~3
and 4 that the conditions of Theorem~\ref{inv} are satisfied for
some nilpotent elements in Lie algebras of types $\mathbf B$ and
$\mathbf D$. Subsection~3.9 illustrates the behavior of the simple
Lie algebras $\g$ of types other than $\bf A$ and $\bf C$ by
producing a nilpotent element $e\in\g$ for which $(\g_e^*)_{\rm
sing}$ has codimension $1$ in $\g_e^*$.

\smallskip

\noindent 0.6. In Section~5 we study the null-cone $\cN(e)$ of
$\g_e^*$, that is the subvariety of $\g_e^*$ consisting of all
linear functions $\xi$ such that $\varphi(\xi)=0$ for all
$\varphi\in\cS(\g_e)^{\g_e}$ with $\varphi(0)=0$. Here we have to
assume that $\g=\mathfrak{gl}_n$. Working with the good generating
set $\{F_1,\ldots, F_n\}\subset \cS(\g)^\g$ mentioned in (0.5) we
show that the zero locus $\cN(e)$ of the ideal $({^e\!}F_1,\ldots,
{^e\!}F_n)$ has codimension $n$ in $\g_e^*$ and hence
${^e\!}F_1,\ldots, {^e\!}F_n$  is a regular sequence in $\cS(\g_e)$.
As a consequence, we describe the tangent cone at $e$ to the variety
of all nilpotent $n\times n$ matrices over $\mathbb K$; see
Corollary~\ref{cone}. Although the variety $\cN(e)$ is irreducible
in some interesting cases, in general it has many irreducible
components. The problem of describing the irreducible components of
$\cN(e)$ for $\g=\mathfrak{gl}_n$ is wide open.

\smallskip

\noindent 0.7. Let $\tilde{e}\in\co_{\min},$ where $\co_{\min}$ is
the minimal (nonzero) nilpotent orbit in $\g$. The element
$\tilde{e}$ is $G$-conjugate to a highest root vector in $\g$.
Recall that outside type $\mathbf A$ the orbit $\co_{\min}$ is {\it
rigid}, i.e., cannot be obtained by Lusztig--Spaltenstein induction
from a nilpotent orbit in a Levi subalgebra of $\g$. We put
Conjecture~\ref{conj:main} to the test by investigating the
invariant algebra $\cS(\g_{\tilde{e}})^{\g_{\tilde{e}}}$. Here our
result is as follows:
\begin{thm}\label{min-nil}
Suppose $\rk\g\ge 2$. Then the singular locus
$(\g_{\tilde{e}}^*)_{\rm sing}$ has codimension $\ge 2$ in
$\g_{\tilde{e}}^*$. If $\g$ is not of type ${\bf E}_8$, then
$\tilde{e}$ admits a good generating system in $\cS(\g)^\g$  and the
invariant algebra $\cS(\g_{\tilde{e}})^{\g_{\tilde{e}}}$ is
isomorphic to a graded polynomial algebra in $\rk \g$ variables. The
degrees of basic invariants of
$\cS(\g_{\tilde{e}})^{\g_{\tilde{e}}}$ are given in Table~\ref{T1}.
\end{thm}
\begin{table}
\begin{tabular}{|l|l|}
\hline \ {\rm Type of} $\g$\ \ &\ {\rm Degrees of basic invariants}\ \ \\
\hline
$\ {\mathbf A}_n,\, n\ge 1$ & $\ 1$, $2$, $\ldots\,$, $n$\\
\hline

${\mathbf B}_n,\, n\ge 3$ & $\ 1$, $3$, $4$, $\ldots\,$,  $2n-2$\\
\hline
${\mathbf C}_n,\, n\ge 2$ & $\ 1$, $3$, $\ldots\,$,  $2n-1$\\
\hline
${\mathbf D}_n,\, n\ge 4$ & $\ 1$, $3$, $4$, $\ldots\,$,  $2n-4$, $n-1$\\
\hline
${\mathbf E}_6$ & $\ 1$, $4$, $4$,  $6$, $7$, $9$\\
\hline
${\mathbf E}_7$ & $\ 1$, $4$, $6$,  $8$, $9$, $11$, $14$\\
\hline
${\mathbf F}_4$ & $\ 1$, $4$, $6$,  $9$\\
\hline
${\mathbf G}_2$ & $\ 1$, $4$\\
\hline
\end{tabular}

\bigskip

\caption{} \label{T1}
\end{table}

\smallskip

In order to prove Theorem~\ref{min-nil} for Lie algebras of types
${\mathbf E}_7$ we have to use the explicit system of basic
invariants for the Weyl group of type $\mathbf E_7$ constructed in
\cite{KM}. In type ${\bf E}_8$, we reduce Conjecture~\ref{conj:main}
for $\g_{\tilde{e}}$ to a specific problem on polynomial invariants
of the Weyl group of type ${\bf E}_7$; see Theorem~\ref{E8-4}. In
principle, this problem can be tackled by computational methods.

\smallskip

We adopt the Vinberg--Onishchik numbering of simple roots and
fundamental weights in simple Lie algebras; see \cite[Tables]{VO}.
The $i$-th fundamental weight is denoted by $\varpi_i$.

\smallskip

After posting this paper on the {\tt arXiv} we have learned from
Jonathan Brundan that he and Jonathan Brown also proved, for
$\g=\mathfrak{gl}_n$, that the invariant algebra $\cS(\g_e)^{\g_e}$
is free for any nilpotent element $e\in\g$; see \cite{BB}. The
approach in \cite{BB} is different from ours; it relies on the
earlier work of Brundan--Kleshchev \cite{BK} and employs an argument
in the spirit of \cite{rt}. The main goal of \cite{BB} is to
construct an explicit set of {\it elementary invariants} that
generate the centre of the universal enveloping algebra $U(\g_e)$.
Using this generating set it is not difficult to confirm
Conjecture~\ref{explicit} of this paper.

\section{Some general results}
\noindent 1.1. Our goal in this section is twofold: to prove an
extended characteristic-zero version of Skryabin's theorem
\cite{sk} on invariants of finite group schemes and to obtain a
slight generalisation of a result of Odesskii--Rubtsov~\cite{or}
on polynomial Poisson algebras. We first recall some basics on the
classical duality between differential forms and polyvector
fields.

Let $\mathbb A^n=\mathbb A^n_{\mathbb K}$ be the $n$-dimensional
affine space with the algebra of regular functions ${\mathcal
A}=\mathbb K[x_1,\ldots,x_n]$. Let $W$ denote the derivation algebra
of $\mathcal A$. This is a free $\mathcal A$-module with basis
consisting of partial derivatives $\partial_1,\ldots,
\partial_n$ with respect to $x_1,\ldots, x_n$.
Let $\Omega^1=\,\mathrm{Hom}_{\mathcal A}(W,{\mathcal A})$ and let
$\Omega=\,\bigoplus_{k=0}^n\,\Omega^k$ be the exterior $\mathcal
A$-algebra on $\Omega^1$. The exterior differential $d\colon
{\mathcal A}\rightarrow\Omega^1,\,\, (df)(D)=D(f),$ extends uniquely
up to a zero-square graded derivation of the $\mathcal A$-algebra
$\Omega$. We identify $\Omega^0$ with $\mathcal A$ and regard
$\Omega^1$ as the $\mathcal A$-module of global sections on the
cotangent bundle $T^*\mathbb A^n$. Note that $\Omega^k$ is a free
$\mathcal A$-module with basis $\{dx_{i_1}\wedge\ldots\wedge
dx_{i_k}\,\vert\,\,1\le i_1<\cdots< i_k\le n\}$.

We view the exterior powers $\Omega^k=\bigwedge_{\mathcal
A}^k\Omega^1$ and $\bigwedge_{\mathcal A}^kW$ as dual $\mathcal
A$-modules by using the nondegenerate $\mathcal A$-pairing
$$
\langle\alpha_1\wedge\ldots\wedge\alpha_k,D_1\wedge\ldots\wedge
D_k\rangle=\det \big(\alpha_i(D_j)\big).
$$
For $\eta\in\Omega^k$,
set $\eta(D_1\wedge\ldots\wedge D_k):=\langle\eta,D_1\wedge\ldots\wedge D_k\rangle$. For $D\in
\bigwedge_{\mathcal A}^k W$, set
$D(\alpha_1\wedge\ldots\wedge\alpha_k):=\langle\alpha_1\wedge\ldots\wedge\alpha_k,
D\rangle$. Then for $D\in\bigwedge_{\mathcal A}^pW=(\Omega^p)^*$
and $D'\in\bigwedge_{\mathcal A}^qW=(\Omega^q)^*$ we have
\begin{eqnarray*}(D\wedge
D')(\alpha_1\wedge\ldots\wedge\alpha_{p+q})&=& \langle
\alpha_1\wedge\ldots\wedge\alpha_{p+q},D\wedge D'\rangle\\
&=& \sum(\mathrm{sgn}\,\sigma)
D(\alpha_{\sigma(1)},\ldots,\alpha_{\sigma(p)})D'(\alpha_{\sigma(p+1)},\ldots,\alpha_{\sigma(p+q)}),
\end{eqnarray*}
where the summation runs over the set of all permutations $\sigma$
of $\{1,\ldots, p+q\}$ which are increasing on $\{1,\ldots, p\}$
and $\{p+1,\ldots, p+q\}$.

For $X\in\bigwedge_{\mathcal A} W$ and $\xi\in\mathbb A^n$, the
specialisation $X_\xi$ is a well-defined element of the exterior
algebra $\bigwedge T_\xi(\mathbb A_n)$ on the tangent space
$T_\xi(\mathbb A^n)$. For $X\in W$, the {\it left interior product}
$i_X$ is the unique $\mathcal A$-linear endomorphism of degree $-1$
on $\Omega$ such that
$$i_X(\eta)(D_1\wedge\ldots\wedge D_{k})=\eta(X\wedge
D_1\ldots\wedge D_k)\qquad\
\big(\forall\,\eta\in\Omega^{k+1}\big).$$ For $\omega\in
\Omega^1$, the {\it right interior product} $j_\omega$ is the
unique $\mathcal A$-linear endomorphism of degree $-1$ on
$\bigwedge_{\mathcal A}W$ such that
$$j_\omega(D)(\alpha_1\wedge\ldots\wedge \alpha_{k})=D(\alpha_1\wedge
\ldots\wedge \alpha_k\wedge\omega)\qquad\
\big(\forall\,D\in\textstyle{\bigwedge_{\mathcal
A}^{k+1}W}\big).$$ Using the above discussion it is easy to
observe that the endomorphisms $i_X$ and $j_\omega$ are graded
derivations (a.k.a. super-derivations) of $\Omega$ and
$\bigwedge_{\mathcal A}W$, respectively. More generally, given
$X\in\bigwedge_{\mathcal A}^pW$ and $\omega\in\Omega^p$ one
defines the {\it right interior product} $i_X$ and the {\it left
interior product} $j_\omega$ to be the unique endomorphisms of
degree $-p$ on $\Omega$ and $\bigwedge_{\mathcal A}\,W$,
respectively, such that
$$\langle i_X(\eta), D \rangle=\langle\eta,X\wedge D\rangle\ \
\mbox{  and  }\ \, \langle \eta, j_\omega(D)
\rangle=\langle\eta\wedge\omega, D\rangle\qquad
\big(\forall\,D\in\textstyle{\bigwedge_{\mathcal A}^pW},\,\,
\forall\,\eta\in\Omega\big).$$ The mappings $X\mapsto i_X$ and
$\omega\mapsto j_\omega$ then give rise to $\mathcal A$-algebra
homomorphisms $i\colon \bigwedge_{\mathcal
A}W\rightarrow\rm{End}(\Omega)^{\rm op}$ and
$j\colon\Omega\rightarrow \mathrm{End}\big(\bigwedge_{\mathcal
A}\,W\big)$. In other words, we have $i_X\circ\, i_{Y}=i_{Y\wedge
X}$ and $j_\alpha\circ j_\beta=j_{\alpha\wedge\beta}$ for all
$X,Y\in\bigwedge_{\mathcal A}W$ and all $\alpha,\beta\in\Omega$.
Finally, $i_X(\omega)=j_\omega(X)=\langle\omega,X\rangle$ whenever
$X\in\bigwedge_{\mathcal A}^p W$ and $\omega\in\Omega^p$.

The top components $\Omega^n$ and $\bigwedge_{\mathcal A}^nW$ are
free modules of rank $1$ over $\mathcal A$ generated by
$dx_1\wedge\ldots\wedge dx_n$ and
$\partial_1\wedge\ldots\wedge\partial_n$, respectively. The
mappings $X\mapsto i_X(dx_1\wedge\ldots\wedge dx_n)$ and
$\omega\mapsto j_\omega(\partial_1\wedge\ldots\wedge
\partial_n)$ induce canonical $\mathcal A$-module isomorphisms
$\bigwedge_{\mathcal A}^p W\cong \,\Omega^{n-p}$ and $\Omega
^p\,\cong \bigwedge_{\mathcal A}^{n-p}W$.

\smallskip

\noindent 1.2. For $g_1,\ldots, g_m\in\mathcal A,$ the {\it Jacobian
locus} ${\mathcal J}(g_1,\ldots,g_m)$ consists of all
$\xi\in{\mathbb A}^n$ for which the differentials $d_\xi
g_1,\ldots,d_\xi g_m$ are linearly dependent. The set ${\mathcal
J}(g_1,\ldots,g_m)$ is Zariski closed in $\mathbb A^n$ and it
coincides with $\mathbb A^n$ if and only if $g_1,\ldots, g_m$ are
algebraically dependent. Our interpretation of Skryabin's result
\cite[Theorem~5.4]{sk} will be based on the following theorem which
is of independent interest:

\begin{thm}\label{skryabin} Suppose ${\mathcal A}=\mathbb K[x_1,\ldots, x_n]$ is
graded in such a way that each $x_i$ is homogeneous of positive
degree. Let $R$ be the subalgebra of $\mathcal A$ generated by
homogeneous elements $f_1,\ldots, f_s$ and assume further that
${\mathcal J}(f_1,\ldots, f_s)$ has codimension $\ge 2$ in $\mathbb
A^n$. Then $R$ is algebraically closed in $\mathcal A$. In other
words, if $\tilde{f}\in \mathcal A$ is algebraic over $R$, then
$\tilde{f}\in R$.
\end{thm}
\begin{proof}
For $t\in\mathbb K^{^\times},$ we denote by $\rho(t)$ the
automorphism of $\mathcal A$ such that $\rho(t){\cdot f}=t^kf$ for
all $f\in{\mathcal A}(k)$, where  ${\mathcal A}(k)$ is the $k$-th
graded component of $\mathcal A$. Let ${\mathcal Q}(R)$ be the field
of fractions of $R$, a subfield of ${\mathbb K}(x_1,\ldots, x_n),$
and denote by $\tilde{R}$ the algebraic closure of $R$ in $\mathcal
A$. Since $\tilde{R}$ is nothing but the intersection of $\mathcal
A$ with the algebraic closure of ${\mathcal Q}(R)$ in ${\mathbb
K}(x_1,\ldots, x_n)$, it is a subalgebra of $\mathcal A$. Since all
$f_i$ are homogeneous, the subalgebra $R$ is $\rho(\mathbb
K^{^\times})$-stable. But then so is $\tilde{R}$. As a consequence,
$\tilde{R}$ is a homogeneous subalgebra of $\mathcal A$. Thus, in
order to prove the theorem it suffices to show that if a {\it
homogeneous} element $\tilde{f}\in\mathcal A$ is algebraic over $R$,
then $\tilde{f}\in R$.

We shall argue by induction on the degree of $\tilde{f}$. So assume
that the statement holds for all homogeneous elements of degree less
than $\deg\tilde{f}$ (when $\deg \tilde{f}=1$, this is a valid
assumption).

\smallskip

\noindent (a) The grading of $\mathcal A$ induces that on the
$\mathbb K$-algebra $\Omega$ where we impose that $\deg dx_i=\deg
x_i$. Note that $a\in\mathcal A$ is algebraic over $R$ if and only
if $da\wedge df_1\wedge\ldots\wedge df_s=0$ in $\Omega$. Since
$\mathcal{J}(g_1,\ldots, g_m)$ consists of all $\xi\in\mathbb A^n$
for which $d_\xi g_1\wedge\ldots\wedge d_\xi g_m=0$, our
assumption on $f_1,\ldots, f_s$ implies that for every subset
$\{i_1,\ldots, i_k\}$ of $\{1,\ldots, s\}$ the locus ${\mathcal
J}(f_{i_1},\ldots, f_{i_k})$ has codimension $\ge 2$ in $\mathbb
A^n$. From this it follows that passing to smaller subsets of
$\{f_1,\ldots,f_s\}$ and renumbering if necessary we can reduce
our proof to the situation where for each $i$ the polynomials
$\{f_1,\ldots,f_{i-1},f_{i+1},\ldots,f_s,\tilde f\}$ are
algebraically independent. So let us assume from now that this is
the case, and put
$$T:=df_1\wedge\ldots\wedge df_s,\ \qquad
T_i:=df_1\wedge\ldots\wedge df_{i-1}\wedge d\tilde f\wedge
df_{i+1}\wedge\ldots\wedge df_s\qquad(1\le i\le s).$$ By our
assumption, $T$ and the $T_i$ are {\it nonzero} homogeneous
elements of $\Omega$.

\smallskip

 \noindent (b) If
$\xi\not\in\mathcal{J}(f_1,\ldots, f_s)$, then $d_\xi f_1,\ldots,
d_\xi f_s$ are linearly independent and $d_\xi\tilde f$ is a
linear combination of $d_\xi f_1,\ldots, d_\xi f_s$. It follows
that the specialisation of $T_i$ at $\xi$ is a scalar multiple of
$d_\xi f_1\wedge\ldots\wedge d_\xi f_s$. As $\Omega$ is a free
$\mathcal A$-module, this yields that $T$ and $T_i$ are linearly
dependent as elements of the ${\mathbb K}(x_1,\ldots, x_n)$-vector
space $\mathbb K(x_1,\ldots, x_n)\otimes_{\mathcal A}\Omega$.
Combined with our discussion in part~(a) this implies that
$a_iT_i=b_iT$ for some nonzero coprime $a_i,b_i\in\mathcal A$. As
${\mathcal J}(f_1,\ldots, f_m)$ has codimension $\ge 2$ in
$\mathbb A^n$, the function $a_i$ must be constant. Thus, $T_i=p_i
T$ where $p_i$ is a nonzero homogeneous element of the graded
algebra $\mathcal A$.

\smallskip

\noindent (c) Since $d^2=0$, we have $dp_i\wedge T=d(p_i
T)=dT_i=0$. Our remarks in part~(a) now show that all $p_i$ are
algebraic over $R$. Let $$F=S_k(X_1,\ldots,X_s)Y^k+
S_{k-1}(X_1,\ldots,X_s)Y^{k-1}+\cdots+ S_0(X_1,\ldots,X_s)$$ be a
nonzero polynomial in ${\mathbb K}[X_1,\ldots, X_s, Y]$ of minimal
possible degree in $Y$  such that $F(f_1,\ldots,f_s,\tilde f)=0$.
Assume further that $S_k$ has minimal possible total degree in
$\mathbb K[X_1,\ldots, X_s]$ and that all $S_i(f_1,\ldots, f_s)$
are homogeneous in the graded algebra $\mathcal A$. Applying the
exterior differential we get $0=\,dF(f_1,\ldots,f_s,\tilde
f)=\,\tilde\psi d\tilde f+\sum\psi_i df_i$\,\, where
\begin{eqnarray*}
\tilde\psi&=&k{\tilde f}^{k-1}S_k(f_1,\ldots, f_s)+\cdots+S_1(f_1,\ldots, f_s),\\
\psi_i&=&{\tilde f}^k\frac{\partial S_k}{\partial X_i}(f_1,\ldots,
f_s)+ {\tilde f}^{k-1}\frac{\partial S_{k-1}}{\partial
X_i}(f_1,\ldots, f_s)+
  \cdots+\frac{\partial S_0}{\partial X_i}(f_1,\ldots, f_s)\qquad\ (1\le i\le m).
\end{eqnarray*}
As $\tilde{\psi}\ne 0$ by our choice of $F$, we have $d\tilde
f=-\sum({\psi_i}/{\tilde\psi})df_i$. This  forces
$T_i=-({\psi_i}/{\tilde\psi}) T$ for all $i$. Then
$\psi_i=-p_i\tilde\psi$ by our concluding remark in part~(b).

\smallskip

\noindent (d) Part~(b) also shows that each $p_i$ is homogeneous
with $\deg p_i=\deg\tilde f-\deg f_i<\deg\tilde f$. Since all $p_i$
are algebraic over $\mathcal A$ by part~(c), our inductive
hypothesis implies that $p_i\in R$ for all $i$. We now look again at
the formulae displayed in part~(c), this time keeping in mind that
$\psi_i+p_i\tilde\psi=0$ and $p_i\in\mathbb K[f_1,\ldots,f_s]$.

If at least one of the partial derivatives $\partial S_k/\partial
X_i$ was nonzero, we would have a nontrivial polynomial relation
for $\tilde{f}, f_1,\ldots, f_s$ with a smaller total degree of
$S_k$. Due to our choice of $F$ this is impossible, however. So
$S_k$ is a nonzero constant, and there will be no harm in assuming
that $S_k=1$. Note that each equality $\psi_i+p_i\tilde\psi=0$ now
induces a polynomial relation for $\tilde{f}, f_1,\ldots, f_s$ of
degree $\le k-1$ in $Y$.  Since such a relation is trivial by our
choice of $F$, the coefficient $(\partial S_{k-1}/\partial
X_i)(f_1,\ldots, f_s)+kp_i$ of $\tilde{f}^{k-1}$ in the relation
has to be zero. In view of our remarks in part~(c) we thus obtain
$$dS_{k-1}(f_1,\ldots, f_s)=-\sum kp_i df_i=-kd\tilde f.$$ Then
$\tilde f=-S_{k-1}/k+\lambda$ for some $\lambda\in\mathbb K$,
which shows that $\tilde f\in R$.
\end{proof}

\smallskip

\noindent {1.3.} Now suppose that ${\mathcal A}$ possesses a
Poisson structure $\{\,\,,\,\}\colon\,{\mathcal A}\times{\mathcal
A}\to {\mathcal A}$ and let $\pi$ denote the corresponding {\it
Poisson bivector}, the element of $\mathrm{Hom}_{\mathcal
A}(\Omega^2,{\mathcal A})$ satisfying $\pi(df\wedge dg)=\{f,g\}$
for all $f,g\in{\mathcal A}$. In view of the duality described in
(1.1) we may assume that $\pi\in \bigwedge^2_{\mathcal A} W$, that
is
$$
\langle df\wedge dg,\pi\rangle=\,\{f,g\}\qquad \ (\forall\,
f,g\in{\mathcal A}).
$$
Let $\rk \pi(\xi)$ denote the rank of the skew-symmetric matrix
$\big(\{x_i,x_j\}\big)_{1\le i,j\le n}$ at $\xi\in \mathbb A^n$.
The {\it index} of the Poisson algebra ${\mathcal A}$, denoted
$\ind\mathcal A$, is defined as
$$\ind\mathcal A:=n-\max_{\xi\in \mathbb A^n}\,\, \rk\pi(\xi).$$
Let $Z({\mathcal A})$ denote the Poisson centre of ${\mathcal A}$
and put $\mathrm{Sing}\,\pi:=\,\{\xi\in\mathbb A^n\,|\,\,
\rk\pi(\xi)<n-\ind\mathcal A\}$. Clearly, $\mathrm{Sing}\, \pi$ is
a proper Zariski closed subset of $\mathbb A^n$. Note that
$\langle df\wedge dg,\pi\rangle=0$ for all $f\in Z({\mathcal A})$
and all $g\in{\mathcal A}$. Hence the linear subspace $\{d_\xi
f\mid f\in Z({\mathcal A})\}$ lies in the kernel of $\pi(\xi)$ and
we have $$ \mathrm{tr.\, deg}_{\mathbb K}\, Z(\mathcal
A)\le\ind\mathcal A.
$$

We say that a subset $\{Q_1,\ldots\, Q_l\}\subset Z({\mathcal A})$
is {\it admissible\/} if $l=\ind\mathcal A$ and the locus
${\mathcal J}(Q_1,\ldots,Q_l)$ has codimension $\ge 2$ in $\mathbb
A^n$. It is clear from the definition that any admissible subset
of $Z({\mathcal A})$ is algebraically independent.

\begin{df} We call a Poisson algebra $({\mathcal A}, \pi)$  {\it quasi-regular}
if the Poisson centre of $\mathcal A$ contains an admissible
subset and $\mathrm{Sing}\, \pi$ has codimension $\ge 2$ in
${\mathbb A}^n$.
\end{df}

Given $k\in\mathbb N$ we set
$$\pi^k:=\,\underbrace{\pi\wedge\pi\wedge\ldots\wedge \pi}_{k \
\scriptstyle{\mathrm{factors}}}\,,$$ an element of
$\bigwedge_{\mathcal A}^{2k}W$. The following is a slight
modification of \cite[Theorem 3.1]{or}.
\begin{thm} \label{or}
Let ${\mathcal A}=\mathbb K[x_1,\ldots,x_n]$ be a quasi-regular
Poisson algebra of index $l$ and  let $\{Q_1,\ldots,Q_l\}\subset
Z({\mathcal A})$ be an admissible set in $Z(\mathcal A)$. Then
\[
\pi^{(n-l)/2}=\, \lambda j_{dQ_1\wedge\ldots\wedge
dQ_l}(\partial_1\wedge\ldots\wedge\partial_n)
\]
for some nonzero $\lambda\in \mathbb K$.
\end{thm}
\begin{proof}
Set $w:=j_{dQ_1\wedge\ldots\wedge
dQ_l}(\partial_1\wedge\ldots\wedge\partial_n)$, an element of
$\bigwedge_{\mathcal A}^{n-l}W$. Since $j\colon\,\Omega\rightarrow
\mathrm{End}(\bigwedge_{\mathcal A}W)$ is an exterior algebra
homomorphism, it must be that $$j_{dQ_i}(w)=j_{dQ_i\wedge
dQ_1\wedge\ldots\wedge
dQ_l}(\partial_1\wedge\ldots\wedge\partial_n)=\,0\qquad\,(1\le i\le
l).$$ Since $Q_i\in Z({\mathcal A})$, we also have
$$\langle df, j_{dQ_i}(\pi)\rangle=\langle df\wedge dQ_i, \pi \rangle=
\{f,Q_i\}=\,0\qquad\, (\forall f\in\mathcal A).$$ Hence
$j_{dQ_i}(\pi)=0.$ Since $j_{dQ_i}$ is a graded derivation of
$\bigwedge_{\mathcal A}W$, it follows that $j_{dQ_i}(\pi^{k})=\,0$
for all $k\in\mathbb N$. As a consequence,
$j_{dQ_i}\big(\pi^{(n-l)/2}\big)=\,j_{dQ_i}(w)=0$ for all $i\le
l$. As $l=\ind\mathcal A$, we have $\pi^{(n-l)/2}\ne 0$.

Given $\xi\in\mathbb A^n$  put $V_\xi:=\bigcap_{\,i=1}^{\,l}\{v\in
T_\xi(\mathbb A^n)\,|\,\,j_{d_\xi Q_i}(v)=0\}$. Suppose $\xi\not\in
{\mathcal J}(Q_1,\ldots, Q_l)$. Then $d_\xi Q_1\wedge\ldots\wedge
d_\xi Q_l\ne 0$ and $\dim V_\xi=n-l$. Since the exterior algebra
$\bigwedge T_\xi({\mathbb A}^n)$ is a free module over its
subalgebra $\bigwedge V_\xi$, it is straightforward to see that
$\bigcap_{\,i=1}^{\,l}\mathrm{Ker}\,j_{d_\xi Q_i}=\bigwedge V_\xi$.
As $\dim\bigwedge^{n-l}V_\xi=1$, our earlier remarks now imply that
$\pi^{(n-l)/2}$ and $w$ are linearly dependent as elements of the
vector space ${\mathbb K}(x_1,\ldots, x_n)\otimes_{\mathcal
A}\big(\bigwedge_{\mathcal A}W\big)$.

Since $d_\xi Q_1\wedge\ldots\wedge d_\xi Q_l\ne 0$, the above
argument also shows that $w\ne 0$. It follows that there exist
nonzero coprime $f_1,f_2\in{\mathcal A}$ such that
$f_1\pi^{(n-l)/2}=f_2 w.$ As the set $\{Q_1,\ldots, Q_l\}$ is
admissible, the function $f_1$ must be constant. As
$\mathrm{Sing}\,\pi$ has codimension $\ge 2$ in $\mathbb A^n$, the
function $f_2$ must be constant as well. Therefore,
$\pi^{(n-l)/2}=\lambda w$ for some nonzero $\lambda\in\mathbb K$,
as stated.
\end{proof}

\smallskip

\noindent 1.4. Next we are going to apply Theorem~\ref{skryabin} to
determine the Poisson centre of certain quasi-regular polynomial
Poisson algebras.
\begin{cl}\label{poisson} Let ${\mathcal
A}=\mathbb K[x_1,\ldots,x_n]$ be a quasi-regular Poisson algebra of
index $l$ and suppose that $\mathcal A=\bigoplus_{k\ge 0}{\mathcal
A}(k)$ is graded in such a way that $x_i\in {\mathcal A}(r_i)$ for
some $r_i>0$, where $1\le i\le n$. Suppose further that $Z({\mathcal
A})$  contains an admissible set $\{Q_1,\ldots,Q_l\}$ consisting of
homogeneous elements of $\mathcal A$. Then $Z(\mathcal A)=\,{\mathbb
K}[Q_1,\ldots, Q_l]$.
\end{cl}
\begin{proof}
By our assumption, $R:={\mathbb K}[Q_1,\ldots,Q_l]$ is a graded
subalgebra of $\mathcal A$ contained in $Z(\mathcal A)$. Let $z$ be
an arbitrary element of $Z(\mathcal A)$. We need to show that $z\in
R$. Our discussion in (1.3) shows that
$$l= \mathrm{tr.\, deg}_{\mathbb K}\, {\mathbb K}(Q_1,\ldots, Q_l)\le
\mathrm{tr.\, deg}_{\mathbb K}\, Z(\mathcal A)\le\ind\mathcal A=l,$$
implying that $z$ is algebraic over $R$. Since ${\mathcal
J}(Q_1,\ldots, Q_l)$ has codimension $\ge 2$ in $\mathbb A^n$, we
can apply Theorem~\ref{skryabin} to complete the proof.
\end{proof}

\smallskip

\noindent 1.5. Let $A=\bigoplus_{k\ge 0} A_k$ be a graded integral
domain over a field $F$. Given $a\in A$ we denote by $\tilde a$ the
initial (lowest) component of $a$. Given an $F$-subalgebra $R$ of
$A$ we let $\tilde R$ denote the $F$-span of all $\tilde{r}$ with
$r\in R$. Clearly, $\tilde{R}$ is a graded $F$-subalgebra of $A$.

\begin{prop}\label{trdeg} Let $A=\bigoplus_{k\ge 0} A_k$ be an
affine graded integral domain over a field $F$ and suppose that
$A_0=F$. Then for any $F$-subalgebra $R$ of $A$ we have ${\rm
tr.\,deg}_F\,\tilde R={\rm tr.\,deg}_F\,R$.
\end{prop}

\begin{proof}
Since the fields of fractions of $R$ and $\tilde{R}$ are isomorphic
to subfields of the field of fractions of $A$, both ${\rm
tr.\,deg}_F\, R$ and ${\rm tr.\,deg}_F\,\tilde R$ are finite. It
follows from \cite[Ch.~II, \S~12, Corollary~2]{ZS} that the field of
fractions of $\tilde{R}$ contains a transcendence basis consisting
of homogeneous elements of $\tilde{R}$. From this it is immediate
that ${\rm tr.\,deg}_F\,\tilde R\le{\rm tr.\,deg}_F\,R$.

Put $m:={\rm tr.\,deg}_F\,\tilde R$ and assume for a contradiction
that $m<{\rm tr.\,deg}_F\,R$. As every algebraically independent
subset of $R$ is contained in a transcendence basis of $R$, our
earlier remarks then show that there exist algebraically
independent elements $a_1,\ldots,a_{m+1}\in R$ such that
${\rm tr.\,deg}_F\,F(\tilde{a}_1,\ldots,
   \tilde{a}_{m+1})=m$. Let $J\subset F[X_1,\ldots,X_{m+1}]$
be the ideal of all polynomial relations between $\tilde{a}_1,
\ldots, \tilde{a}_{m+1}$. Since $F[\tilde{a}_1,\ldots,
\tilde{a}_{m+1}]\subset A$ is a domain of Krull dimension $m$, one
observes easily that $J$ is a prime ideal of codimension $1$ in
the polynomial algebra $F[X_1,\ldots, X_{m+1}]$. As a consequence,
$J$ is generated by one polynomial of positive degree, say $H$.

Let $\tilde{R}_0\subseteq\tilde R$ denote the subalgebra of
initial components of $R_0:=F[a_1,\ldots, a_{m+1}]$. We claim that
$\tilde{R}_0$ is generated by the $\tilde{a}_i$'s and the initial
component $\tilde h$ of $H(a_1,\ldots,a_{m+1})$. To prove the
claim we let $f(a_1,\ldots,a_{m+1})$ be an arbitrary element of
$R_0$. If $\tilde f:=f(\tilde{a}_1,\ldots,\tilde{a}_{m+1})$ is not
zero, then $\tilde{f}$ is the initial component of
$f(a_1,\ldots,a_{m+1})$. If $\tilde f=0$, then $f\in I$ implying
that $f=f_0H$ for some polynomial $f_0$ of smaller degree. Since
$A$ is a domain, the initial component of $f(a_1,\ldots,a_{m+1})$
is nothing but $\tilde{f}_0\tilde{h}$, where $\tilde{f}_0$ is the
initial component of $f_0(a_1,\ldots,a_{m+1})$. Since $\deg
f_0<\deg f$, our claim follows by induction on the degree of $f\in
F[X_1,\ldots, X_{m+1}]$. As a result, the algebra $\tilde{R}_0$ is
finitely generated over $F$.

Next we note that the grading of $A$ induces a descending
filtration ${\mathcal F}=(I_k)_{k\ge 0}$ of $R_0$, where
$I_k=R_0\cap\bigoplus_{i\ge k}A_i$ for all $k$. Furthermore,
$\tilde{R}_0\cong\text{gr}_{\mathcal F}\,R_0$, the corresponding
graded algebra. Consequently, the algebra $\text{gr}_{\mathcal
F}\,R_0$ is Noetherian. Since $A_0=F$, we now apply
\cite[Theorem~4.4.6(b)]{BH} to deduce that $R_0\cong F[X_1,\ldots,
X_{m+1}]$ and $\text{gr}_{\mathcal F}\, R_0\cong\tilde{R}_0$ have
the same Krull dimension. However, $\dim\,R_0=m+1$ whilst
$\dim\,\tilde{R}_0={\rm tr.\,deg}_{F}\,\tilde{R}_0=m$. By
contradiction, the result follows.
\end{proof}

\section{Slodowy slices and symmetric invariants of centralisers}

\noindent 2.1. Let $\chi=(e,\,\cdot\,)$ and $r=\dim \g_e$. The
action of $\ad h$ gives $\g$ a graded Lie algebra structure,
$\g=\bigoplus_{i\in\mathbb Z}\, \g(i)$, where $\gt g(i)=\{x\in
\g\,|\,\,[h,x]=ix\}$. It is well-known that $\gt g_e$ is a graded
Lie subalgebra of the parabolic subalgebra $\gt p:=\bigoplus_{i\ge
0}\,\g(i)$ of $\g$, that is $\gt g_e=\bigoplus_{i\ge 0}\,\g_e(i)$
where $\g_e(i)=\gt g_e\cap\g(i)$. Choose a $\mathbb K$-basis
$x_1,\ldots, x_m$ of $\gt p$ with $x_i\in\g(n_i)$ for some
$n_i\in\mathbb Z_+$, such that $x_1,\ldots, x_r$ is a basis of
$\g_e$ and $x_i\in [f,\gt g]$ for all $i\ge r+1$. Such a basis
exists because $\g=\gt g_e\oplus [\g,f]$ and $\gt p$ contains
$\g_e$.

Define a skew-symmetric bilinear form
$\langle\,\cdot\,,\,\cdot\rangle$ on the subspace $\g(-1)$ by
setting $\langle x , y\rangle=(e,[x,y])$ for all $x,y\in\gt
g(-1)$. As $\g_e\subset\gt p$, this form is nondegenerate. Choose
a basis $z_1,\ldots,z_s,z_{s+1},\ldots,z_{2s}$ of $\gt g(-1)$ such
that
$$\langle z_{i+s},z_j\rangle=\,\delta_{ij},
\qquad\langle z_{i},z_j\rangle=\langle z_{i+s},z_{j+s}
\rangle=0\qquad\ \,(1\le i,j\le r)$$ and denote by $\g(-1)^0$
the linear span of $z_{s+1},\ldots, z_{2s}$. Let $\gt m=\gt
g(-1)^0\oplus\sum_{i\le-2}\,\g(i)$, a nilpotent Lie subalgebra
of dimension $(\dim G\cdot e)/2$ in $\g$.

Given a Lie algebra $\gt s$ over $\mathbb K$ denote by $U(\gt s)$
the universal enveloping algebra of $\gt s$. As $\chi$ vanishes on
the derived subalgebra of $\gt m$, the ideal $N_\chi$ of $U(\gt
m)$ generated by all $x-\chi(x)$ with $x\in\gt m$ has codimension
$1$ in $U(\gt m)$. Let ${\mathbb K}_\chi=U(\gt m)/N_\chi$, a
one-dimensional $U(\gt m)$-module, and denote by $1_\chi$ the
image of $1$ in ${\mathbb K}_\chi$. Set
$$Q_\chi=\,U(\gt
g)\otimes_{U(\gt m)\,}{\mathbb K}_\chi\quad \mbox{and}\quad
H_\chi=\,\mathrm {End}_{\g}(Q_\chi)^{\text{op}}.$$ According to
\cite{Sasha1} and \cite{gg} the associative algebra $H_\chi$ is a
noncommutative filtered deformation of the coordinate algebra
${\mathbb K}[\cS_e]$ endowed with its Slodowy grading
\cite[7.4]{Slo}.

\smallskip

\noindent 2.2. Given $({\bf a},{\bf b})\in \mathbb Z_+^m\times
\mathbb Z_+^s$ we set $x^{\bf a}z^{\bf b}=x_1^{a_1}\cdots
x_m^{a_m}z_1^{b_1}\cdots z_s^{m_s},$ an element of $U(\g)$. By the
PBW theorem, the monomials $x^{\bf a}z^{\bf b}\otimes 1_\chi$,
where $({\bf a},{\bf b})\in \mathbb Z_+^m\times \mathbb Z_+^s$,
form a $\mathbb K$-basis of the induced $U(\g)$-module $Q_\chi$.
For $k\in\mathbb Z_+$ we denote by $Q_\chi^k$ the $\mathbb K$-span
of all $x^{\bf a}z^{\bf b}\otimes 1_\chi$ with
$$|({\bf a},{\bf b})|_e:=\sum_{i=1}^m a_i(n_i+2)+\sum_{i=1}^s
b_i\le k.$$ Any element $h\in H_\chi$ is uniquely determined by
its effect on the canonical generator $1_\chi$. We let $H_\chi^k$
denote the subspace of $H_\chi$ spanned by all $h\in H_\chi$ with
$h(1_\chi)\in Q_\chi^k$. Then $H_\chi=\bigcup_{k\ge 0}\,H_\chi^k$
and $H_\chi^i\cdot H_\chi^j\subseteq H_\chi^{i+j}$ for all
$i,j\in\mathbb Z_+$; see \cite{Sasha1} or \cite{gg}. The
increasing filtration $\{H_\chi^i\,|\,\,i\in \mathbb Z_+\}$ of the
associative algebra $H_\chi$ is often referred to as the {\it
Kazhdan filtration} of $H_\chi$. The corresponding graded algebra
$\gr H_\chi$ is commutative. The elements $x$ from
$Q_\chi^k\setminus Q_\chi^{k-1}$ and $H_\chi^k\setminus
H_\chi^{k-1}$ are said to have {\it Kazhdan degree} $k$, written
$\deg_e(x)=k$.

According to \cite[Theorem~4.6]{Sasha1} the algebra $H_\chi$ has a
distinguished generating set $\Theta_1,\ldots,\Theta_r$ such that
$$\Theta_k(1_\chi)\,=\,\Big(x_k+\sum_{1\le |({\bf i},{\bf j})|_e\le
n_k+2}\lambda_{{\bf i},{\bf j}}^k\,\, x^{\bf i}z^{\bf j}\Big)\otimes
1_\chi,\qquad 1\le k\le r,$$ where $\lambda_{{\bf i},{\bf
j}}^k\in\mathbb K$ and $\lambda_{{\bf i},{\bf j}}^k=0$ if either
$|({\bf i},{\bf j})|_e=n_k+2$ and $|{\bf i}|+|{\bf j}|=1$ or ${\bf
j}= {\bf 0}$ and $i_t=0$ for $t\ge r+1$. The monomials
$\Theta_1^{k_1}\cdots \Theta_r^{k_r}$ and $(\gr\Theta_1)^{k_1}\cdots
(\gr \Theta_r)^{k_r}$ with $(k_1,\ldots,k_r)\in\mathbb Z_+^r$ form
$\mathbb K$-bases of $H_\chi$ and $\gr H_\chi$, respectively.
Furthermore,
$[\Theta_i,\Theta_j]=\Theta_j\circ\Theta_i-\Theta_i\circ\Theta_j\in
H_\chi^{n_i+n_j+2}$ for all $1\le i,j\le r$ (recall that the product
in $H_\chi$ is {\it opposite} to the composition product).

As explained in \cite[Sect.~2]{Sasha2}, there exists a linear map
$\Theta:\g_e\rightarrow H_\chi,\,$ $x\mapsto\Theta_x$
such that $\Theta_{x_i}=\Theta_i$ for all $i$ and
\begin{eqnarray}
[\Theta_{x_i},\Theta_{x_j}]\,\equiv\,\Theta_{[x_i,x_j]}+q_{ij}(\Theta_1,\ldots,
\Theta_r)\quad\ \big(\mathrm{mod}\ H_\chi^{n_i+n_j}\big)\qquad\,
(1\le i,j\le r),
\end{eqnarray}
where $q_{ij}$ is a polynomial in $r$ variables such that $\deg_e
\big(q_{ij}(\Theta_1,\ldots,\Theta_r)\big)=n_i+n_j+2$ and $\deg
\mathrm{in}(q_{ij})\ge 2$ whenever $q_{ij}\ne 0$. Moreover, the map
$\Theta$ has the property that $\Theta_{[x,y]}=[\Theta_x,\Theta_y]$
for all $x\in \g_e(0)$ and $y\in \g_e$. In particular,
$\Theta(\g_e(0))$ is a Lie subalgebra of $H_\chi$ with respect to
the commutator product.

\smallskip

\noindent 2.3. Let $m_1,\ldots, m_l$ be the exponents of the Weyl
group of $\g$. By the Chevalley Restriction Theorem, there exist
algebraically independent elements $F_1,\ldots, F_l\in \cS(\g)^G$
such that $F_i\in S^{m_i+1}(\g)$ for all $i$ and
$\cS(\g)^G={\mathbb K}[F_1,\ldots, F_l]$. Let
$$
\varphi: \g\longrightarrow \,{\mathbb A}^l,\qquad\
x\mapsto \big(\kappa(F_1)(x),\ldots, \kappa(F_l)(x)\big),
$$
be the
adjoint quotient map of $\g$, and let $\varphi_e$ denote its
restriction to the Slodowy slice $\cS_e=e+\g_f$. Composing
$\varphi_e$ with the translation $\tau: \gt
g_f\stackrel{\sim}{\rightarrow}\, \cS_e,\,$ $x\mapsto e+x,$ one
obtains a morphism
$$
\psi:=\varphi_e\circ\tau:\g_f\longrightarrow\, {\mathbb
A}^l,\qquad\ x\mapsto (\psi_1(x),\ldots,\psi_l(x)).
$$
According to
\cite[5.2 \& 7.4]{Slo}, the morphism $\psi$ is faithfully flat with
normal fibres, while in \cite[Sect.~5]{Sasha1} it is proved that all
fibres of $\psi$ are irreducible complete intersections of dimension
$r-l$. It should be mentioned here that each $\psi_i$ is homogeneous
of degree $2m_i+2$ with respect to the Slodowy grading of ${\mathbb
K}[\g_f]\cong {\mathbb K}[\cS_e]$.

Let $U^k$ be the $k$th component of the standard filtration of
$U(\g)$. In view of the PBW theorem, the corresponding graded
algebra $\gr U(\g)$ identifies with the symmetric algebra
$\cS(\g)$. We let $Z(\g)$ denote the centre of of $U(\g)$. It
is well-known that there exist algebraically independent elements
$\tilde{F}_1,\ldots,\tilde{F}_l$ in $Z(\g)$ such that
$\tilde{F}_i\in U^{m_i+1}$ and $\gr\tilde{F}_i=F_i$ for all $i$;
see \cite[7.4]{Di} for example. Moreover, the map taking each
$F_i$ to $\tilde{F}_i$ extends uniquely to an algebra isomorphism
between $\cS(\g)^G$ and $Z(\g)$. Given $F\in \cS(\g)^G$ we shall
denote by $\tilde{F}$ the image of $F$ under this isomorphism.
Note that when $F\in S^k(\g)^G\setminus\{0\}$, we have
$\tilde{F}\in U^k\setminus U^{k-1}$.

Each $\tilde{F}\in Z(\g)$ maps into the centre of $H_\chi$ via
$\tilde{F}\mapsto \tilde{F}(1_\chi)$.
 By \cite[6.2]{Sasha1}, this map is injective.
 To keep the notation simple we shall
 identify the elements of $Z(\g)$ with their images in $Z(H_\chi)$.
 Note that $\tilde{F}_i\in H_\chi^{2m_i+2}\setminus H_\chi^{2m_i+1}$; see
\cite[6.2]{Sasha1}. For $1\le i\le r$, we denote by $\xi_i$ the
restriction of  $\kappa(x_i)$ to $\g_f$, which we regard as a
homogeneous polynomial function of degree $n_i+2$ on $\g_f$. We
denote by $\tilde{\psi}_i$ the image of $\tilde{F}_i$ in the Poisson
algebra $\gr H_\chi$. Clearly, each $\tilde{\psi}_i$ lies in the
Poisson centre of $\gr H_\chi$.

\smallskip

\noindent 2.4. Let $M$ denote the subspace of $\g$ spanned by
$z_1,\ldots, z_s$ and $x_1,\ldots, x_m$. We say that the monomial
$x^{\bf a}z^{\bf b}\in \cS(M)$ has Kazhdan degree $\sum_{i=1}^m
a_i(n_i+2)+\sum_{i=1}^s b_i$. By \cite[6.3]{Sasha1}, the map
$\delta'$ which takes $\gr \Theta_k$ to
$x_k+\sum_{|({\bf i},{\bf j})|_e = n_k+2}\lambda_{{\bf i},{\bf j}}^k\,\, x^{\bf i}z^{\bf j}$
for all $1\le k\le r$ extends to a graded algebra embedding
$\gr H_\chi\hookrightarrow \cS(M)$. Let $\nu:\,\cS(M)\twoheadrightarrow
\cS(\g_e)$ be the graded algebra epimorphism with the property that
$z_i, x_j\in{\rm Ker}\,\nu$ for $1\le i\le s,\,$ $r+1\le j\le m$ and
$\nu(x_k)=x_k$ for $1\le k\le r$. As in \cite[6.3]{Sasha1} we denote
by $\delta''$ the restriction of $\nu\circ\delta'$ to $\gr H_\chi$,
and set $\delta:=\kappa\circ\delta''$.

By \cite[Prop.~6.3]{Sasha1}, the map $\delta\colon\,\gr
H_\chi\rightarrow {\mathbb K}[\gt g_f]$ is a graded algebra
isomorphism satisfying $\delta(\gr \Theta_k)=\xi_k$ for all $k\le
r$ and $\delta(\tilde{\psi}_i)=\psi_i$ for all $i\le l$. This
implies that $\delta''\colon\,\gr
H_\chi\stackrel{\sim}{\rightarrow} \cS(\gt g_e)$ is a graded
algebra isomorphism with the following properties:
\begin{eqnarray}\label{eq2}
\delta''(\tilde{\psi}_i)=\kappa_e^{-1}(\psi_i) \qquad (1\le i\le
l);\qquad\ \delta''(\gr\Theta_i)=\kappa_e^{-1}(\xi_i)=x_i\qquad
(1\le i\le r).
\end{eqnarray}
We use $\delta''$ to transport the Poisson algebra structure of $\gr
H_\chi$ to the symmetric algebra $\cS(\gt g_e)$. Combining (2) and
(\ref{eq2}) we observe that the new Poisson bracket of $\cS(\gt g_e)$
satisfies the following condition:
\begin{eqnarray}\label{eq3}
\{x_i,x_j\}\,=\,[x_i,x_j]+q_{ij}(x_1,\ldots, x_r)\qquad\quad\
(1\le i,j\le r).
\end{eqnarray}
Furthermore, each $\kappa_e^{-1}(\psi_i)$ is in the Poisson centre
of $\cS(\gt g_e)$.

\smallskip

\noindent 2.5. With these preliminaries at hand we are in a position
to prove Proposition~\ref{invariants}.

\smallskip

\noindent{\it Proof of Proposition~\ref{invariants}.}  Let
$F=g(F_1,\ldots,F_l)$ be a homogeneous element of $\cS(\g)^G$ and
let $\tilde{F}=g(\tilde{F}_1,\ldots,\tilde{F}_l)$ be the
corresponding element of $Z(\g)\hookrightarrow H_\chi$; see our
discussion in (2.3). Since each $\tilde{F}_i$ commutes with $h$, the
definition of the Kazhdan filtration and (\ref{eq2}) yield
$$\delta''(\gr \tilde{F})=\,
\delta''\big(g(\tilde{\psi}_1,\ldots,\tilde{\psi}_l)\big)=\,
g\big(\kappa_e^{-1}(\psi_1),\ldots,\kappa_e^{-1}(\psi_l)\big)\,=\,
\kappa_e^{-1}\big(g(\psi_1,\ldots,\psi_l)\big),$$ see
\cite[6.2]{Sasha1} for more detail. Note that $\delta''(\gr
\tilde{F})=\kappa_e^{-1}\big(g(\psi_1,\ldots,\psi_l)\big)$ belongs
to the Poisson centre of $\cS(\g_e)$, that is $\{x,\delta''(\gr
\tilde{F})\}=0$ for all $x\in\g_e$. Abusing notation we denote by
$\ad x$ the derivation of the algebra $\cS(\gt g_e)$ induced by
the inner derivation of $x\in\g_e$. Then
$$
0=\{x,\delta''(\gr \tilde{F})\}\,=\,(\ad
x)\big(\mathrm{in}(\delta''(\gr \tilde{F}))\big)+ \mbox{ terms of
higher standard degree, }
$$ in view of (\ref{eq3}). (One should also keep
in mind that $q_{ij}\ne 0$ implies $\deg\mathrm{in}(q_{ij})\ge 2$).
Since this holds for all $x\in\g_e$, we deduce that
$\mathrm{in}(\delta''\big(\gr \tilde{F}))\in \cS(\g_e)^{\gt g_e}$.
But then
$$^{e\!}F:=\,\kappa_e^{-1}\big(\mathrm{in}(\tau^*(\kappa(F)_{\vert\,\cS_e}))\big)
=\,\kappa_e^{-1}\big(\mathrm{in}(g(\psi_1,\ldots,\psi_l))\big)=\,
\mathrm{in}(\delta''\big(\gr \tilde{F}))\in \cS(\g_e)^{\gt g_e}.$$
We thus obtain $^{e\!} F\in \cS(\g_e)^{\g_e}=\cS(\gt
g_e)^{G_e^\circ}.$

Now let $C_e=G_e\cap G_f$. It is well-known that $C_e$ is a
reductive subgroup of $G_e$, and $G_e$ is generated by $C_e$ and the
unipotent radical $R_u G_e$; see \cite[3.7]{CG} for example.
Clearly, both $\g_f$ and $\cS_e=e+\g_f$ are $C_e$-stable, and the
mappings $\kappa$ and $\kappa_e$ are $C_e$-equivariant. Since $F\in
\cS(\g)^G$, this entails $^{e\!}F\in \cS(\g_e)^{C_e}$. But then
$^{e\!}F\in \cS(\g_e)^{C_e\cdot\, G_e^\circ}=\,\cS(\gt g_e)^{G_e}$,
completing the proof. \qed

\smallskip

\noindent 2.6. Theorem~\ref{or} will enable us to obtain a
differential criterion for regularity of linear functions applicable
to a large class of centralisers in $\mathfrak g$. Recall that a
linear function $\gamma\in{\gt g}_e^*$ is called {\it regular} if
$\dim{\gt g}_e^\gamma=\ind \gt g_e$, where $\gt
g_e^\gamma=\{x\in{\gt g}_e\,|\,\,\gamma([x,{\gt g}_e])=0\}$ is the
stabiliser of $\gamma$ in ${\gt g}_e$.

\begin{thm}\label{summa}
Suppose $\ind \gt g_e=l$. Then the following are true for any
homogeneous generating system $F_1,\ldots, F_l$ of the invariant
algebra $\cS(\gt g)^G$:

\smallskip

\begin{itemize}
\item[(i)\ ] $\sum_{i=1}^l\deg\, ^{e\!}F_i\le (r+l)/2$ where
$r=\dim \gt g_e$.

\medskip

\item[(ii)\ ] The elements $^{e\!}F_1,\ldots, {^{e\!}F}_l$ are
algebraically independent if and only if $\sum_{i=1}^l\deg\,
{^{e\!}}F_i = (r+l)/2$.

\medskip

\item[(iii)\ ] Suppose $\sum_{i=1}^l\deg\, ^{e\!}F_i = (r+l)/2$.
Then the differentials $d_\gamma ({^{e\!}}F_1),\ldots,
d_\gamma({^{e\!}}F_l)$ are linearly independent at $\gamma\in\gt g_e^*$
if and only if $\gamma$ is  regular in $\gt g_e^*$.
\end{itemize}
\end{thm}
\begin{proof}
We are going to apply Theorem\,\ref{or} to the Poisson algebra
$\mathrm{gr}\,H_\chi$. Let $\pi_e$ denote the Poisson bivector of
$\mathrm{gr}\,H_\chi$ and let $\pi_e^{PL}$ be the Poisson bivector
of the polynomial algebra ${\mathcal A}:=\cS(\gt g_e)$ regarded with
its standard Poisson structure. We identify $\mathrm{gr}\,H_\chi$
with $\mathcal A$ by using the recipe described in (2.4) and set
$f_i:=\kappa^{-1}_e(\tau^*(\kappa(F_i)_{|\,\mathcal S_e})),$ $1\le
i\le l$. It follows from \cite[Theorem~5.4]{Sasha1} that the ideal
$(f_1,\ldots, f_l)\subset \mathcal A$ is radical and its zero locus
in $\gt g_e^*$ is normal. This implies that ${\mathcal
J}(f_1,\ldots, f_l)$ has codimension $\ge 2$ in ${\gt g}_e^*$.

From the alternative description of the Poisson structure on
$\mathrm{gr}\,H_\chi$ given in \cite[Sect.~3]{gg} it follows that
$$
\rk\pi_e(\gamma)=\,\dim
(\Ad\,G)\big(e+(\kappa_e^*)^{-1}(\gamma)\big)-\dim (\Ad\,G)\,e \
\qquad \big(\forall\,\gamma\in \g_e^*\big).
$$
Consequently, $\gamma\in\mathrm{Sing}\,\pi_e$ if and only if the
adjoint orbit $(\Ad\,G)\big(e+(\kappa_e^*)^{-1}(\gamma)\big)$ is not
of maximal dimension. By Kostant's criterion for regularity, this
happens if and only if $e+(\kappa_e^*)^{-1}(\gamma)\in{\mathcal
J}\big(\kappa(F_1),\ldots,\kappa(F_l)\big)$. Chasing through the
definitions it is easy to see that the latter happens if and only
of $\gamma\in{\mathcal J}(f_1,\ldots, f_l)$.
 Thus, $\mathrm{Sing}\,\pi_e=\,{\mathcal J}(f_1,\ldots, f_l)$.
Our earlier remarks now show that ${\rm Sing}\,\pi_e$ has
codimension $\ge 2$ in $\g_e^*$. As
$\ind(\mathrm{gr}\,H_\chi)=\ind\gt g$, we conclude that the subset
$\{f_1,\ldots, f_l\}$ is admissible and the Poisson algebra
$(\mathrm{gr}\,H_\chi,\pi_e)$ is quasi-regular.

The standard grading of $\mathcal A$ (by total degree) induces
gradings of the $\mathbb K$-algebras $\Omega$ and
$\bigwedge_{\mathcal A} W$ where we impose that $\deg\,dx_i=0$ and
$\deg\,\partial_i=0$ for all $i$. Our assumption that $\ind \gt
g_e=l$ yields $\big(\pi_e^{PL})^{(r-l)/2}\ne 0$ whereas
(\ref{eq3}) entails that $\mathrm{in}(\pi_e)=\pi_e^{PL}$.
Consequently,
\begin{eqnarray}\label{eq4}
\mathrm{in}(\pi_e^{(r-l)/2})=(\pi_e^{PL})^{(r-l)/2}\ne 0.
\end{eqnarray}
As
$\mathrm{in}(f_i)={^e\!}F_i$ for all $i$, we also have that
\begin{eqnarray}\label{eq5}
\deg\big(\mathrm{in}(df_1\wedge\ldots\wedge
df_r)\big)\ge\deg\,d({^e\!}F_1)\wedge\ldots\wedge d({^e\!}F_l).
\end{eqnarray}
Combining (\ref{eq4}) and (\ref{eq5}) with Theorem~\ref{or} we now
conclude that
$$
\frac{r-l}{2}=\,\deg\big((\pi_e^{PL})^{(r-l)/2}\big)=\,
\deg\big(\mathrm{in}(j_{df_1\wedge\ldots\wedge
df_l}(\partial_1\wedge\ldots\wedge\partial_r)\big)\ge-l+\textstyle{\sum_{i=1}^l}\,\deg
\,{^e\!}F_i.
$$
Statement~(i) follows. Now $^{e\!}F_1,\ldots, {^{e\!}F}_l$ are
algebraically independent in $\cS(\gt g_e)$ if and only if
$d({^e\!}F_1)\wedge\ldots\wedge d({^e\!}F_l)\ne 0$. Since the latter
happens if and only if $\deg\big(\mathrm{in}(df_1\wedge\ldots\wedge
df_r)\big)=\,-l+\textstyle{\sum_{i=1}^l}\deg \,{^e\!}F_i$, the above
argument also yields (ii).

\smallskip

Finally, suppose $\sum_{i=1}^l\deg\, { ^{e\!}}F_i = (r+l)/2$. Then
$\mathrm{in}(df_1\wedge\ldots\wedge
df_r)=d({^e\!}F_1)\wedge\ldots\wedge d({^e\!}F_l)$, and
Theorem~\ref{or} forces
\begin{eqnarray}\label{eq6}
(\pi_e^{PL})^{(r-l)/2}\,=\,\lambda j_{d({^e\!}F_1)\wedge\ldots\wedge
d({^e\!}F_l)}(\partial_1\wedge\ldots\wedge\partial_r),\quad\
\lambda\in\mathbb K^\times.
\end{eqnarray}
Since $\ind\gt g_e=l$, the specialisation of
$(\pi_e^{PL})^{(r-l)/2}$ at $\gamma$ is nonzero if and only if
$\gamma$ a regular linear function of $\gt g_e$. On the other hand,
the RHS of (\ref{eq6}) is nonzero at $\gamma$ if and only if the
differentials $d_\gamma ({^{e\!}}F_1),\ldots, d_\gamma({^{e\!}}F_l)$
are linearly independent. This completes the proof.
\end{proof}

\smallskip

\noindent 2.7. Suppose Elashvili's conjecture holds for $\g_e$.
Simple examples show that the sum of the degrees of
${^e\!}F_1,\ldots, {^e\!}F_l$ depends on the choice of homogeneous
generators $F_1,\ldots,F_l$ of $\cS(\g)^G$. We say that a
homogeneous generating system $\{F_1,\ldots,F_l\}\subset \cS(\g)^G$
is {\it good} for $e$ if
$$\textstyle{\sum}_{i=1}^l\deg{^e\!}F_i = (\dim\g_e+\rk\g)/2.$$
For any generating system $\{F_1,\ldots, F_l\}\subset \cS(\g)^G$
which is good for $e$ the Jacobian locus ${\mathcal
J}({^e\!}F_1,\ldots, {^e\!}F_l)$ is a {\it proper} Zariski closed
subset of $\g_e^*$; see Theorem~\ref{summa}. We say that a
homogeneous generating system $\{F_1,\ldots,F_l\}\subset
\cS(\g)^G$ is {\it very good} for $e$ if the Jacobian locus
${\mathcal J}({^e\!}F_1,\ldots, {^e\!}F_l)$ has codimension $\ge
2$ in $\g_e^*$. It follows from Theorem~\ref{summa}(ii) that for
any very good generating system $\{F_1,\ldots,F_l\}\subset
\cS(\g)^G$ we have the equality
$\textstyle{\sum}_{i=1}^l\deg{^e\!}F_i = (\dim\g_e+\rk\g)/2.$ This
shows that very good systems are good.

We are now in a position to prove the main result of this section:

\begin{thm}\label{verygood}
 Suppose $e$ admits a very good generating system $\{F_1,\ldots, F_l\}\subset \cS(\g)^G$
 and assume further that Elashvili's conjecture holds for
 $\g_e$, that is
 $\ind \gt g_e=l$. Then
$$\cS(\g_e)^{\g_e}=\,\cS(\g_e)^{G_e}=\,{\mathbb
K}[{^e\!}F_1,\ldots, {^e\!}F_l].$$ In particular,
$\cS(\g_e)^{\g_e}$ is a graded polynomial algebra in $l=\rk\g$
variables.
\end{thm}
\begin{proof}
By Theorem~\ref{invariants}, the elements ${^e\!}F_1,\ldots,
{^e\!}F_l$ are in $\cS(\g_e)^{G_e}$. Since $\ind \g_e=l$ and
${\mathcal J}({^e\!}F_1,\ldots, {^e\!}F_l)$ has codimension $\ge 2$
in $\g_e^*$ by our assumption, the set $\{{^e\!}F_1,\ldots,
{^e\!}F_l\}$ is an admissible for the Poisson algebra $\cS(\g_e)$.
Moreover, Theorem~\ref{summa}(iii) shows that the Poisson algebra
$\cS(\g_e)$ is quasi-regular. Applying Corollary~\ref{poisson} to
the Poisson algebra $\cS(\g_e)$ regarded with its standard grading
we now obtain that $S(\g_e)^{\g_e}$ coincides with ${\mathbb
K}[{^e\!}F_1,\ldots, {^e\!}F_l].\,$ Since ${\mathbb
K}[{^e\!}F_1,\ldots, {^e\!}F_l]\subseteq
\cS(\g_e)^{G_e}\subseteq\cS(\g_e)^{\g_e}$, the result follows.
\end{proof}

\begin{rmk}\label{slodowy}
As explained in the proof of Theorem~\ref{summa}, the Poisson
algebra $(\gr H_\chi,\pi_e)$ is quasi-regular and $\{f_1,\ldots,
f_l\}$ is an admissible set for $\gr H_\chi$. Applying
Corollary~\ref{poisson} to the Poisson algebra $\gr H_\chi$
(regarded with its Slodowy grading) we are able to deduce that the
Poisson centre $Z(\gr H_\chi)$ of $\gr H_\chi$ is generated by
$f_1,\ldots, f_l$. In particular, $Z(\gr H_\chi)$ is a polynomial
algebra in $l$ variables. This, in turn, implies that
$Z(H_\chi)=Z(\g)$. We thus recover a result of Victor Ginzburg;
see the footnote in \cite{Sasha2}.
\end{rmk}

\smallskip

\noindent 2.8. Let ${^e\!}Z$ denote the  $\mathbb K$-span of all
${^e\!}F$ with $F\in\cS(\g)^G$, a subalgebra of $\cS(\g_e)^{G_e}$.
For later applications we put on record the following consequence
of Proposition~\ref{trdeg}:
\begin{cl} For any nilpotent element $e\in\g$ we have the equality
$\,\,{\rm tr.deg}_{\mathbb K}({^e\!}Z)=\,\rk\g$.
\end{cl}
\begin{proof}  Recall
that ${^e\!}Z$ coincides with the algebra of initial components of
the subalgebra $\kappa^{-1}_e(\tau^*(\kappa(\cS(\g)^G)_{|\mathcal
S_e}))$ of $\cS(\g_e)$, where the latter is regarded with its
standard grading. Since ${\mathcal S}(\gt g)^G$ is spanned by
homogeneous elements, Proposition~\ref{trdeg} implies that $${\rm
tr.deg}_{\mathbb K}({^e\!}Z)=\,{\rm tr.deg}_{\mathbb
K}\big(\kappa({\mathcal S}(\gt g)^G)_{|\cS_e}\big)=\,{\rm
tr.deg}_{\mathbb K}\, {\mathcal S}(\gt g)^G=\rk\g,$$ as stated (one
should also keep in mind that $\cS_e$ is a slice to the adjoint
orbit of $e$).
\end{proof}

\begin{qn}
Is it true that ${^e\!}Z$ is always finitely generated over $\mathbb
K$?
\end{qn}

\section{Regular linear functions on centralisers}

\noindent 3.1. Given a finite dimensional Lie algebra $\gt q$ and
a linear function $\gamma$ on $\gt q$ we let $\gt q^\gamma$ denote
the stabiliser of $\gamma$ in $\gt q$. Recall that $\ind\gt
q=\,\min_{\gamma\in\gt q^*} \dim\gt q^\gamma$. We set
$$\gt q^*_{\rm sing}:=\{\gamma\in\gt q^*\mid \dim\gt
q^\gamma>\ind\gt q\}.$$ The set $\gt q^*_{\rm reg}:=\gt
q^*\setminus \gt q^*_{\rm sing}$ consists of all regular linear
fuctions of $\gt q$. The main goal of this section is to prove
that $(\g_e^*)_{\rm sing}$ has codimension $\ge 2$ in $\g_e^*$ for
any nilpotent element $e$ in $\g=\gt{gl}_n$ and $\g=\gt{sp}_{2n}$,
where $n\ge 2$. When dealing with $\g=\gt{gl}_n$ we do not impose
any restrictions on the characteristic of $\mathbb K$, whilst for
$\g=\gt{sp}_{2n}$ we require that $\text{char}\,{\mathbb K}\ne 2$.

\smallskip

\noindent 3.2. Let $\VV$ be an $n$-dimensional vector space over
$\mathbb K$ and let $e$ be a nilpotent element in
$\g=\gt{gl}(\VV)$. Let $k$ be the number of Jordan blocks of $e$
and $W\subseteq \VV$ a ($k$-dimensional) complement of $\Im e$ in
$\VV$. Let $d_i+1$ denote the size of the $i$-th Jordan block of
$e$. We always assume that the Jordan blocks are ordered such that
$d_1\ge d_2\ge\ldots\ge d_k$. Choose a basis $w_1, w_2, \ldots,
w_k$ in $W$ such that the vectors $e^{j}{\cdot}w_i$ with $1\le
i\le k$, $0\le j\le d_i$ form a basis for $\VV$, and put
$\VV[i]:=\textrm{span}\{e^j{\cdot}w_i\,|\,\, j\ge 0\}$. Note that
$e^{d_i+1}{\cdot}w_i=0$ for all $i\le k$. When $k=1$, the element
$e$ is regular in $\g$, so that $\g_e$ is abelian of dimension $n$
and $(\g_e^*)_{\rm sing}=\varnothing$. So we assume from now on that
$k\ge 2$.

If $\xi\in\g_e$, then $\xi(e^j{\cdot}w_i)=e^j{\cdot} \xi(w_i)$, hence
$\xi$ is completely determined by its values on $W$. Each vector
$\xi(w_i)$ can be written as
\begin{equation}
\xi(w_i)\,=\,\sum_{j,s} c_i^{j,s}e^s{\cdot}w_j,\ \qquad
c_i^{j,s}\in\mathbb K.
\end{equation}
Thus, $\xi$ is completely determined by the coefficients
$c_i^{j,s}=\,c_i^{j,s}(\xi)$. This shows that $\g_e$ has a basis
$\{\xi_i^{j,s}\}$ such that
$$
\left\{
\begin{array}{l}
\xi_i^{j,s}(w_i)=e^s{\cdot}w_j, \\
\xi_i^{j,s}(w_t)=0 \enskip \mbox{for } t\ne i, \\
\end{array}\right.
\quad 1\le i,j\le k, \ \mbox{ and }\ \max\{d_j-d_i, 0\} \le s\le d_j
\ .
$$
Note that $\xi\in\g_e$ preserves each $\VV[i]$ \iff $c_i^{j,s}(\xi)=0$
for $i\ne j$.

\smallskip

\noindent 3.3 Given a collection $a_1,\ldots, a_k$ of scalars in
$\mathbb K$ we consider the linear function $\alpha$ on $\g_e$
defined by the formula
\begin{equation}          \label{p-ap}
\alpha(\xi)=\sum\limits_{i=1}^k a_i c_i^{i,d_i}\qquad\ \
(\forall\,\xi\in\g_e),
\end{equation}
where $c_i^{j,s}$ are the coefficients of $\xi\in\gt g_e$. Let $\gt
g_e^\alpha$ denote the stabiliser of $\alpha$ in $\gt g_e$. By
aesthetic reasons we prefer it to $(\gt g_e)^\alpha$.

\begin{prop}[\cite{fan}]   \label{ind1}
If the scalars $a_1,\ldots, a_k$ are nonzero and pairwise
distinct, then the stabiliser of $\alpha=\alpha(a_1,\ldots, a_k)$
in $\g_e$ consists of all elements in $\g_e$ preserving the
subspaces $\VV[i]$, where $1\le i\le k$. In other words,
$\g_e^\alpha$ is the linear span of the basis elements
$\xi_i^{i,s}$, and $\dim \g_e^\alpha=n$. In particular,
$\alpha\in(\g_e^*)_{\rm reg}$.
\end{prop}
A direct computation shows that the following commutator relation
holds in $\g_e$:
\begin{eqnarray}\label{commutator}
[\xi_i^{j,s},\xi]=\sum_{t,\ell}c_t^{i,\ell}(\xi)\xi_t^{j,\ell+s}-
\sum_{t,\ell}c_j^{t,\ell}(\xi)\xi_i^{t,\ell+s}\qquad\
(\forall\,\xi\in\g_e);
\end{eqnarray}
see \cite{fan} for more detail. To show that $(\g^*_e)_{\rm sing}$
has codimension $\ge 2$ in $\g_e^*$, for $\g=\gt{gl}(\VV)$, we have
to produce more regular elements in $\g_e^*$.
\begin{prop} \label{beta0}
Define $\beta\in\g_e^*$ by setting
$\beta(\xi)=\sum\limits_{i=1}^{k-1} c_{i+1}^{i,d_i}(\xi)$ for all
$\xi\in\g_e$. Then $\dim \g_e^\beta=n$, so that
$\beta\in(\g_e^*)_{\rm reg}$.
\end{prop}
\begin{proof}
From (\ref{commutator}) and the definition of $\beta$ it follows
that $\beta([\xi_i^{j,s},\xi])=c_{j+1}^{i,d_j-s}(\xi)
-c_j^{i-1,d_{i-1}-s}(\xi)$ for all $\xi\in\g_e$. Suppose $({\rm
ad}^*\xi)\beta=0$. Then $\beta([\xi,\g_e])=0$ forcing
$c_{j+1}^{i,d_j-s}(\xi)=c_j^{i-1,d_{i-1}-s}(\xi)$ for all
$i,j\in\{1,\ldots,k\}$ and all $s$ such that $\max(0,d_j-d_i)\le
s\le d_j$.

We claim that $c_j^{i,s}(\xi)=0$ for $i<j$. Suppose for a
contradiction that this is not the case and take the maximal $j$ for
which there are $i<j$ and $d_i-d_j\le t\le d_i$ such that
$c_j^{i,t}(\xi)\ne 0$. Recall that, according to our convention,
$d_i\ge d_j$. Moreover, $d_{i+1}\ge d_j$, because $i+1\le j$. Set
$s:=d_i-t$. Then $0\le s\le d_j$ and
$c_{j+1}^{i+1,d_j-s}(\xi)=c_j^{i,d_i-s}(\xi)$. As $j+1>j$ and
$i+1<j+1$, the coefficients $c_{j+1}^{i+1,d_j-s}(\xi)$ and
$c_j^{i,t}(\xi)$ are both zero, hence the claim.

Now take $\xi_{i+1}^{i,s}\in\g_e$ with $d_{i}-d_{i+1}\le s\le
d_{i}$. Since $\beta([\xi,\xi_{i+1}^{i,s}])=0$, we have
$c_{i+1}^{i+1,d_{i}-s}(\xi)=c_{i}^{i,d_i-s}(\xi)$. Therefore,
$c_i^{i+1,t}(\xi)=c_{i}^{i,t}(\xi)= c_1^{1,t}(\xi)$ for $0\le t\le
d_{i+1}$. In the same way one can show that
$c_{i}^{i+\ell,t}(\xi)=c_{i-1}^{i+\ell-1,t}(\xi)=
c_{1}^{1+\ell,t}(\xi)$ for $0\le t\le d_{i+\ell}$. It follows that
$\xi$ is completely determined by its effect on $w_1$. So
$\dim\g_e^\beta\le n$ simply because $\xi(w_1)\in\VV$. On the other
hand, $\dim\g_e^\beta\ge \ind \g_e\ge \ind \g=n$ by Vinberg's
inequality. The result follows.
\end{proof}

\smallskip

\noindent 3.4. Let $a\colon\, {\mathbb K}^\times\rightarrow
\GL(\VV)_e$ be the cocharacter such that $a(t)¿w_i=t^iw_i$ for all
$i\le k$ and $t\in{\mathbb K}^{^\times}$, and define a rational
linear action $\rho\colon\,{\mathbb
K}^{^\times}\rightarrow\GL(\g_e^*)$ by the formula
\begin{equation}          \label{ro}
\rho(t)\gamma=\,t(\Ad^*\,a(t))^{-1}\gamma \qquad\ \ \ \ \ \ \ \,
\big(\forall\,\gamma\in\g_e^*, \ \, \forall\,t\in \mathbb
K^{^\times}\big).
\end{equation}
\begin{prop}     \label{cod1}
$(\mathbb K\alpha\oplus\mathbb K\beta)\cap (\g_e^*)_{\rm sing}=0$.
\end{prop}
\begin{proof}
Since $(\Ad\, a(t))(\xi_i^{j,s})=\,t^{j-i}\xi_i^{j,s}$, we have
$(\Ad^*\,a(t))(\alpha)=\alpha$ and $(\Ad^*\,a(t))(\beta)=t\beta$.
Hence $\rho(t)\alpha=t\alpha$ and $\rho(t)\beta=\beta$. So
$\mathbb K\alpha\oplus\mathbb K\beta$ is $\rho(\mathbb
K^{^\times})$-stable and  the induced action of $\rho({\mathbb
K}^{^\times})$ on this plane is a contraction to $\mathbb K\beta$.
Since $\dim(\g_e)^{\rho(t)\gamma}=\,\dim \g_e^\gamma$ and $\beta\in
(\g_e^*)_{\rm reg}$, all linear functions $x\alpha+y\beta$ with
$y\ne 0$ are regular. The linear functions $x\alpha$ with $x\ne 0$
are regular by Proposition~\ref{ind1}.
\end{proof}

\begin{thm} \label{cd2-gl}
Suppose $\dim \VV\ge 2$ and $\g=\mathfrak{gl}(\VV)$. Then for any
nilpotent element $e\in \g$ the locus $(\g_e^*)_{\rm sing}$ has
codimension $\ge 2$ in $\g_e$.
\end{thm}
\begin{proof}
Since $(\g_e^*)_{\rm sing}$ is conical and Zariski closed, the
assertion follows immediately from Proposition~\ref{cod1}.
\end{proof}

\smallskip

\noindent 3.5. Using similar ideas we prove below a symplectic
analogue of Theorem~\ref{cd2-gl}. Our argument in the symplectic
case is more involved. We also provide an example showing that
Theorem~\ref{cd2-gl} does not extend to all nilpotent elements in
orthogonal Lie algebras. We begin with some useful facts on
$\mathbb Z_2$-graded Lie algebras.

Let $\gt q=\gt q_0\oplus\gt q_1$ be a symmetric decomposition
(i.e., a $\mathbb Z_2$-grading) of a Lie algebra
$\gt q$. Then $\gt q^*=\gt q_0^*\oplus\gt q_1^*$. If $\alpha\in\gt q^*$,
then $\tilde{\alpha}$ denotes its
restriction to $\gt q_0$.

\begin{prop}\label{1}
Suppose $\alpha\in\gt q^*$ and $\alpha(\gt q_{1})=0$.
Then $(\gt q_0)^{\tilde{\alpha}}=\gt q^\alpha\cap\gt q_0$.
\end{prop}
\begin{proof}
Take $\xi\in\gt q_0$. Since $[\xi, \gt q_{1}]\subset \gt q_{1}$,
we have that
$\tilde{\alpha}([\xi, \gt q_0])=0$ \iff $\alpha([\xi, \gt q])=0$.
Hence $(\gt q_0)^{\tilde{\alpha}}=(\gt q_0)_\alpha$, where
$(\gt q_0)_\alpha$ is the stabiliser of $\alpha$ in $\gt q_0$.
Clearly $(\gt q_0)_\alpha=\gt q^\alpha\cap\gt q_0$.
\end{proof}
Each $\gamma\in\gt q_0^*$ gives rise to a skew-symmetric bilinear form
$\hat\gamma$ on $\gt q_1$ by
$\hat\gamma(x,y)=\gamma([x,y])$ for all $x,y\in\gt q_1$.
The following assertion is taken from \cite{fan}.

\begin{prop}  \label{prop2}
In the above notation we have  $\ind\gt q\le
\ind\gt q_0+\min\limits_{\gamma\in\gt q_0^*}\dim(\Ker\hat\gamma)$.
\end{prop}
\begin{proof} Take any $\gamma\in\gt q_0^*$ and extend it
to a linear function on $\gt q$ by
setting $\gamma(\gt q_1)=0$.
Then $\gt q^{\gamma}=(\gt q_0)^{\gamma}\oplus(\gt q^\gamma\cap\gt q_{1})=
(\gt q_0)^{\gamma}\oplus\Ker\hat\gamma$. There exists a nonempty
Zariski open subset $U_1\subset\gt q_0^*$ such that
$\dim(\gt q_0)^\gamma=\ind\gt q_0$ for all $\gamma\in U_1\subset\gt q_0^*$.
The linear functions $\gamma$ on $\gt q_0$ for which
$\Ker\hat\gamma$ has the minimal possible dimension form another
nonempty Zariski open subset in $\gt q_0^*$, call it $U_2$. For
each $\gamma\in U_1\cap U_2\neq\varnothing$, the dimension of
$\gt q^\gamma$ equals the required sum, hence the result.
\end{proof}

\begin{lm}\label{ogr}
Suppose $\alpha\in\gt q^*$ is such that $\alpha(\gt q_1)=0$ and
$\dim\gt q^\alpha=\ind\gt q$. Then $\dim(\gt
q_0)^{\tilde{\alpha}}=\ind\gt q_0$.
\end{lm}
\begin{proof}
By Proposition~\ref{prop2} we have:
$$
\ind\gt q_0\ge\ind\gt q -\min\limits_{\gamma\in\gt
q_0^*}\dim(\Ker\hat\gamma)\ge \dim\gt q^\alpha
 -\dim(\Ker\hat\alpha)=\dim(\gt q^\alpha\cap \gt q_0).
$$
Applying Proposition~\ref{1} yields the assertion.
\end{proof}

\smallskip

\noindent 3.6. Let $(\ \,,\ )$ be a nondegenerate symmetric or
skew-symmetric form on $\VV$ and let $J$ be the matrix of $(\ \,,\
)$ with respect to a basis $B$ of $\VV$. Let $X$ denote the matrix
of $x\in\gt{gl}(\VV)$ relative to $B$. The linear mapping $x\mapsto
\sigma(x)$ sending each $x\in\gt{gl}(\VV)$ to the linear
transformation $\sigma(x)$ whose matrix relative to $B$ equals
$-JX^t J^{-1}$ is an involutory automorphism of $\gt{gl}(\VV)$
independent of the choice of $B$. The elements of $\gt {gl}(\VV)$
preserving $(\ \, ,\ )$ are exactly the fixed points of $\sigma$. We
now set $\widetilde{\g}:=\gt{gl}(\VV)$ and let
$\widetilde{\g}=\,\widetilde{\g}_0\oplus\widetilde{\g}_1$ be the
symmetric decomposition of $\widetilde{\g}$ with respect to
$\sigma$.  The elements $x\in\widetilde{\g}_1$ have the property
that $( x{\cdot}v,w) =(v, x{\cdot}w) $ for all $v,w\in \VV$.

Set $\g:=\widetilde{\g}_0$ and let $e$ be a nilpotent element of
$\g$. Since $\sigma(e)=e$, the centraliser $\widetilde{\g}_e$ of
$e$ in $\widetilde{\g}$ is $\sigma$-stable and
$(\widetilde{\g}_e)_0=\,\widetilde{\g}_e^{\sigma}=\,\g_e$. This
yields the $\g_e$-invariant symmetric decomposition
$\widetilde{\g}_e=\,(\widetilde{\g}_e)_0\oplus(\widetilde{\g}_e)_1$.

Suppose that $\dim \VV=2n\ge 4$ and our form is skew-symmetric.
Then $\widetilde{\gt g}_0\cong\gt{sp}_{2n}$. Since $e$ is a nilpotent
transformation of $\VV$, we recycle the notation introduced in
(3.2). Note that in the present case if $d_i$ is even, that is if
the dimension of $\VV[i]=\,\textrm{span}\{e^j{\cdot}w_i\,|\,\,
j\ge 0\}$ is odd, then the restriction of $(\ \, ,\ )$ to $\VV[i]$
is identically zero. By the same reason as in (3.2) it can be
assumed that $k\ge 2$.

\begin{lm}\cite[Sect.~1]{ja}     \label{ab}
The vectors $\{w_i\}_{i=1}^k$ can be chosen such that the
following conditions are satisfied:
\begin{itemize}
\item[(i)\,] if $d_i$ is odd, then the restriction of $(\ \,,\ )$
to $\VV[i]$ is nondegenerate and $(\VV[i],\VV[j])=0$ for any $j\ne
i$.

\item[(ii)\,] if $d_i$ is even, then there is a unique $i'\ne i$
such that $(\VV[i'],\VV[i])\ne 0$.
\end{itemize}
\end{lm}
We thus obtain a decomposition of the set of Jordan blocks of odd
size (i.e., those with $d_i$ even) into pairs $\{i,i'\}$. Note that
$d_{i'}=d_i$ necessarily holds and the restriction of $(\ \,,\ )$ to
$\VV[i]\oplus \VV[i']$ is nondegenerate. For $i\le k$ such that
$d_i$ is odd we put $i'=i$.

Choose vectors $\{w_i\}_{i=1}^k$ according to Lemma~\ref{ab}.
Since the form $(\ \,,\ )$ is $\g$-invariant,
$(e^{d_i}{\cdot}w_i,v)=(-1)^{d_i}(w_i,e^{d_i}{\cdot}v)$ and
$e^{d_i}{\cdot}w_i$ is orthogonal to all $e^{s}{\cdot}w_{j}$ with
either $j\ne i'$ or $s>0$. Since  $(\ \,,\ )$ is nondegenerate, we
also have that $(e^{d_i}{\cdot}w_i, w_{i'})\ne 0$ for all $i$.

\smallskip

\smallskip

\noindent 3.7. Let $\alpha=\alpha(a_1,\ldots,
a_k)\in\widetilde{\g}_e^*$ be as in (\ref{p-ap}) and assume that
$\{a_i\}\subset\mathbb K$ are nonzero and pairwise distinct.
Assume further that $a_{i'}=-a_i$ whenever $i\ne i'$. Then
$\alpha$ vanishes on $(\widetilde{\g}_e)_1$; see \cite[Lemma
2]{fan}. By Lemma~\ref{ogr} and Proposition~\ref{ind1},
$\tilde{\alpha} \in (\g_e^*)_{\rm reg}$. Unfortunately, the linear
function $\beta$ defined in Proposition~\ref{beta0} does not
always vanish on $(\widetilde{\g}_e)_1$. For this reason, we need
a more sophisticated construction.

Renumbering the $\VV[i]$'s if necessary we may assume without loss
of generality that $i'=i\pm 1$ for each pair $\{i,i'\}$. As
$d_i=d_{i'}$, our assumption that $d_1\ge d_2\ge \ldots\ge d_k$
will not be violated. Note that if $i'\ne i+1$, then $i'\le i$ and
$(i+1)'\ge i+1$.

For each $i\le k-1$ with $i'\ne i+1$ we now define a linear
function $\gamma_i$ on $\widetilde{\g}_e$ by setting
$$
\gamma_i(\xi):=-\frac{ (w_{i+1}, e^{d_{i+1}}{\cdot}w_{(i+1)'})   } { (e^{d_i}{\cdot}w_i,w_{i'})  }
c_{i'}^{(i+1)',d_{i+1}}(\xi)   \qquad\ \ (\forall\,\xi\in\widetilde{\g}_e),
$$
and put $\beta':=\sum_{i\le k-1,\,i'\ne\, i+1}\gamma_i$. Recall
from (3.4) that the map $\rho$ gives us a rational action of
$\mathbb K^{^\times}$ on $\widetilde{\g}_e^*$. From Lemma~\ref{ab}
and the definition of $\beta$ it is immediate that $\beta+\beta'\ne 0$.

\begin{lm}\label{beta}
For all $i\le k-1$ with $i'\ne i+1$ we have
$\rho(t)\gamma_i=t^{s_i}\gamma_i$ where $s_i\ge 2$. Moreover,
$\beta+\beta'$ vanishes on $(\widetilde{\g}_e)_1$.
\end{lm}
\begin{proof}
Recall that $(w_i,e^{d_i}{\cdot}w_{i'})\ne 0$ for all $i$. Take
any $\xi\in(\widetilde{\g}_e)_1$. Then
$$
c_{i+1}^{i,d_i}(\xi)(e^{d_i}{\cdot}w_{i},w_{i'})= (\xi(w_{i+1}),w_{i'})=
(w_{i+1},\xi(w_{i'}))=c_{i'}^{(i+1)',d_{i+1}}(\xi) (w_{i+1},e^{d_{i+1}}{\cdot}w_{(i+1)'}).
$$
For $i'\ne i+1$ this yields $c_{i+1}^{i,d_i}(\xi)=-\gamma_i(\xi)$.
Suppose $i'=i+1$. Then also $(i+1)'=i$ and $d_i=d_{i+1}$ is even,
hence
$$
c_{i+1}^{i,d_i}(\xi)(e^{d_i}{\cdot}w_{i},w_{i+1})= c_{i+1}^{i,d_{i}}(\xi)
(w_{i+1},e^{d_i}{\cdot}w_{i})=-c_{i+1}^{i,d_{i}}(\xi) (e^{d_i}{\cdot}w_{i},w_{i+1}).
$$
So $c_{i+1}^{i,d_i}(\xi)=0$ (recall that ${\rm char}\,\mathbb K\ne 2$). But then
$$
(\beta+\beta')(\xi)=\sum_{i\le k-1,\,i'\ne \,i+1}
\big(c_{i+1}^{i,d_i}(\xi)+\gamma_i(\xi)\big)=0.
$$

It follows from \eqref{ro} that $\rho(t)\gamma_i=
t^{s_i}\gamma_i$, where $s_i=(i+1)'-i'+1$. Since $i'\ne i+1$, we
have $i'\le i$ and $(i+1)'\ge i+1$. Then $s_i\ge i+1-i+1=2$.
\end{proof}

Combining Lemma~\ref{beta} with \cite[Lemma 2]{fan}, we observe that
any $\gamma\in\mathbb K\alpha\oplus\mathbb K(\beta+\beta')$ vanishes
on $(\widetilde{\g}_e)_1$. Let $E$ denote the $\mathbb K$-span of
$\tilde{\alpha}$ and $\widetilde{\beta+\beta'}$ in $\g_e^*$.

\begin{prop}\label{regsp}
Under the above assumptions, $\dim E=2$ and $E\cap(\g_e^*)_{\rm
sing}=0$.
\end{prop}
\begin{proof}
Let $\gamma=x\alpha+y(\beta+\beta')$ with $x,y\in\mathbb K$. By
Lemma~\ref{beta}, $\rho(t)\gamma_i=t^{s_i}\gamma_i$, where $s_i\ge
2$, while in (3.4) it is shown that $\rho(t)\alpha=t\alpha$ and
$\rho(t)\beta=\beta$. Since $\alpha$ and $\beta+\beta'$ are nonzero
and $\rho({\mathbb K}^{^\times})$ is diagonalisable, it follows that
$\alpha$ and $\beta+\beta'$ are linearly independent. As both
$\alpha$ and $\beta+\beta'$ vanish on $(\widetilde{\g}_e)_1$, this
yields that $\dim E=2$.

The above discussion also shows that $\lim_{t\to
0}\rho(t)\gamma=y\beta$. If $y\ne 0$, then Proposition~\ref{beta0}
gives $\gamma\in(\widetilde{\g}_e^*)_{\rm reg}$. By
Proposition~\ref{ind1}, $\alpha\in(\widetilde{\g}_e^*)_{\rm reg}$ as
well. Then $\dim(\widetilde{\g}_e)^\gamma=2n=\ind\widetilde{\g}_e$
for any nonzero $\gamma\in \mathbb K\alpha\oplus\mathbb
K(\beta+\beta')$. As any such $\gamma$ vanishes on
$(\widetilde{\g}_e)_1$, applying Lemma~\ref{ogr} we now conclude
that $E\setminus\{0\}\subset(\g_e^*)_{\rm reg}$. Equivalently,
$E\cap(\g_e^*)_{\rm sing}=0$.
\end{proof}

The following is an immediate consequence of Proposition~\ref{regsp}.

\begin{thm} \label{cd2-sp} Let  $e$ be a
nilpotent element of $\g=\gt{sp}_{2n},\,$
$n\ge 2$. Then $(\g_e^*)_{\rm sing}$ has codimension $\ge 2$ in
$\g_e^*$.
\end{thm}
\begin{proof}
Straightforward (see the proof of Proposition~\ref{cd2-gl}).
\end{proof}

\smallskip

\noindent 3.8. We shall see in a moment that there are nilpotent
elements $e$ in the orthogonal Lie algebra $\gt g=\gt{so}(\VV)$ for
which $(\gt g_e^*)_{\mathrm{sing}}$ has codimension $1$ in $\g_e^*$.
But first we would like to give two positive examples.

Suppose $\dim\VV$ is odd and let $e$ be a nilpotent element in
$\gt{so}(\VV)$ with $2m+1$ Jordan blocks indexed by the integers
ranging from $-m$ to $m$, where $m\ge 1$. Similar to the symplectic
case we may assume that there is an involution $i\to i'$ on the set
of indices such that $i'=i$ if and only if $d_i$ is even and
$(\VV[i],\VV[j])=0$ whenever $j\ne i'$. Recall that $d_{i'}=d_i$
necessarily holds.

Suppose that $i'=-i$ and $d_i\le d_j$ for $i>j\ge 0$. Then $d_0$ is
even and the other $d_i$ are odd. Choose ${\mathbb K}[e]$-generators
$w_i\in\VV[i]$ such that $(w_i,e^{d_i}w_{-i})=1$ for $i\ge 0$, and
let $\tilde\alpha$ denote the restriction to $\g_e$ of the linear
function $\alpha$ on $\gt{gl}(V)_e$ given by
$$
\alpha(\xi)=\sum\limits_{i=-m+1}^m c_{i-1}^{i,d_i}(\xi) \qquad \ \
\big(\forall\,\xi\in\gt{gl}(V)_e\big).
$$
By \cite[Section 4]{fan}, this linear function is regular. Let
$g\in{\rm GL}(\VV)$ be such that
$$
g(w_i)=w_{-i} \enskip \text{ for } i\ge 0,\  \ \enskip
g(w_i)=-w_{-i} \text{ for } i<0,\  \enskip  \
g(e^s{\cdot}w_i)=(-1)^se^s{\cdot}w_i \enskip \text{ for } s\ge 1.
$$
Then $g\in{\rm O}(\VV)$ and $(\Ad\,g)e=-e$, i.e., $g$  normalises
${\mathbb K}e$. Hence $\Ad\, g$ acts on $\g_e$ as a Lie algebra
automorphism. Set $\tilde\beta:=(\Ad^*\,g)\tilde\alpha$. In
coordinates,
$$
\tilde\beta(\xi)=\sum\limits_{i=-m}^{m-1}
c_{i+1}^{i,d_i}(\xi)-2c_1^{0,d_0}(\xi) \qquad \ \
\big(\forall\,\xi\in\gt g_e\big).$$ Set $E':=\mathbb
K\tilde\alpha+\mathbb K\tilde\beta$, a subspace of $\g_e^*$.
\begin{lm}\label{codim2-so1}
The subspace $E'$ is $2$-dimensional and $E'\cap(\gt g_e^*)_{\rm
sing}=0$. The singular locus $(\gt g_e^*)_{\rm sing}$ has
codimension $\ge 2$ in $\g_e^*$.
\end{lm}
\begin{proof} By \cite[Section~4]{fan}, the function $\tilde\alpha$ is regular in $\g_e^*$. Hence
so is $\tilde\beta=(\Ad^*g)\tilde\alpha$. In particular, both
$\tilde{\alpha}$ and $\tilde{\beta}$ are nonzero. There exists a
cocharacter $a\colon\,{\mathbb K}^{^\times}\rightarrow {\rm
SO}(\VV)_e$ such that $a(t)w_i=t^i w_i$ for all $i$. It has the
property that $(\Ad^* a(t))\tilde\alpha=t^{-1}\tilde\alpha$ and
$(\Ad^* a(t))\tilde\beta=t\tilde\beta$. This implies that $\dim
E'=2$. Since the Zariski closed set $E'\cap(\gt g_e^*)_{\rm sing}$
is conical and $(\Ad^* a({\mathbb K}^{^\times}))$-stable and both
$\tilde\alpha$ and $\tilde\beta$ are regular, it also follows that
$E'\cap(\gt g_e^*)_{\rm sing}=0$.
\end{proof}

Suppose now that $\tilde\VV=\VV\oplus\mathbb K v$ is an even
dimensional vector space such that $(\VV,v)=0$ and $(v,v)=1$. Let
$e\in\gt{so}(\VV)$ be the same nilpotent element as above (one with
$2m+1$ Jordan blocks and with $i'=-i$ for all $i$). We regard $e$ as
a nilpotent element of $\gt{so}(\tilde\VV)$ by setting $e(v)=0$.
Then $e\in\gt{so}(\tilde\VV)$ has $2m+2$ Jordan blocks. Assume that
the new Jordan block of size $1$ is indexed by $M$ with $M>m$ and
that $v$ is its generator.

\begin{lm}\label{codim2-so2}
For $e$ as above, the singular locus $(\gt g_e^*)_{\mathrm{sing}}$
has codimension $\ge 2$ in $\g_e^*=\gt{so}(\tilde\VV)^*_e$.
\end{lm}
\begin{proof} Note that $\gt{so}(\VV)$ is a symmetric subalgebra
of $\g=\gt{so}(\tilde\VV)$. Let
$\gt{so}(\tilde\VV)=\gt{so}(\VV)\oplus\gt p$ be the corresponding
symmetric decomposition. Then we can identify the dual space of the
centraliser $\gt{so}(\VV)_e$ with the annihilator of $\gt p_e:=\gt
p\cap\gt g_e$ in $\gt g_e^*$. Let $\tilde\alpha$ and $\tilde\beta$
be the same linear functions as in Lemma~\ref{codim2-so1}. We view
them as linear functions of $\gt g_e$ vanishing on $\p_e$. As ${\rm
O}(\VV)\hookrightarrow{\rm O}(\tilde\VV)$ and ${\rm
SO}(\VV)_e\hookrightarrow G_e$, we still have that
$\tilde\beta=(\Ad^* g)\tilde\alpha$ and $(\Ad^*
a(t))\tilde\alpha=t^{-1}\tilde\alpha$, $(\Ad^*
a(t))\tilde\beta=t\tilde\beta$ for the same cocharacter
$a\colon\,{\mathbb K}^{^\times}\rightarrow G_e$ as in
Lemma~\ref{codim2-so1}. Therefore, in order to prove the statement
it suffices to show that $\tilde\alpha\in(\gt g_e^*)_{\rm reg}$. By
construction,
$$
\dim\gt g_e^{\tilde\alpha}=(\dim\tilde\VV)/2-1+
     \dim\,\{\xi\in\gt p_e\,|\,\, \tilde\alpha([\xi,\gt p_e])=0\}.
$$
The linear space $\gt p_e$ has a basis
$\{\xi_i^{M,0}+\epsilon(i)\xi_M^{-i,d_i}\,|\, -m\le i\le m\}$ where
$\epsilon(i)=-1$ for $i\ge 0$ and $\epsilon(i)=1$ for $i<0$. Using
(\ref{commutator}) we get
\begin{eqnarray*}
\tilde\alpha\big([\xi_i^{M,0}+\epsilon(i)\xi_M^{-i,d_i},\,\xi_j^{M,0}+\epsilon(j)\xi_M^{-j,d_j}]\big)&=&0
\enskip \mbox{for }\  j\ne -i-1;\\
\tilde\alpha\big([\xi_i^{M,0}+\epsilon(i)\xi_M^{-i,d_i},\,
\xi_{-i-1}^{M,0}+\epsilon(-i-1)\xi_M^{i+1,d_{i+1}}]\big)&=&
2\epsilon(i),\ \quad -m\le i\le m.
\end{eqnarray*}
It follows that $\tilde\alpha$ induces on $\gt p_e$ a skew-symmetric
bilinear form of rank $2m$. But then
$$
\dim\,\{\xi\in\gt p_e\,|\,\, \tilde\alpha([\xi,\gt p_e])=0\}=1
$$
and the statement follows.
\end{proof}

\noindent 3.9. For any simple Lie algebra $\g$ of type different
from $\mathbf A$ and $\mathbf C$ we provide in this subsection a
uniform construction of $e\in\cN(\g)$ for which
$(\g_e^*)_{\mathrm{sing}}$ has codimension $1$ in $\g_e^*$. We
assume for simplicity that ${\rm char}\,\mathbb K=0$. The Lie
algebras $\sln$ and $\spn$ are distinguished by the property that
their highest root is not a fundamental dominant weight. This
seemingly insignificant fact is a source of many structural
differences between $\sln$ and $\spn$, and the other finite
dimensional simple Lie algebras. In our situation, it manifests
itself as follows.

Let $G{\cdot}\tilde e=\co_{\min}$ be the minimal nilpotent orbit in
$\g$ and $\{\tilde e,\tilde h,\tilde f\}$ an
$\mathfrak{sl}_2$-triple. Consider the $\mathbb Z$-grading
determined by $\tilde h$
$$\g=\bigoplus_{i=-2}^2 \g(i).$$
Here $\g(2)=\mathbb K\tilde e$ and $\g(-2)=\mathbb K\tilde{f}$.
Let $G(0)$ denote the stabiliser of $\tilde{h}$ in $G$. This is a
Levi subgroup of $G$ which acts on $\g(1)$ with finitely many
orbits. If $\g\ne\sln$, then the centre of $\g(0)$ is
one-dimensional and $\g(1)$ is a simple $\g(0)$-module.
Furthermore, if $\g\ne\spn$, then the open $G(0)$-orbit in $\g(1)$
is {\it affine\/}. Let $e\in\g(1)$ be a point in this orbit.

From now on we assume in this subsection that $\g$ is not isomorphic
to $\sln$ or $\spn$. Our goal is to prove that $(\g_e^*)_{\rm sing}$
has codimension $1$ in $\g_e^*$. Set $\el=[\g(0),\g(0)]$ and let $K$
denote the stationary subgroup of $e$ in $G(0)$. Then $\ka:=\Lie K$
is a Lie subalgebra of $\el$ acting trivially on $\g(2)$. The
centraliser $\g_e$ is graded and has the following structure. Its
component of degree 0 is $\ka$ and its component of degree 1 is
isomorphic as a $\ka$-module to $\mathbb K e \oplus W\oplus W^*$,
where $W$ is a $\ka$-module of dimension $\frac{\dim\g(1)}{2}-1$.
The component of degree 2 is still $\mathbb K\tilde e$. Consider the
hyperplane $\cH=\{\gamma\in\g_e^*\,|\,\, \gamma(\tilde{e})=0\}$. We
wish to prove that $\cH\subset (\g_e^*)_{\mathrm{sing}}$. Because
$\tilde e$ acts trivially on $\cH$, the representation of
$\g_e/\mathbb K\tilde e$ in $\cH$ is equivalent to the coadjoint
representation of $\g_e/\mathbb K\tilde e$. That is, we have to
compute the index of this Lie algebra. Modulo the trivial direct
summand $\mathbb Ke$, this algebra is the semi-direct product of
$\ka$ and $W\oplus W^*$, denoted $\ka\ltimes(W\oplus W^*)$. For such
semi-direct products, one can use Ra\"\i s' formula for the index
\cite{rais}. As the generic stabiliser for the representation of
$\ka$ on $W\oplus W^*$, say $\mathfrak s$, is reductive, Ra\"\i s'
formula yields
$$
\ind(\ka\ltimes(W\oplus W^*))=\,\rk \mathfrak s+\dim (W\oplus W^*)\md K .
$$
It turns out that in all cases of interest for us this number
equals $\rk\g$. Taking into account the direct summand
$\mathbb Ke$ and the passage to $\cH$, we see that generic $G_e$-orbits in
$\cH$ are of codimension $\rk\g+2$ in $\g_e^*$. On the other hand,
it is straightforward to see that for any linear function
$\gamma\in\g_e^*\setminus{\mathcal H}$ satisfying
$\gamma_{\,|\g_e(1)}=0$ one has $\dim(\g_e)_\gamma=\dim\ka_\gamma+2$. As
$\ka$ is reductive with $\rk\ka=\rk \g-2$, this implies that
$\ind\g_e=\rk\g$. Then $\cH\subset (\g_e^*)_{\mathrm{sing}}$, as
wanted.

In Tables~\ref{data-exc} and \ref{data-so}, we provide the
necessary information related to these computations. In
Table~\ref{data-exc}, $W$ is always a simple $\ka$-module which is
represented by its highest weight.

\begin{table}[htb]
\caption{\sc Data for the exceptional Lie algebras}
\label{data-exc}
\begin{tabular}{cccccccc}
$\g$ & $\el$ & $\ka$ & $W$ & $\dim W$ & $\mathfrak s$ & $\dim(W\oplus W^*)\md K$ &
$\ind(\ka\ltimes(W\oplus W^*))$ \\ \hline
$\GR{E}{6}$ & $\GR{A}{5}$ & $2\GR{A}{2}$ & $\varpi_1+\varpi'_1$ & 9 & $T_2$ & 4 & 6 \\
$\GR{E}{7}$ & $\GR{D}{6}$ & $\GR{A}{5}$ & $\varpi_2$ & 15 & $(\GR{A}{1})^3$ & 4 & 7 \\
$\GR{E}{8}$ & $\GR{E}{7}$ & $\GR{E}{6}$ & $\varpi_1$ & 27 & $\GR{D}{4}$ & 4 & 8 \\
$\GR{F}{4}$ & $\GR{C}{3}$ & $\GR{A}{2}$ & $2\varpi_1$ & 6 & 0 & 4 & 4 \\
$\GR{G}{2}$ & $\GR{A}{1}$ & 0 & 0 & 1 & 0 & 2 & 2 \\ \hline
\end{tabular}
\end{table}

\begin{table}[htb]
\caption{\sc Data for $\son$, $n\ge 7$}   \label{data-so}
\begin{tabular}{cccccccc}
$\g$ & $\el$ & $\ka$ &  $\dim W$ & $\mathfrak s$ & $\dim(W{\oplus} W^*)\md K$ &
$\ind(\ka{\ltimes}(W{\oplus} W^*))$ \\ \hline
$\mathfrak{so}_7$ & $\mathfrak{so}_{3}{\times}\mathfrak{sl}_2$ &
$\mathfrak t_1$ & 2 & 0 & 3 & $3$ \\
$\begin{array}{c} \son \\ (n{\ge} 8)\end{array}$ &
$\mathfrak{so}_{n-4}{\times}\mathfrak{sl}_2$ &
$\mathfrak{so}_{n-6}{\times}\mathfrak t_1$ & $n{-}5$ &
$\mathfrak{so}_{n-8}$ & 4 & $[n/2]$
\\ \hline
\end{tabular}
\end{table}

\smallskip

\noindent 3.10. Adopt the notations and conventions of (3.9) and let
$\tilde{e}$ be an element in the minimal nilpotent orbit ${\mathcal
O}_{\min}$. Then $\ind\gt g_{\tilde{e}}=\rk\gt g$ by \cite{dima2}.
We now wish to investigate the singular locus of $\g_{\tilde{e}}^*$.
\begin{thm}\label{min-nilp:sing}
If $\rk \g\ge 2$ and $\tilde{e}\in{\mathcal O}_{\min}$, then
$(\g_{\tilde{e}}^*)_{\rm sing}$ has codimension $\ge 2$ in
$\g_{\tilde{e}}^*$.
\end{thm}
\begin{proof} In view of Theorems~\ref{cd2-gl} and~\ref{cd2-sp}
 the statement holds when $\g$ is of type $\mathbf A$ or
 $\mathbf C$. So we may assume in this proof that $\g$ is not
 isomorphic to $\sln$ or $\spn$. Then
 $$\g_{\tilde{e}}=\el\oplus\g(1)\oplus\g(2).$$ Since $\dim\g(2)=1$,
 we have a skew-symmetric bilinear form
 $\langle\,\cdot\,,\cdot\,\rangle$ on $\g(1)$ such that
 $[x,y]=\langle x,y\rangle \tilde{e}$ for all $x,y\in\g(1)$. This
 form is nondegenerate.

 Given a subset $X\subset\gt g_{\tilde{e}}$  we denote by
 $\mathrm{Ann}(X)$ the annihilator of $X$ in $\gt g_{\tilde{e}}^*$,
 that is
 $$
 \mathrm{Ann}(X):=\{\gamma\in\gt g_{\tilde{e}}^*\,|\,\,
 \gamma(X)=0\}.
 $$
 Then $\mathrm{Ann}(\tilde{e}):=\mathrm{Ann}(\{\tilde e\})$ is a
 hyperplane in $\gt g_{\tilde{e}}^*$. We claim that
 $\mathrm{Ann}(\tilde{e})\not\subset(\g_{\tilde{e}}^*)_{\mathrm{sing}}$.
 To prove the claim we are going to argue in the spirit of (3.9).

 Let $L$ denote the derived subgroup of $G(0)$. Since $\tilde e$
 acts trivially on $\mathrm{Ann}(\tilde{e})$, the representation of
 $\g_{\tilde{e}}/\mathbb K\tilde e$ in $\mathrm{Ann}(\tilde{e})$ is
 equivalent to the coadjoint representation of
 $\g_{\tilde{e}}/\mathbb K\tilde e$. This Lie algebra is the
 semi-direct product of $\el$ and $\g(1)$, denoted
 $\el\ltimes\g(1)$.  The generic stabiliser for the representation
 of $\el$ on $\g(1)$ is isomorphic to $\mathfrak k$. Since $\ka$ is
 reductive, Ra\"\i s' formula \cite{rais} yields
 $$
 \ind(\el\ltimes\g(1))=\,\rk \mathfrak k+\dim \g(1)\md L.
 $$
 As the complement $\g(1)\setminus G(0){\cdot}e$ is a hypersurface in $\g(1)$
 and the semisimple group $L$ has codimension $1$ in $G(0)$, the
 orbit $L{\cdot}e$ has codimension $1$ in $\g(1)$. This implies
 that $\dim \g(1)\md L=1$, whereas Tables~\ref{data-exc} and
 ~\ref{data-so} yield $\rk \mathfrak k=\rk\g-2$. Therefore,
 $\ind(\el\ltimes\g(1))=\rk\g-1$. Each
 $\gamma\in\mathrm{Ann}(\tilde e)$  may be regarded as a linear function on
 $\el\ltimes\g(1)$. Moreover, it is easy to see that $\gt g_{\tilde{e}}^\gamma\cong\mathbb
 K\tilde{e}\oplus (\el\ltimes\g(1))^\gamma$ for every
 $\gamma\in\mathrm{Ann}(\tilde e)$. This implies that for a generic
 $\gamma\in\mathrm{Ann}(\tilde e)$ we have $\dim\gt
 g_{\tilde{e}}^\gamma=\rk\g=\ind\g_{\tilde e}$. The claim follows.

 It remains to deal with  the affine open set $Y:=\g_{\tilde
 e}^*\setminus\mathrm{Ann}(\tilde e)$. Set $\gt n:=\g(1)\oplus\g(2)$
 and let $N\subset G_{\tilde e}$ be the connected subgroup of $G$ with
 $\Lie N=\gt n$. The derived subgroup $(N,N)$ is $1$-dimensional
 with $\Lie (N,N)=\mathbb K\tilde{e}$, and $N/(N,N)\cong \g(1)$ as varieties.
 Let $\alpha\in\mathrm{Ann}(\el\oplus\g(1))$ be a
 nonzero function. The set $\mathrm{Ann}(\g(1))\cap Y$ is Zariski
 closed in $Y$ and can be identified with $\el^*\oplus\mathbb
 K^{^\times}\!\!\alpha$. Let $\gamma=\beta+a\alpha$ be an element
 of  $\mathrm{Ann}(\g(1))\cap Y$ with $\beta\in\el^*$ and $a\ne 0$.
 Then
 $$
 (\mathrm{Ad}^*N)\gamma\,=\,\Big\{\beta+\frac{a}{2}(\mathrm{ad}^*
 v)^2\alpha+a(\mathrm{ad}^*v)\alpha+a\alpha\,|\,\, v\in\g(1)\Big\}.
 $$
 Since the form $\langle\,\cdot\,,\cdot\,\rangle$ is nondegenerate,
 it follows that the $N$-saturation of $Y\cap\mathrm{Ann}(\g(1))$
 is equal to $Y$, that each $N$-orbit $(\mathrm{Ad}^*N)\gamma$ is
 isomorphic to $N/\mathbb (N,N)\cong \g(1)$, and that
 $\g_{\tilde{e}}^\gamma=\el_\beta\oplus\mathbb K\tilde{e}$. In
 particular, the action morphism
 $$
 \tau\colon\, \big(N/(N,N)\big)\times\big(\mathrm{Ann}(\g(1))\cap Y\big)\to\, Y
 $$
 is an isomorphism. Suppose $g\in N/(N,N)$ and
 $\gamma=\beta+a\alpha$, where $\beta\in\el^*$ and $a\ne 0$. Then
 $\tau((g,\gamma))\in(\gt g_{\tilde e}^*)_{\mathrm{reg}}$ if and
 only if $\beta\in(\el^*)_{\mathrm{reg}}$. Since
 $(\el^*)_{\mathrm{sing}}$  has codimension $3$ in $\el^*$, the
 intersection $(\g_{\tilde{e}}^*)_{\mathrm{sing}}\cap Y$ is of
 codimension $3$ in $Y$ and also in $\gt g_{\tilde e}^*$. The
 result follows.
 \end{proof}

\section{Degrees of basic invariants}

\noindent
4.1. From now on  we assume that ${\rm char}\,\mathbb
K=0$. Let $Q$ be a connected linear algebraic group with Lie
algebra $\gt q$. Suppose we are given a rational linear action of
$Q$ on a vector space $V$. The differential of this action at the
identity element of $Q$ is a representation of the Lie algebra
$\gt q$ in $V$.

\begin{df}\label{generic} A vector $x\in V$ (a
stabiliser $\gt q_x$) is called {\it a generic point} ({\it a
generic stabiliser}), if there exists a Zariski open subset
$U\subset V$ such that $x\in U$ and $\gt q_x$ is $Q$-conjugate to
any $\gt q_y$ with $y\in U$.
\end{df}

Let $e$ be a nilpotent element in $\g=\gt{gl}(\VV)$ and set
$G:=\GL(\VV)$. Let $\alpha=\alpha(a_1,\ldots, a_k)\in\g_e^*$ be as
in (3.3) and put $\gt h:=\g_e^\alpha$.

\begin{prop}[\cite{fan}] \label{stab} If all $a_1,\ldots, a_k$ are nonzero
and pairwise distinct, then
$\gt h$ is a generic stabiliser for the coadjoint representation
of $G_e$.
\end{prop}

For $1\le i\le n$, let $\Delta_i$ denote the sum of the principal
minors of order $i$ of the generic matrix
$(x_{ij})_{1\le i,j\le n}$, a regular function on $\g$, and set
$F_i:=\kappa^{-1}(\Delta_i)$. It is well-known that
$\{F_1,\ldots,F_n\}$ is a generating set of the invariant algebra
${\mathcal S}(\gt g)^G$. Recall from (0.2) the definition of
${^e\!}F_1,\ldots, {^e\!}F_n$. Let $(d_1+1\ge d_2+1\ge\cdots\ge
d_k+1)$ be the partition of $n$ corresponding to $e$ and put
$d_0=0$.

\begin{thm}\label{deg1} Suppose $\dim\VV\ge 2$ and let
$F_1,\ldots F_n$ be as above. Then $\{F_1,\ldots, F_n\}$ is a very
good generating system for $e$ and
$\cS(\g_e)^{\g_e}=\,\cS(\g_e)^{G_e}=\,{\mathbb K}[{^e\!}F_1,\ldots, {^e\!}F_n].$ Moreover,
$$
\deg({^e\!}F_{d_0+\cdots+d_{i}+i+1})=\cdots=
\deg({^e\!}F_{d_0+\cdots+d_i+d_{i+1}+i+1})=i+1\qquad\
\ \, (0\le i\le k-1).
$$
\end{thm}
\begin{proof}
Let $\alpha$ be as in Proposition~\ref{stab} and and let $\gt r$
be the linear span of all $\xi_i^{j,s}$ with $i\ne j$. Let $\gt t$
be the span of all $\xi_{i}^{i,0}$, a maximal toral subalgebra of
$\g_e$. Then the centraliser
$\gt h=\gt c_{\g_e}(\gt t)$ is an abelian Cartan
subalgebra of $\g_e$. Moreover, $\g_e=\gt h\oplus\gt r$ and $[\gt h,\gt r]= \gt r$
(this follows from the formula displayed in the
proof of Proposition~2 in \cite{fan}).
We identify $\gt h^*$ with ${\rm Ann}(\gt r)\subset\g_e^*$.
The above implies that $\gt h^*=\{\gamma\in\g_e^*\,|\,\,({\rm ad}^*\,\gt h)\gamma=0\}$. Since
$\gt h$ is a generic stabiliser, we have
$\overline{G_e\cdot \gt h^*}=\g_e^*$.
Therefore, the restriction map $\varphi\mapsto\varphi_{|\gt h^*}$
induces an embedding $\mathbb K[\g_e^*]^{G_e}\hookrightarrow \mathbb K[\gt h^*]$.
It follows
that each ${^e\!}{F_i}_{|{\gt h^*}}$ is nonzero and hence has the
same degree as ${^e\!}F_i$.

 Let $\gt r^{\perp}\subset\gt g$ be the orthogonal
complement to $\gt r$ with respect to $\kappa$ and $\gt
s:=\cS_e\cap\gt r^{\perp}$. Then $\gt s=e+(\kappa_e^*)^{-1}({\rm
Ann}\,\gt r)$, implying that  ${{^e\!}F_i}_{|\gt
h^*}={{^e\!}F_i}_{|{{\rm Ann}\,\gt r}}$ is equal to the component
of minimal degree of the restriction of $\Delta_i$ to $\gt s$. Let
$\gt g[i]\cong\gt{gl}(\VV[i])$ be the subalgebra of $\gt g$
consisting of all endomorphisms acting trivially on $\VV[j]$ for
$j\ne i$, and $\hat{\gt g}:=\g[1]\oplus\cdots\oplus\gt g[k]$. Then
$\hat{\g}$ is a Levi subalgebra of $\g$ and $\gt
s\subset\hat{\g}$.

For $1\le \ell\le d_j+1$ we denote by $\Delta_{\ell}[j]$ the sum
of all principal minors of order $\ell$ of the generic matrix
$\big(x_{pq}^{(j)}\big)_{1\le p,q\le d_j+1}$, a homogeneous
element of degree $\ell$ in $\mathbb K[\hat{\g}]$, and put
$\Delta_0[j]=1$. Since the characteristic polynomial of a
block-diagonal matrix is the product of the characteristic
polynomials of its blocks, it follows that
$$
{\Delta_\ell}_{\vert\hat{\g}}\,=\sum_{\ell_1+\cdots\ell_k=\ell}\Delta_{\ell_1}[1]
\cdots\Delta_{\ell_k}[k]\qquad\ \qquad \, (1\le \ell\le n).
$$
As $e_{\vert\VV[i]}$ is a regular nilpotent element of
$\gt{gl}(\VV[i])$, for each $\ell_i\ge 1$ we have the inequality
$\deg {^e\!}\big(\kappa^{-1}(\Delta_{\ell_i}[i])\big)\ge 1$. It
follows that $\deg {^e\!}F_\ell\ge q$, where
$$q:=\min\,\{s\,|\,\,\,\ell=t_1+\cdots+t_s,\,\,\ 0<t_i\le d_i+1\}.$$ More
precisely,
\begin{eqnarray*}
\deg{^e\!}F_i&=&1 \quad\mbox{for }\
i=1,\ldots, d_1+1, \\
\deg{^e\!}F_i&\ge& 2 \quad\mbox{for }\
i=d_1+2,\ldots, d_1+d_2+2,\\
\deg{^e\!}F_i&\ge& 3 \quad\mbox{for }\
i=d_1+d_2+3,\ldots,d_1+d_2+d_3+3,
\end{eqnarray*}
and so on. Consequently, $\sum_{i=1}^n\deg{^e\!}F_i\ge \sum_{i=1}^k
i(d_i{+}1)$. On the other hand, using the formula for $\dim\g_e$ in
\cite{ja} we obtain
$$\dim\g_e=\sum_{i=1}^k (2i{-}1)(d_i{+}1)= 2\sum_{i=1}^k
i(d_i{+}1)-n;$$ see also \cite[p.~398]{C}. In view of
Theorem~\ref{summa}(i) we must have equalities throughout, forcing
$\sum_{i=1}^n\deg{^e\!}F_i=(\dim\g_e+n)/2$.

As $\ind\g_e=n$ by \cite{fan}, Theorem~\ref{summa}(i) yields that
the generating set $\{F_i\,|\,\,1\le i\le n\}$ is good for $e$.
Combining Theorem~\ref{summa}(iii) with Theorem~\ref{cd2-gl} shows
that this set is actually {\it very} good for $e$. But then
$\cS(\g_e)^{\g_e}=\,\cS(\g_e)^{G_e}=\,{\mathbb K}[{^e\!}F_1,
\ldots, {^e\!}F_n]$ in view of Theorem~\ref{verygood}.
\end{proof}

\begin{rmk}
The degrees of ${^e\!}F_1,\ldots, {^e\!}F_n$ can be read off from
the Young diagram of $e$, as shown below:
\begin{figure}[htb]
\setlength{\unitlength}{0.017in}
\begin{center}
\begin{picture}(90,80)(0,0)

\put(0,0){\line(1,0){90}}
\put(0,0){\line(0,1){80}}
\put(0,80){\line(1,0){10}}
\put(10,80){\line(0,-1){10}}
\put(10,70){\line(1,0){20}}
\put(30,70){\line(0,-1){20}}
\put(30,50){\line(1,0){30}}
\put(60,50){\line(0,-1){10}}
\put(60,40){\line(1,0){10}}
\put(70,40){\line(0,-1){10}}
\put(70,30){\line(1,0){10}}
\put(80,30){\line(0,-1){10}}
\put(90,0){\line(0,1){20}}
\put(80,20){\line(1,0){10}}

\qbezier[42](10,0),(10,35),(10,70)
\qbezier[42](20,0),(20,35),(20,70)

\qbezier[30](30,0),(30,25),(30,50)
\qbezier[30](40,0),(40,25),(40,50)
\qbezier[30](50,0),(50,25),(50,50)

\qbezier[24](60,0),(60,20),(60,40)
\qbezier[18](70,0),(70,15),(70,30)
\qbezier[12](80,0),(80,10),(80,20)

\put(3,73){$\scriptstyle 1$}
\put(3,63){$\scriptstyle 1$}
\put(3,53){$\scriptstyle 1$}
\put(3,43){$\scriptstyle 1$}
\put(3,33){$\scriptstyle 1$}
\put(3,23){$\scriptstyle 1$}
\put(3,13){$\scriptstyle 1$}
\put(3,3){$\scriptstyle 1$}
\put(3,-7){$\scriptstyle 1$}

\put(13.5,63){$\scriptstyle 2$}
\put(13.5,53){$\scriptstyle 2$}
\put(13.5,43){$\scriptstyle 2$}
\put(13.5,33){$\scriptstyle 2$}
\put(13.5,23){$\scriptstyle 2$}
\put(13.5,13){$\scriptstyle 2$}
\put(13.5,3){$\scriptstyle 2$}
\put(13.5,-7){$\scriptstyle 2$}

\put(23.5,63){$\scriptstyle 3$}
\put(23.5,53){$\scriptstyle 3$}
\put(23.5,43){$\scriptstyle 3$}
\put(23.5,33){$\scriptstyle 3$}
\put(23.5,23){$\scriptstyle 3$}
\put(23.5,13){$\scriptstyle 3$}
\put(23.5,3){$\scriptstyle 3$}
\put(23.5,-7){$\scriptstyle 3$}

\put(33.5,43){$\scriptstyle 4$}
\put(33.5,33){$\scriptstyle 4$}
\put(33.5,23){$\scriptstyle 4$}
\put(33.5,13){$\scriptstyle 4$}
\put(33.5,3){$\scriptstyle 4$}
\put(33.5,-7){$\scriptstyle 4$}

\put(42,-7){$\ldots\ldots\ldots$}

\put(83.5,13){$\scriptstyle k$}
\put(83.5,3){$\scriptstyle k$}
\put(83.5,-7){$\scriptstyle k$}
\end{picture}
\end{center}
\caption{}\label{pikcha_A}
\end{figure}
\end{rmk}
\noindent 4.2. In this subsection we give a description of
${^e\!}F_i$ in terms of $\xi_i^{j,s}$. No generality will be lost by
assuming that $h{\cdot}w_i=-d_iw_i$ for $1\le i\le k$ and
$f(e^s{\cdot}w_i)\in\mathbb K(e^{s-1}{\cdot}w_i)$. Then each
$\xi_i^{j,s}$ is an eigenvector for $\ad h$. More precisely, using
our discussion in (3.2) it easy to observe that
\begin{equation}\label{ksi}
(\ad h)(\xi_i^{j,s})\,=\,(d_i-d_j+2s)\xi_i^{j,s}.
\end{equation}
Given a subset $I\subset\{1,\ldots,k\}$ we denote by $|I|$ the
cardinality of $I$. Given a permutation $\sigma$ of
$I=\{i_1,\ldots,i_m\}$ and a nonnegative function $\bar
s\colon\,I\to\mathbb Z_{\ge 0}$ we associate with the triple
$(I,\sigma, \bar s)$ the monomial
$$
\Xi(I,\sigma,\bar s)\,:=\,\,\xi_{i_1}^{\sigma(i_1),\,\bar
s(i_1)}\xi_{i_2}^{\sigma(i_2),\,\bar s(i_2)}
    \ldots\xi_{i_m}^{\sigma(i_m),\,\bar s(i_m)}\in {\mathcal S}(\gt g_e)
$$
of degree $m=|I|$. If $\bar s(i_j)$ does not satisfies the restriction
on $s$ given in (3.2), then we assume that
$\xi_{i_j}^{\sigma(i_j),\bar s(i_j)}=0$.
For every $\Xi=\Xi(I,\sigma,\bar s)$ we denote by
$\lambda(I,\sigma,\bar s)$ the weight of $\Xi$ with respect to $h$.
Obviously, $\lambda(I,\sigma,\bar s)$ is the sum of the $\ad
h$-eigenvalues ($h$-weights)
of the factors
$\xi_{i_j}^{\sigma(i_j),\bar s(i_j)}$. 
\begin{lm}\label{description} For each $\ell\le n$, we have
$$
{^e\!}F_\ell\,=\,\,\,\sum\limits_{|I|\,=\,m,\,\,\lambda(I,\sigma,\bar
s)\,=\,2(\ell-m)} a(I,\sigma,\bar s)\,\Xi(I,\sigma,\bar s)
$$
for some $a(I,\sigma,\bar s)\in\mathbb K$.
\end{lm}
\begin{proof} 1) Fix a basis
$\{y_1,\ldots,y_n\}=
  \{w_1,e{\cdot}w_1,\ldots,e^{d_1}{\cdot}w_1,w_2,
    \ldots,w_k,\ldots,e^{d_k}{\cdot}w_k\}$
of $V$ and let $E_{ij}\in\gt{gl}(V)$ be such that
$E_{ij}(y_k)=\delta_{jk}y_i$ for all $1\le i,j,k\le n$. View
$F_\ell$ as a polynomial in {\it variables} $E_{ij}$ and let $T$ be
a monomial of $F_\ell$ for which $\deg{{^e}T}=\deg{{^e\!}F}_\ell$.
It can be presented as a product $T=T_1\cdots T_k$, where each $T_q$
involves only only those $E_{ij}$ annihilating $\bigoplus_{t\ne
q}\VV[t]$. If $E_{ij}$ is such a variable with $j\ne i-1$, then the
restriction of $E_{ij}$ to $\kappa(\cS_e)$ is either zero or
proportional to some $\xi_q^{u,s}$. Note also that if
$y_i=e^{d_q}{\cdot}w_q$, i.e., if $y_{i+1}\not\in\VV[q]$, then the
restriction of  $E_{i+1,i}$ to $\kappa(\cS_e)$ equals
$\xi_q^{q+1,d_{q+1}}$. So when we restrict $T$ to $\kappa(\cS_e)$,
nonzero constants (terms of degree $0$) will arise only from those
variables $E_{i+1,i}$ with $y_{i+1}\in\VV[q]$. But all such
variables lie under the main diagonal and the monomial $T$ comes
from a principal minor, hence $T_q$ cannot contain only them. Thus,
if $\deg T_q>0$, then either ${T_q}_{|\kappa(\cS_e)}$ is zero or
$\deg {{^e}T_q}\ge 1$.

On the other hand, $\sum\deg{{^e}T_q}=\deg{{^e\!}F_\ell}$ and each
${T_q}_{|\kappa(\cS_e)}$ is nonzero, by our assumption on $T$. Let
$d(T)$ denote the cardinality of $\{q\le k\,|\,\, \deg T_q>0\}$. The
above discussion shows that $\deg{^e}T\ge d(T)$. Since $\deg T_q\le
d_q+1$ and $\sum\deg T_q=\deg F_\ell$, our discussion in (4.1)
yields $\deg {{{^e\!}F}_\ell}\le d(T)$. Hence $\deg{^e}T=d(T)$,
forcing $\deg{{^e}T_j}\le 1$ for all $1\le j\le k$. This means that
each monomial of ${{^e\!}F}_\ell$, when expressed via
$\{\xi_i^{j,s}\}$, has no factors of the form
$\xi_q^{j,s}\xi_q^{i,t}$.

\smallskip

\noindent 2) Let $\Xi=\xi_{i_1}^{j_1,s_1}\ldots\xi_{i_m}^{j_m,s_m}$
be a monomial involved in ${^e\!}F_\ell$. In part~1) we have proved
that all indices $i_1,\ldots,i_m$ are distinct. Let
$I=\{i_1,\ldots,i_m\}$. Suppose there is $j=j_q$ with $j\not\in I$.
Then $\Xi$ has a positive weight with respect to the semisimple
element $\xi_j^{j,0}\in\gt g_e$. But ${{^e\!}F}_\ell$ is invariant
under $\gt g_e$, hence $\Xi$ must be of weight zero. This
contradiction shows that $j\in I$. Similarly, each $i_q$ must be
among $j_1,\ldots,j_m$. In other words, $(j_1,\ldots,j_m)$ is a
permutation of $(i_1,\ldots,i_m)$.

\smallskip

\noindent 3) Since all $\xi_i^{j,s}$ are eigenvectors for $\ad h$,
each monomial $\Xi$ involved in ${^e\!}F_\ell$ has the same weight
as ${^e\!}F_\ell$ itself. Since $F_\ell$ is $h$-invariant and $f$
has weight $-2$, we see that each $\Xi$ has weight $2(\ell-m)$. This
completes the proof.
\end{proof}

\begin{conj}\label{explicit} Up to a nonzero constant,
$$
{^e\!}F_\ell\,=\,\,\,\sum\limits_{|I|\,=\,m,\,\,\lambda(I,\sigma,\bar
s)\,=\,2(\ell-m)} ({\rm sgn}\,\sigma)\,\Xi(I,\sigma,\bar s),
$$
where $m=\deg{^e\!}F_\ell$.
\end{conj}

\smallskip

\noindent
4.3. In this subsection we use the notation of (3.6) and
consider a nilpotent element $e$ of the symplectic Lie algebra
$\g=\widetilde{\g}_0$ (recall that
$\widetilde{\g}_0={\widetilde{\g}}^\sigma$ where
$\widetilde{\g}=\gt{gl}(\VV)$ and $\dim\VV=2n$). It is well-known
that ${\Delta_{2i}}_{|\g}$ with $1\le i\le n$ generate the
invariant algebra $\mathbb K[\g]^{\g}$ and the regular functions
$\Delta_{2i+1}$ vanish on $\g$. For $1\le i\le n$ we denote by
$\delta_i$ the component of minimal degree of the restriction of
$\Delta_{2i}$ to $\cS_e=e+\g_f$. Since $e+\g_f$ is an affine
subspace of $e+\widetilde{\g}_f$ and $\g_e^*$ identifies with the
linear subspace $\kappa_e^*(\g_f)$ of
$\widetilde{\g}_e^*=\kappa_e^*(\widetilde{\g}_f)$, one observes
easily that either $\deg\delta_i=\deg{^e\!}F_{2i}$ or the
restriction of ${^e\!}F_{2i}$ to $\g_e^*$ is zero and
$\deg\delta_i>\deg{^e\!}F_{2i}$.

For $1\le i\le n$ we denote by $\bar{F}_{2i}\in\cS(\g)^\g$ the preimage
of ${\Delta_{2i}}_{|\g}\in\mathbb K[\g]^\g$ under the Killing
isomorphism $\cS(\g)\stackrel{\sim}{\rightarrow}\mathbb K[\g]$. Note
that $\deg{^e\!}\bar{F}_{2i}= \deg\delta_i$ for all $i\le n$.

\begin{thm}\label{deg2} Suppose $\dim\VV=2n\ge 4$ and let
$F_1,\ldots, F_{2n}$ be as above. Then $\{\bar{F}_{2i}\,|\,\,1\le
i\le n\}$ is a very good generating system for any
$e\in\g\cong\gt{sp}_{2n}$ and
$\cS(\g_e)^{\g_e}=\,\cS(\g_e)^{G_e}=\,{\mathbb K}[{^e\!}\bar{F}_2,
\ldots, {^e\!}\bar{F}_{2n}].$ Furthermore,
$\deg {^e\!}\bar{F}_{2i}=\deg {^e\!}F_{2i}$  for all $i\le n$.
\end{thm}
\begin{proof}
From Theorem~\ref{deg1} and the formula for $\dim\gt g_e$ given in
\cite[3.1(3)]{ja} we deduce that
$$
\begin{array}{l}
\dim\gt g_e=\frac{1}{2}\left(\dim\widetilde{\gt g}_e+
\sum\limits_{i,\,d_i \text{
even}}1\right)=\,\sum\limits_{j=1}^{2n}\deg{^e\!}F_j-n
+\sum\limits_{i,\,i'=i+1}1\,. \\
\end{array}
$$
On the other hand, applying Theorem~\ref{deg1} to
$\widetilde{\g}_e$ yields
$$
\begin{array}{l}
\sum\limits_{j=1}^{n}\deg{{^e\!}F}_{2j}=
 \sum\limits_{i,\,d_i \text{ odd}}
  \frac{i(d_i+1)}{2}+\sum\limits_{i,\,i'=i+1}
   \left( i\frac{d_i}{2}+(i+1)\frac{d_i+2}{2}\right)=
\frac{1}{2}\left(\sum\limits_{j=1}^{2n}
     \deg{{^e\!}F}_j\right)+\sum\limits_{i,\,i'=i+1}\frac{1}{2}\,
  .\\
\end{array}
$$
To check this equality one takes the Young diagram of shape
$(d_1+1\ge\cdots\ge  d_k+1)$ with all boxes in the $j$-th column
labelled $j$ (as shown in Figure~\ref{pikcha_A}) and then sums up
all labels assigned to the {\it even} boxes of the diagram (counted
from bottom to top and from left to right). One should also
keep in mind that $d_i=d_{i'}$ for all $i$ and $d_i+1$ is odd
whenever $i'\ne i$. Using the above formulae one obtains
$$2\sum\deg{^e\!}F_{2i}-\dim\g_e\,=
 \sum\limits_{i\le k,\,\,i'=i}\frac{d_i+1}{2}+
 \sum\limits_{i\le k,\,\,i'=i+1}(d_i+1)=n.
$$
Since $\deg\delta_i\ge\deg{^e\!}F_{2i}$ for all $i\le n$, by our
earlier remarks, we now derive
$$\sum_{i=1}^n\deg{^e\!}\bar{F}_{2i}=
\sum_{i=1}^n\deg\delta_i\ge\sum_{i=1}^n\deg{^e\!}F_{2i}=
(\dim\g_e+n)/2.$$ On the other hand, $\{\bar{F}_{2i}\,|\,\,1\le
i\le n\}$ is a generating system for $\cS(\g)^\g$. As $\ind\g_e=\rk
\g=n$ by \cite{fan}, Theorem~\ref{summa}(i) shows that
$\sum_{i=1}^n\deg{^e\!}\bar{F}_{2i}\le(\dim\g_e+n)/2$. Hence
$\deg{^e\!}\bar{F}_{2i}=\deg{^e\!}F_{2i}$ for all $i$ and
$\{\bar{F}_{2i}\,|\,\,1\le i\le n\}$ is a good generating system for
$e$. Combining Theorem~\ref{summa}(iii) with Theorem~\ref{cd2-sp},
we now see that the generating set $\{\bar{F}_{2i}\,|\,\,1\le i\le
n\}$ is very good for $e$. Then Theorem~\ref{verygood} yields
$\cS(\g_e)^{\g_e}=\,\cS(\g_e)^{G_e}=\,{\mathbb K}[{^e\!}\bar{F}_2,
\ldots, {^e\!}\bar{F}_{2n}]$, completing the proof.
\end{proof}

\smallskip

\noindent 4.4. Now suppose $\gt g=\gt{so}(\VV)$. Recall that $\gt g$
is a symmetric subalgebra of $\widetilde{\gt g}=\gt{gl}(\VV)$. Let
$F_1,\ldots,F_n$ be as in (4.1) and set $\bar F_i:={F_i}_{|\gt
g^*}$. If $n=\dim\VV$ is odd, then the set $\{\bar F_{2i}\mid
0<i<n/2\}$ is a basis of $\cS(\gt g)^G$. If $n$ is even, then $\bar
F_n=P^2$, where $P$ is the {\it pfaffian}. Clearly,
$({^e\!}P)^2={{^e\!}\bar F}_{n}$. Similar to the symplectic case, we
have $\deg {^e\!}\bar F_{2i}\ge\deg {^e\!} F_{2i}$. From
\cite[3.1(3)]{ja} it follows that
\begin{equation}\label{dim-centr}
\begin{array}{l}
\dim\gt g_e=\frac{1}{2}\left(\dim\widetilde{\gt g}_e-
   \sum\limits_{i,\,d_i \text{ even}}1\right).
\\
\end{array}
\end{equation}
Note that $l=\rk\g=[(\dim\VV)/2]$. In order to compute
$\sum_{i=1}^l\deg{^e\!}F_{2i}$ we again consider our labelled Young
diagram  (see Figure~\ref{pikcha_A}) and sum up the labels assigned
to the even boxes. It is important to observe that in the present
case neighbouring columns of the same odd size will always have a
different number of even boxes. This is illustrated in
Figure~\ref{pikcha_D}.

\begin{figure}[htb]
\setlength{\unitlength}{0.017in}
\begin{center}
\begin{picture}(60,50)(0,0)

\put(0,0){\line(1,0){60}}
\put(0,0){\line(0,1){50}}
\put(0,50){\line(1,0){10}}
\put(10,50){\line(0,-1){20}}
\put(10,30){\line(1,0){30}}
\put(40,30){\line(0,-1){20}}
\put(40,10){\line(1,0){20}}
\put(60,0){\line(0,1){10}}

\qbezier[24](0,10),(19.5,10),(39,10)
\qbezier[24](0,20),(19.5,20),(39,20)
\qbezier[6](0,30),(5,30),(10,30)
\qbezier[6](0,40),(4.5,40),(9,40)

\qbezier[18](10,0),(10,15),(10,30)
\qbezier[18](20,0),(20,15),(20,30)
\qbezier[18](30,0),(30,15),(30,30)
\qbezier[6](40,0),(40,5),(40,10)
\qbezier[6](50,0),(50,5),(50,10)

\put(3,33){$\scriptstyle 1$}
\put(3,13){$\scriptstyle 1$}
\put(3,-7){$\scriptstyle 1$}

\put(13.5,23){$\scriptstyle 2$}
\put(13.5,3){$\scriptstyle 2$}
\put(13.5,-7){$\scriptstyle 2$}

\put(23.5,13){$\scriptstyle 3$}
\put(23.5,-7){$\scriptstyle 3$}

\put(33.5,23){$\scriptstyle 4$}
\put(33.5,3){$\scriptstyle 4$}
\put(33.5,-7){$\scriptstyle 4$}

\put(43.5,-7){$\scriptstyle 5$}

\put(53.5,3){$\scriptstyle 6$}
\put(53.5,-7){$\scriptstyle 6$}

\end{picture}
\end{center}
\caption{}\label{pikcha_D}
\end{figure}

\noindent Taking into account (\ref{dim-centr}) and the equality
$\sum_{j=1}^n\deg {^e\!}F_j=(\dim\widetilde{\gt g}_e+n)/2$ we now
arrive at the following:
\begin{equation}\label{summa3}
\begin{array}{l}
\sum\limits_{j=1}^{l}\deg
{^e\!}F_{2j}\,=\,\,\sum\limits_{i'=i+1}(2i+1)\frac{d_i+1}{2}\,\,+
\sum\limits_{i=i',\, i \text{ odd}} i\frac{d_i}{2}\,\,
+\sum\limits_{i=i',\, i \text{ even}} i\frac{d_i+2}{2}\\
=\frac{1}{2}\left(\sum\limits_{j=1}^n\deg {^e\!}F_j\,-
     \sum\limits_{i=i',\, i \text{ odd}} i\, +
                   \sum\limits_{i=i',\, i \text{ even}} i \right)\\
    \  \ \ \ \ \ = \frac{1}{2}\left(\dim\gt g_e+\frac{n}{2}\,+
              \sum\limits_{i,\,d_i \text{ even}}\frac{1}{2}\,\,-
     \sum\limits_{i=i',\, i \text{ odd}} i\, +
          \sum\limits_{i=i',\, i \text{ even}} i \right).\\
\end{array}
\end{equation}
\begin{lm}\label{deg-so1}
Let $e$ be a nilpotent element in $\g=\mathfrak{so}(\VV)$ such that
\begin{itemize}
\item[1)] $d_1$ is even;

\smallskip

\item[2)] if $d_{i-1}$ is even for $i$ odd, then $d_{i}$ is
even.
\end{itemize}
Then either $\bar F_2,\bar F_4,\ldots, \bar F_{n-1}$ or $\bar
F_2,\bar F_4,\ldots,\bar F_{n-2},P$ (depending on the parity of $n$)
is a good generating system for $e\in\g$.
\end{lm}
\begin{proof} Let $t_1,\ldots,t_q$ be the indices of the odd-sized
Jordan blocks of $e$.
Recall that there is a decomposition of the set of Jordan blocks of even
size (i.e., those with $d_i$ odd) into pairs $\{i,i'\}$ such that $i'\ne i$
and $d_{i'}=d_i$, see, for example, \cite[Sect.~1]{ja}.
Therefore $t_j$ and $j$ has the same parity. Note also that
$q$ and $n$ have the same parity for any nilpotent element in $\gt{gl}_n$.

First suppose $n$ and $q$ are odd. Recall that $\deg {^e\!}\bar F_{2i}\ge\deg
{^e\!}F_{2i}$ for all $i$. By Theorem~\ref{summa}, we have
$\sum_{i=1}^l\deg {^e\!}\bar F_{2i}\le(\dim\gt g_e+\rk\g)/2$. Due to
(\ref{summa3}) it suffices to prove that
$$
\begin{array}{l}
\frac{n}{2}+\frac{q}{2}-t_1+t_2-t_3+\cdots-t_q\,=\,\frac{n-1}{2}\,.\\
\end{array}
$$
By the assumptions of the lemma, $t_1=1$, $t_3=t_2+1$, $t_5=t_4+1$,
and so on. Thus, $\sum_{i=1}^q (-1)^{i}t_i=-1-(q-1)/2$, which is
exactly what we wanted.

Now suppose $n$ is even. Then $q$ is also even,
and $\deg {^e\!}P\ge(\deg {^e\!}F_n)/2$.
Moreover, since $d_{t_q+1}$ cannot be even and odd at the same time,
we have $t_q=k$, that is the last Jordan block has odd size. As
above, $t_{j+1}=t_j+1$ for all even $1<j<q$. Therefore,
$$
\begin{array}{l}
\sum\limits_{j=1}^{n/2-1}\deg {^e\!}\bar F_{2j}+\deg {^e\!}P\,\ge
\, \frac{1}{2}\left(\dim\gt g_e+\frac{n+q}{2}+\sum\limits_{i=1}^q (-1)^{i}t_i\right)-\frac{k}{2}\\
\enskip\ \ \ =\,\frac{1}{2}\left(\dim
\g_e+\frac{n+q}{2}-1-\frac{q-2}{2}+k-k \right)=
\frac{1}{2}\left(\dim\gt g_e+\frac{n}{2}\right), \\
\end{array}
$$
and we are done.
\end{proof}

\begin{lm}\label{deg-so2} Suppose $\dim \VV=2l$ and let $e$ be a nilpotent element in
$\g=\mathfrak{so}(\VV)$ such that \begin{itemize}
\item[1)] $d_1$ is odd and $d_2=d_1$;

\smallskip

\item[2)] $d_{i}$ is even for $i\ge 3$.
\end{itemize}
Then $e$ admits a good generating system in
$\cS(\g)^\g$.
\end{lm}
\begin{proof}
Recall that $\deg {^e\!}\bar F_{2i}\ge\deg {^e\!}F_{2i}$ for all $i$
and $\deg {^e\!}P\ge(\deg {^e\!}F_{2l})/2$. The even-sized Jordan blocks
of $e$ can be
decomposed into pairs $\{i,i'\}$ with $i'\ne i$. Hence
$k$ is even and it follows from (\ref{summa3}) that
$$
\begin{array}{l}
\sum\limits_{i=1}^{l-1}\deg {^e\!}\bar F_{2i}+\deg {^e\!}P \,\ge\,
  \frac{1}{2}\left(\dim\gt g_e+l+\frac{k-2}{2}+\frac{k-2}{2}\right)-\frac{k}{2}
  \,=\frac{1}{2}\left(\dim\gt g_e+l\right)
 -1.
\\
\end{array}
$$
Thus, the system $\bar F_2,\ldots,\bar F_{r-2},P$ is ``almost
good''. Applying Lemma~\ref{description} we see that in the present
case
$${^e\!}F_{2d_1+2}\,=\,a_1\xi_1^{1,d_1}\xi_2^{2,d_2}+a_2\xi_1^{2,d_1}\xi_2^{1,d_1}$$
for some $a_1,a_2\in\mathbb K$. Since ${^e\!}F_{2d_1+2}$ is
irreducible, being a generator of the polynomial algebra ${\mathcal
S}(\widetilde{\gt g})^{\widetilde{\gt g}}$, it must be that
$a_1a_2\ne 0$. Both $\xi_1^{2,d_1}$ and $\xi_2^{1,d_2}$ vanish on
$\gt g_e$ and so does $\xi_1^{1,d_1}-\xi_2^{2,d_1}$. Therefore,
${^e\!}\bar F_{2d_1+2}=a_1\xi_1^{1,d_1}\xi_2^{2,d_2}$. For
$d=(d_1+1)/2$ we have ${^e\!}\bar F_{2d}=\xi_1^{1,d_1}$, up to a
nonzero scalar. Consequently, ${^e\!}\bar F_{2d_1+2}=c({^e\!}\bar
F_{2d})^2$ for some $c\in\mathbb K^{^\times}$.

If $k>2$, then $d_1+1<l$, and we can replace $\bar F_{2d_1+2}$ by
$\bar F_{2d_1+2}':=\bar F_{2d_1+2}-c\bar F_{d_1+1}^2$. Since
$\deg{^e\!}\bar F_{2d_1+2}'\ge\deg{^e\!}\bar F_{2d_1+2}+1$,
Theorem~\ref{summa}(i) implies that $\bar F_2,\ldots,\bar
F_{2d_1+2}',\ldots,P$ is a good generating system for $e$.

If $k=2$, then ${^e\!}P=c_0 {^e\!}\bar F_{2d}$ for some
$c_0\in\mathbb K^{^\times}$. In this case we can replace $P$ by
$P':=P-c_0\bar F_{d_1+1}$. Then $\deg {^e\!}P'\ge \deg {^e\!}P+1$,
implying that $\bar F_2,\ldots,\bar F_{2r-2},P'$ is a good
generating system for $e$.
\end{proof}

Combining Lemmas~\ref{codim2-so1}, \ref{codim2-so2}, and
\ref{deg-so1} we obtain the following result:

\begin{thm}\label{vggs-so} Let
$\VV$ be an $n$-dimensional vector space over $\mathbb K$, where $n$
is odd, and let $\tilde{\VV}=\VV\oplus{\mathbb K}v$ be as in (3.8).
Let $(d_1+1\ge d_2+1\ge\cdots\ge d_k+1)$ be a partition of $n$ such
that
\begin{itemize}
\item[1)] $d_1$ is even and $k>1$;

\smallskip

\item[2)] $d_{i}$ is odd for all $i\ge 2$.
\end{itemize}Let $e$ and $\hat{e}$ be
nilpotent elements in $\g=\gt{so}(\VV)$ and
$\hat{\g}=\gt{so}(\tilde{\VV})$, respectively, corresponding to the
partitions $(d_1+1,\ldots,d_k+1)$ and $(d_1+1,\ldots,d_k+1,1)$. Then
$e$ and $\hat{e}$ admit very good generating systems in $\cS(\g)^\g$
and $\cS(\hat{\g})^{\hat{\g}}$, respectively, and the invariant
algebras $\cS(\g_e)^{\g_e}$ and
$\cS(\hat{\g}_{\hat{e}})^{\hat{\g}_{\hat{e}}}$ are free.
\end{thm}

\begin{rmk} Conditions of Lemmas~\ref{deg-so1} and \ref{deg-so2}
are only sufficient for the existence of a good generating system.
But we conjecture that the other nilpotent elements in
$\g=\gt{so}(\VV)$ do not possess good generating systems in
$\cS(\g)^\g$.
\end{rmk}

\begin{ex}\label{ex-so}
Now we wish to exhibit a nilpotent element $e$ in $\gt
g=\gt{so}(\VV)$ {\it without} a good generating system in ${\mathcal
S}(\gt g)^\g$. Some details will be left to the reader. Let
$e\in\gt{so}_{12}$ be a nilpotent element with partition
$(5,3,2,2)$. Then $\dim\gt g_e=18$, $\,\ind\gt g_e=6,$ but
$$\sum_{i=1}^5\deg {^e\!}F_{2i}+(\deg
{^e\!}F_{12})/2=11<(18+6)/2=12.$$ One can show that $\deg
{^e\!}\bar F_{2i}=\deg {^e\!}F_{2i}$ and $\deg {^e\!}P=2$. We have
only two ${^e\!}\bar F_{2i}$'s of degree one, but the centre of
$\gt g_e$ is $3$-dimensional and ${^e\!}\bar F_8=a^2$, where $a$
is a central element of $\g_e$ linear independent of ${^e\!}\bar
F_2$ and ${^e\!}\bar F_4$. Moreover, up to a scalar ${^e\!}\bar
F_{10}=a\!\cdot\!{^e\!}P$. We see that ${^e\!}\bar F_{2i}$'s and
${^e\!}P$ are algebraically dependent. On the other hand,
computations show that there is no good way to modify the system
of generators $\bar F_2,\bar F_4,\bar F_6,\bar F_{10},P$ of
$\cS(\g)^\g$.
\end{ex}

\noindent 4.5. Suppose that $\rk \g\ge 2$. Our next goal in this
section is to attack Conjecture~\ref{conj:main} for the elements of
the minimal nilpotent orbit $\co_{\min}=G{\cdot}\tilde{e}$ in $\g$.
More precisely, we are going to show that if $\g$ is not of type
${\bf E}_8$, then $\tilde{e}$ admits a good generating system in
$\cS(\g)^\g$. Thanks to Theorem~\ref{min-nilp:sing} and
Theorem~\ref{verygood} this will reduce verifying
Conjecture~\ref{conj:main} for the elements in $\co_{\min}$ to the
case where $\g$ is of type ${\mathbf E}_8$. Some partial results on
the ${\bf E}_8$ case are obtained in (4.8) where
Conjecture~\ref{conj:main} for $\co_{\min}$ is reduced to a
computational problem on polynomial invariants for the Weyl group of
type ${\bf E}_7$.

We adopt the notation introduced in (3.9) and (3.10), choose a
Cartan subalgebra $\tilde{\gt t}$ of $\g$ contained in $\g(0)$,
and denote by $\Phi$ the root system of $\g$ with respect to
$\tilde{\gt t}$. Choose a positive system $\Phi_+$ in $\Phi$ such
that for every $\gamma\in \Phi_+$ the root subspace
$\g_\gamma=\mathbb K e_\gamma$ is contained in the parabolic
subalgebra $\p:=\g(0)\oplus\g(1)\oplus\g(2)$. Note that
$\Phi=\sqcup_{-2\le i\le 2}\,\Phi_i$ where
$\Phi_i:=\{\gamma\in\Phi\,|\,\,\g_\gamma\subset\g(i)\}$. Clearly,
$\Phi_2=\{\tilde{\alpha}\}$ where $\tilde{\alpha}$ is the longest
root in $\Phi_+$. No generality will be lost by assuming that
$\tilde{e}=e_{\tilde{\alpha}}$ and
$\tilde{f}=e_{-\tilde{\alpha}}$. Set $\gt t:=\Ker\tilde{\alpha}$.
It is well-known (and easy to see) that $\gt t$ is a Cartan
subalgebra in $\g_{\tilde{e}}$ and $\g_{\tilde{e}}\,=\,\gt t\oplus
\bigoplus_{\gamma\in\Phi_i,\,i\ge 0}\, \g_\gamma.$ For
$\beta\in\sqcup_{i\ge 0}\,\Phi_i$ we denote by $\xi_\beta$ the
linear function on $\g_{\tilde{e}}$ that vanishes on $\gt t$ and
has the property that $\xi_\beta(e_\gamma)=\delta_{\beta\gamma}$
for all $\gamma\in\sqcup_{i\ge 0}\,\Phi_i$. The dual space $\gt
t^*$ will be identified with the subspace of $\g_{\tilde{e}}^*$
consisting of all linear functions $\xi$ such that
$\xi(e_\gamma)=0$ for all $\gamma\in\sqcup_{i\ge 0}\,\Phi_i$. Set
$\gt h:=\gt t\oplus \,\mathbb K\tilde{e}$, an abelian subalgebra
of $\g_{\tilde{e}}$. We regard $\gt h^*=\gt t^*\oplus \,\mathbb
K\xi_{\tilde{\alpha}}$ as a subspace of $\g_{\tilde{e}}^*$.

Choose $\xi_0\in\gt t^*$ such that
$\xi_0([\g_\gamma,\g_{-\gamma}])\ne 0$ for all $\gamma\in\Phi_0$
and put $\eta:=\xi_0+\xi_{\tilde{\alpha}}$, an element of $\gt
h^*$. Since $\eta$ vanishes on $\g(1)$, it is immediate from our
discussion in (3.10) that $\g_{\tilde{e}}^\eta=\gt h$. In
particular, $\eta\in(\g_{\tilde{e}}^*)_{\rm reg}$.  Moreover, our
earlier remarks show that
$$\gt h^*=\{\xi\in\g_{\tilde{e}}^*\,|\,\,({\rm ad}^*\,\gt h)\xi=0\} \quad
\mbox{and}\quad \gt h^*\cap ({\rm ad}^*\,\gt g_{\tilde{e}})\gt
h^*=0.$$ It follows that the differential of the coadjoint action
morphism $G_{\tilde{e}}\times {\gt h}^*\rightarrow
\g_{\tilde{e}}^*$ is surjective at $1\times\eta$. Then
$\overline{G_{\tilde{e}}\cdot \gt h^*}=\g_{\tilde{e}}^*$, implying
that the restriction map $\varphi\mapsto \varphi_{|\gt h^*}$
induces an embedding $\mathbb
K[\g_{\tilde{e}}^*]^{G_{\tilde{e}}}\hookrightarrow \mathbb K[\gt h^*]$.
Hence, for every nonzero homogeneous $F\in \cS(\g)^\g$ the regular
function $^{\tilde{e}\!}F_{|{\gt h^*}}$ is nonzero and thus has
the same degree as $^{\tilde{e}\!}F$.

\smallskip

\noindent 4.6. The Weyl group $W=N_G(\tilde{\gt t})/Z_G(\tilde{\gt
t})$ is generated by the orthogonal reflections $s_\gamma$ in the
hyperplanes $\Ker\gamma$, where $\gamma\in\Phi$. Let
$C_W(\tilde{h})$ be the stabiliser of $\tilde{h}$ in $W$. It is
well-known that $C_W(\tilde{h})=\langle
s_\gamma\,|\,\,\gamma(\tilde{h})=0\rangle$. Obviously,
$C_W(\tilde{h})$ preserves $\gt t$. We denote by $\rho_0$ the
corresponding representation of $C_W(\tilde{h})$ and put
$W_0:=\rho_0(C_W(\tilde{h}))$. Note that $W_0$ is a finite
reflection subgroup of $\GL(\gt t)$. Since $\gt
t=\Ker\,\tilde{\alpha}$ and
$s_{\tilde{\alpha}}(\tilde{h})=-\tilde{h}$, any nonzero
$\varphi\in \cS(\tilde{\gt t})^W$ has the form
\begin{eqnarray}\label{w0}\varphi\,=\,\sum_{i=0}^\nu\,\varphi^{(i)}\tilde{h}^{2i}\ \qquad
\quad\ \quad  \big(\varphi^{(i)}\in \cS(\gt t)^{W_0},\ \ \,
\varphi^{(\nu)}\ne 0,\ \ \,  \nu=\nu(\varphi)\big).\end{eqnarray}
For $\psi\in\cS(\tilde{\gt t})$ we denote by $\tilde{\gt t}_\psi$
the set of all $h\in\tilde{\gt t}$ such that $\psi(x+\lambda
h)=\psi(x)$ for all $x\in \tilde{\gt t}$ and all $\lambda\in
\mathbb K$. If $\psi\in \cS(\tilde{\gt t})^W$, then $\tilde{\gt
t}_\psi$ is a $W$-invariant subspace of $\tilde{\gt t}$. As
$\tilde{\gt t}$ is an irreducible $W$-module, then for $\varphi$
as in (\ref{w0}) we must have $\tilde{\gt t}_\varphi=0$.
Consequently, $\nu(\varphi)\ge 1$.
\begin{prop}\label{weyl}
Let $\{\varphi_1,\varphi_2,\ldots, \varphi_l\}$ be a homogeneous
generating set in $\cS(\tilde{\gt t})^W$ with $\deg\varphi_1=2$, and
$\nu_i=\nu(\varphi_i)$. Then $\sum_{i=2}^l\,\nu_i\ge
\frac{1}{2}\dim\g(1)$ and if $\sum_{i=2}^l\,\nu_i=
\frac{1}{2}\dim\g(1)$, then
$\varphi_2^{(\nu_2)},\ldots,\varphi_l^{(\nu_l)}$ are algebraically
independent and $\cS(\g)^\g$ admits a good generating system for
$\tilde{e}$.
\end{prop}
\begin{proof}
Consider the Levi subalgebra $\tilde{\gt s}=\mathbb
K\tilde{f}\oplus\tilde{\gt t}\oplus\mathbb K\tilde{e}$ of $\g$ and
put $c:=\tilde{h}^2+4\tilde{e}\tilde{f}$, an element of
$\cS(\tilde{\gt s})$. Since $\gt z(\tilde{\gt s})=\gt t$  and
$[\tilde{\gt s},\tilde{\gt s}]=\mathbb K\tilde{f}\oplus\mathbb K
\tilde{h}\oplus\mathbb K\tilde{e}$, we have that $\cS(\tilde{\gt
s})^{\tilde{\gt s}}\cong \cS(\gt t)\otimes_{\mathbb K}\mathbb
K[c]$ as algebras. We identify $\tilde{\gt s}^*$ with $\mathbb
K\xi_{-\tilde{\alpha}}\oplus\gt t^*\oplus \mathbb
K\xi_{\tilde{\alpha}}$. Then $\gt h^*=\gt t^*\oplus\mathbb
K\xi_{\tilde{\alpha}}\subset \tilde{\gt s}^*$. Since $\tilde{e}$
is regular nilpotent in $\tilde{\gt s}$, the restriction map
$F\mapsto F_{|\xi_{-\tilde{\alpha}}+\gt h^*}$ induces an algebra
isomorphism $\iota\colon\,\cS(\tilde{\gt s})^{\tilde{\gt
s}}\stackrel{\sim}{\longrightarrow} S(\gt h)$ such that
$\iota_{|\gt t}=\id$ and $\iota(c)=4\tilde{e}$.

By the Chevalley Restriction Theorem, there exists a homogeneous
generating system $\{F_1,\ldots, F_l\}\subset\cS(\g)^\g$ such that
${F_k}_{|\tilde{\gt t}^*}=\varphi_k$ for all $k$. Since
$\varphi_k=\,\sum_{i=0}^{\nu(\varphi_k)}\varphi_{k}^{(i)}\tilde{h}^{2i}$
by (\ref{w0}), it follows that ${F_k}_{|\tilde{\gt
s}^*}=\,\sum_{i=0}^{\nu(\varphi_k)}\varphi_{k}^{(i)}c^{i}$. But
then $\iota({F_k}_{|\tilde{\gt
s}^*})=\,\sum_{i=0}^{\nu(\varphi_k)}4^i\varphi_{k}^{(i)}\tilde{e}^{i}$. It
is now immediate from the definition of $^{\tilde{e}\!}F$  that
\begin{eqnarray}\label{12}(^{\tilde{e}\!}{F_k})_{|\gt
h^*}=\,4^{\nu_k}\varphi_{k}^{(\nu_{k})}\tilde{e}^{\nu_k}\qquad\quad\qquad
(1\le k\le l).\end{eqnarray} Since $\deg {^{\tilde{e}\!}F_k}=\deg
{(^{\tilde{e}\!}{F_k})_{|\gt h^*}}$ by our concluding remark in
(4.5), Theorem~\ref{summa}(i) in conjunction with (\ref{12}) gives
$$\sum_{i=1}^l\,\deg
\varphi_i^{(\nu_i)}+\sum_{i=1}^l\,\nu_i\le
(\dim\g_{\tilde{e}}+l)/2.$$ On the other hand, (\ref{w0}) shows that
$$\sum_{i=1}^l\,\deg
\varphi_i^{(\nu_i)}+2\sum_{i=1}^l\,\nu_i\,=\sum_{i=1}^l\,\deg\varphi_i\,=\sum_{i=1}^l\,\deg
F_i =(\dim\g+l)/2.$$ As $\dim\g-\dim \g_{\tilde{e}}=2+\dim\g(1)$ and
$\nu_1=1$ by our assumption on $\deg \varphi_1$, we are now able to
conclude that $\sum_{i=2}^l\,\nu_i\ge \frac{1}{2}\dim\g(1)$.

If $\sum_{i=2}^l\,\nu_i= \frac{1}{2}\dim\g(1)$, then the above shows
that
$$\sum_{i=1}^l\deg{^{\tilde{e}\!}F_i}\,=\sum_{i=1}^l\,\deg
\varphi_i^{(\nu_i)}+\sum_{i=1}^l\,\nu_i= (\dim\g_{\tilde{e}}+l)/2.$$
Hence $\{F_1,\ldots, F_l\}\subset \cS(\g)^\g$ is a good generating
system for $\tilde{e}$, implying that $^{\tilde{e}\!}F_1,
{^{\tilde{e}\!}F_2},\ldots, {^{\tilde{e}\!}F_l}$ are algebraically
independent; see Theorem~\ref{summa}(ii). As $\varphi_1^{(\nu_1)}$
is a nonzero constant, our discussion in (4.5) together with
(\ref{12}) shows that $e,\,
\varphi_2^{(\nu_2)}\tilde{e}^{\nu_2},\,\ldots,\,\varphi_l^{(\nu_l)}\tilde{e}^{\nu_l}$
are algebraically independent in $\cS(\g_{\tilde{e}})$. Then
$\varphi_2^{(\nu_2)},\ldots,\varphi_l^{(\nu_l)}$ must be
algebraically independent in $\cS(\gt t)$. This completes the proof.
\end{proof}

\smallskip

\noindent 4.7. Proposition~\ref{weyl} in conjunction with
Theorems~\ref{summa}(iii), ~\ref{cd2-gl}, ~\ref{cd2-sp},
~\ref{min-nilp:sing} and~\ref{verygood} will enable us to show that
$\cS(\g_{\tilde{e}})^{\g_{\tilde{e}}}$ is a graded polynomial
algebra in $\rk\g$ variables in all cases except when $\g$ is of
type ${\mathbf E}_8$. We shall identify $\cS(\tilde{\gt t})$ with
$\cS(\tilde{\gt t}^*)$ by means of the $W$-invariant scalar product
$(\,\cdot\,|\,\cdot\,)$ used in \cite{B} and \cite{VO}. Note that
$\tilde{h}=\tilde{\alpha}^\vee$ identifies with a nonzero multiple
of $\tilde{\alpha}$. The basis of simple roots contained in $\Phi_+$
will be denoted by $\Delta$.

\smallskip

\noindent (1) Suppose $\g$ is of type ${\mathbf A}_n$, $n\ge 2$.
Then $\tilde{\gt t}^*$ is spanned by
$\varepsilon_1,\varepsilon_2,\ldots,\varepsilon_{n+1}$ subject to
the relation
$\varepsilon_1+\varepsilon_2+\cdots+\varepsilon_{n+1}=0$. The Weyl
group $W={\mathfrak S}_{n+1}$ permutes the $\varepsilon_i$'s. Put
$$s_k:=\sum_{\sigma\in{\mathfrak
S}_{n+1}}\,\varepsilon_{\sigma(1)}\varepsilon_{\sigma
(2)}\cdots\varepsilon_{\sigma(k)}\qquad\quad\quad\, (2\le k\le
n+1).$$ Since $\tilde{\alpha}=\varepsilon_1-\varepsilon_{n+1}$ and
$(\tilde{\alpha}|\varepsilon_i)=0$ for $2\le i\le n$, it is
routine that $\nu(s_k)=1$ for $2\le k\le n+1$. Now set
$\varphi_k:=s_{k+1},\,$ $1\le k\le n$. Then
$\{\varphi_1,\varphi_2,\ldots,\varphi_n\}$ is a homogeneous
generating set in $\cS(\tilde{\gt t}^*)^W$ with $\deg
\varphi_1=2$. Since
$\sum_{i=2}^{n}\deg\nu(\varphi_i)=n-1=\frac{1}{2}\dim \g(1),$ we
derive that $\cS(\g_{\tilde{e}})^{\g_{\tilde{e}}}$ is a graded
polynomial algebra in $n$ variables. The degrees of basic
invariants are $1$, $2$, $\ldots\,$, $n$. Since $\g=\gt{sl}_{n+1}$
and the partition of $\tilde{e}$ is $(2,1^{n-1})$, this is
consistent with the combinatorial description in
Theorem~\ref{deg1}.

\smallskip

\noindent (2) Suppose $\gt g$ is of type ${\mathbf C}_n$, $n\ge
2$. Then $\tilde{\alpha}=2\varepsilon_1$, and we can assume that
$\varphi_k=\tilde{s}_k$, where
\begin{eqnarray}\label{sp-weyl}\tilde{s}_k:=\sum_{\sigma\in{\mathfrak
S}_{n+1}}\,\varepsilon_{\sigma(1)}^2\varepsilon_{\sigma
(2)}^2\cdots\varepsilon_{\sigma(k)}^2\qquad\quad\quad\, (1\le k\le
n).\end{eqnarray} As $(\tilde{\alpha}|\varepsilon_i)=0$ for $2\le
i\le n$, it is clear that $\nu(\varphi_k)=1$ for all $k$. Then
$\sum_{i=2}^{n}\deg\nu(\varphi_i)=n-1=\frac{1}{2}\dim \g(1),$
which shows that $\cS(\g_{\tilde{e}})^{\g_{\tilde{e}}}$ is a
graded polynomial algebra in $n$ variables. The degrees of basic
invariants are $1$, $3$, $\ldots\,$, $2n-1$. Since
$\g=\gt{sp}_{2n}$ and the partition of $\tilde{e}$ is
$(2,1^{2n-2})$, this is consistent with our description in
Theorem~\ref{deg2}.

\smallskip

\noindent (3) Suppose $\gt g$ is of type ${\mathbf B}_n$, $n\ge
3$. Then $\tilde{\alpha}=\varepsilon_1+\varepsilon_2$. For
$k\in\{1,3,\ldots,n\}$ put $\varphi_k:=\tilde{s}_k$, where
$\tilde{s}_k$ is as in (\ref{sp-weyl}), and set
$\varphi_2:=\tilde{s}_2-\frac{1}{4}\tilde{s}_1^2$. As
$(\tilde{\alpha}|\varepsilon_i)=0$ for $3\le i\le n$, it is
straightforward to see that $\nu(\varphi_2)=1$  and
$\nu(\varphi_k)=2$ for $3\le k\le n$. Then
$\sum_{i=2}^{n}\deg\nu(\varphi_i)=1+2(n-2)=\frac{1}{2}\dim \g(1).$
Hence $\cS(\g_{\tilde{e}})^{\g_{\tilde{e}}}$ is a graded
polynomial algebra in $n$ variables, and the degrees of basic
invariants are $1$, $3$, $4$, $\ldots\,$, $2n-2$.

\smallskip

\noindent (4) Suppose $\gt g$ is of type ${\mathbf D}_n$, $n\ge
4$. Then again $\tilde{\alpha}=\varepsilon_1+\varepsilon_2$. For
$k\in\{1,3,\ldots,n-1\}$ put $\varphi_k:=\tilde{s}_k$ and set
$\varphi_2:=\tilde{s}_2-\frac{1}{4}\tilde{s}_1^2$. Finally, set
$\varphi_n:=p$ where $p=\prod_{i=1}^n\varepsilon_i$. As in
part~(3) we obtain $\nu(\varphi_2)=1$ and $\nu(\varphi_k)=2$ for
$3\le k\le n-1$. Since $\nu(\varphi_n)=1$, we have
$\sum_{i=2}^{n}\deg\nu(\varphi_i)=1+2(n-3)+1=\frac{1}{2}\dim
\g(1).$ Thus, $\cS(\g_{\tilde{e}})^{\g_{\tilde{e}}}$ is a graded
polynomial algebra in $n$ variables, and the degrees of basic
invariants are $1$, $3$, $4$, $\ldots\,$, $2n-4$, $n-1$.

\smallskip

\noindent (5) Suppose $\gt g$ is of type ${\mathbf G}_2$ and
assume that $\Delta=\{\alpha,\beta\}$ where $\beta$ is a short
root. Then $\tilde{\alpha}=2\alpha+3\beta$ and
$(\beta|\tilde{\alpha})=0$. The degrees of basic invariants in
$\cS(\tilde{\gt t}^*)^W$ are $2$, $6$. There exists
$\varphi_1\in\cS(\tilde{\gt t}^*)^W$ such that
$\varphi_1=\tilde{\alpha}^2+\lambda_0\beta^2$ for some
$\lambda_0\in \mathbb K$. Since $\deg \varphi_1^3=6$, we can find
a basic $W$-invariant $\varphi_2$ in $\cS^{6}(\tilde{\gt t}^*)$
such that
$\varphi_2=\lambda_1\tilde{\alpha}^4\beta^2+\lambda_2\tilde{\alpha}^2\beta^4+
\lambda_3\beta^6$ for some
$\lambda_1,\lambda_2,\lambda_3\in\mathbb K$. Then
$\nu(\varphi_2)\le 2=\frac{1}{2}\dim\g(1)$. Applying
Proposition~\ref{weyl} yields
$\nu(\varphi_2)=\frac{1}{2}\dim\g(1)$. Then
$\cS(\g_{\tilde{e}})^{\g_{\tilde{e}}}$ is a graded polynomial
algebra in two variables, and the degrees of basic invariants are
$1$, $4$.

\smallskip

\noindent (6) Suppose $\gt g$ is of type ${\mathbf F}_4$. In this
case $\tilde{\alpha}=\varepsilon_1+\varepsilon_2$ and basic
invariants in $\cS(\tilde{\gt t}^*)^W$ have degrees $2$, $6$, $8$,
$12$. Let $W'$ denote the subgroup of $W$ generated all
reflections $s_\alpha$ corresponding to long roots in $\Phi$. The
reflection group $W'$ has type ${\mathbf D}_4$ and acts on the
$\varepsilon$-basis of $\tilde{\gt t}^*$ in the standard way.
Therefore, $\cS(\tilde{\gt t}^*)^{W'}=\mathbb
K[\tilde{s}_1,\tilde{s}_2,\tilde{s}_3,p]$ where $\tilde{s}_1,
\tilde{s}_2, \tilde{s}_3, p$ are as in part~(4). Note that $W'$ is
a normal subgroup of $W$ and $W/W'\cong\,{\mathfrak S}_3$.

Set $\varphi_1=\tilde{s}_1$. It is easy to see that
$\varphi_1\in\cS(\tilde{\gt t}^*)^W$. Since
$\varphi_1^3\in\cS^6(\tilde{\gt t}^*)^{W}$ and $\nu(\varphi_1^3)=3$,
there exists a basic invariant $\varphi_2\in \cS^6(\tilde{\gt
t}^*)^{W}$ for which $\nu(\varphi_2)\le 2$. Next observe that
$M:=\cS^4(\tilde{\gt t}^*)^{W'}$ is a $W/W'$-module with basis
$\{\tilde{s}_2,p,\tilde{s}_1^2\}$. We denote by $M'$ the submodule
of $M$ spanned by all $(w-1){\cdot}m$ with $w\in W$ and $m\in M$.
Let $\beta:=\varepsilon_1$ and
$\gamma:=\frac{1}{2}(\varepsilon_1+\varepsilon_2+\varepsilon_3+\varepsilon_4)$,
short roots in $\Phi$, and put
$p':=\textstyle{\prod_{i=1}^4}\,(\varepsilon_i-\gamma).$ Since
$p'=s_\gamma(\prod_{i=1}^4\,\varepsilon_i)=s_\gamma(p)$ and
$s_\beta(p)=-p$, we have $p,p'\in M'$. Since
$s_\beta(\tilde{s}_2)=\tilde{s}_2$ and $\cS^4(\tilde{\gt
t}^*)^{W}=\mathbb K \tilde{s}_1^2$, this shows that $M'$ is
isomorphic to the reflection module for $W/W'\cong{\mathfrak S}_3$,
and $p$ and $p'$ form a basis for $M'$.

The above discussion implies that there exist homogeneous
polynomials $q_2,q_3\in\mathbb K[X,Y]$ of degree $2$ and $3$,
respectively, such that $q_2(p,p')$ and $q_3(p,p')$ generate the
invariant algebra $\mathbb K[M']^{{\mathfrak S}_3}\subset
\cS(\tilde{\gt t}^*)^{W}$. As
$(\tilde{\alpha}|\gamma)=(\tilde{\alpha}|\varepsilon_1)=(\tilde{\alpha}|\varepsilon_2)=1$,
one checks easily that $\nu(p)=\nu(p')=1$. Hence $\nu(q_2(p,p'))\le
2$ and $\nu(q_3(p,p'))\le 3$. Since $\cS^6(\tilde{\gt t}^*)^{W'}$ is
spanned by $\tilde{s}_3$, $\tilde{s}_1\tilde{s}_2$, $\tilde{t}_1p$,
$\tilde{s}_1^3$, there are $\lambda_1,
\lambda_2,\lambda_3,\lambda_4\in\mathbb K$ such that
$\varphi_2=\lambda_1\tilde{s}_3+\tilde{s}_1(\lambda_2p+\lambda_3p')+
\lambda_4\tilde{s}_1^3$. As $(M')^W=0$, it must be that
$\lambda_1\ne 0$. From this it is immediate that $$\cS(\tilde{\gt
t}^*)^{W'}\cong \,{\mathbb K}[\varphi_1,\varphi_2]\otimes_{\mathbb
K}{\mathbb K}[M']$$ as $W$-modules. But then we can set
$\varphi_3:=q_2(p,p')$ and $\varphi_4:=q_3(p,p')$ to obtain a
generating set
$\{\varphi_1,\varphi_2,\varphi_3,\varphi_4\}\subset\cS(\tilde{\gt
t}^*)^{W}$ with $\deg\varphi_1=2$ and
$\sum_{i=2}^4\,\nu(\varphi_i)\le 2+2+3=7$. Since in the present case
$\dim\g(1)=14$, Proposition~\ref{weyl} shows that
$\nu(\varphi_2)=\nu(\varphi_3)=2$ and $\nu(\varphi_4)=3$. Hence
$\cS(\g_{\tilde{e}})^{\g_{\tilde{e}}}$ is a graded polynomial
algebra in four variables, and the degrees of basic invariants are
$1$, $4$, $6$, $9$.

\smallskip

\noindent (7) Now suppose $\gt g$ is of type ${\mathbf E}_6$ and
let $\sigma$ denote the outer involution in ${\rm Aut}(\Phi)$
preserving $\Delta$. In the present case, the degrees of basic
invariants in $\cS(\tilde{\gt t}^*)^{W}$ are $2$, $5$, $6$, $8$,
$9$, $12$. The reflection group $W_0$ has type ${\mathbf A}_5$ and
basic invariants in $\cS(\gt t^*)^{W_0}$ have degrees $2$, $3$,
$4$, $5$, $6$. We choose a homogeneous generating system
$\{\psi_1,\ldots, \psi_5\}\subset\cS(\gt t^*)^{W_0}$ with
$\deg\psi_i=i+1$ for $1\le i\le 5$. Since
$\sigma(\tilde{\alpha})=\tilde{\alpha}$, both $W_0$ and $\gt t$
are $\sigma$-stable. Set ${\tilde{\gt t}}^\sigma:=\{t\in\tilde{\gt
t}\,|\,\,\sigma(t)=t\}$ and $\gt t^\sigma:=\tilde{\gt
t}^\sigma\cap {\gt t}$.  The groups $W^\sigma=\{w\in
W\,|\,\,\sigma w=w\sigma\}$ and $W_0^\sigma=W^\sigma\cap W_0$ act
on $\tilde{\gt t}^\sigma$ and $\gt t^\sigma$, respectively, and we
shall denote by $\tilde{\rho}$ and $\rho$ the corresponding
representations. It is well-known that $\tilde{\rho}(W^\sigma)$
and $\rho(W_0^\sigma)$ are reflection groups of type ${\mathbf
F}_4$ and ${\mathbf C}_3$, respectively.

Note that $\tilde{\gt t}^\sigma=\gt t^\sigma\oplus\mathbb
K\tilde{h}$. To make use of the results obtained in part~(6) we
shall restrict functions from $\cS(\tilde{\gt t}^*)^W$ to
$\tilde{\gt t}^\sigma$. Let $\bar{\psi}_i$ denote the restriction
of $\psi_i$ to $\gt t^\sigma$. Since $\rho(W_0^\sigma)$ is a
reflection group of type ${\mathbf C}_3$, we have that
$\bar{\psi_2}=\bar{\psi_4}=0$ and $\mathbb K[\gt
t^\sigma]^{W_0^\sigma}= \mathbb
K[\bar{\psi_1},\bar{\psi_3},\bar{\psi}_5]$.

Observe that $\dim\cS^5(\tilde{\gt t}^*)^{W}=1$. Let
$\tilde{\varphi}_2$ be a nonzero element in $\cS^5(\tilde{\gt
t}^*)^{W}$. By our remarks in (4.6) we have
$\nu(\tilde{\varphi}_2)\ge 1$. Thus, it can be assumed that
$\tilde{\varphi}_2=\tilde{\alpha}^2\psi_2+\tilde{\varphi}_2^{(0)}$
where $\tilde{\varphi}_2^{(0)}\in\cS^5({\gt t}^*)^{W_0}$. Clearly,
$\nu(\tilde{\varphi}_2)=1$. Next note that $\dim\cS^9(\tilde{\gt
t}^*)^{W}=2$. As $\nu(\tilde{\varphi}_2\tilde{\varphi}_1^2)=3$, we
can find $\tilde{\varphi}_5\in\cS^9(\tilde{\gt
t}^*)^{W}\setminus{\mathbb K}\tilde{\varphi}_2\tilde{\varphi}_1^2$
for which $\nu(\tilde{\varphi}_5)\le 2$. This element is a basic
invariant of $\cS(\tilde{\gt t}^*)^{W}$.

Let $\{\varphi_1,\varphi_2,\varphi_3,\varphi_4\}\subset \mathbb
K[\tilde{\gt t}^\sigma]^{W^\sigma}$ be the generating set obtained
in part~(6). Choose $\tilde{\varphi}_1\in \cS^2(\tilde{\gt
t}^*)^W$ such that
$\tilde{\varphi}_1=\tilde{\alpha}^2+\tilde{\varphi}_1^{(0)}$ where
$\tilde{\varphi}_1^{(0)}\in\cS^2({\gt t}^*)^{W_0}.$ As
$\tilde{\varphi}_1^3\in \cS^6(\tilde{\gt t}^*)^W$ we can find a
nonzero $\tilde{\varphi}_3\in\cS^6(\tilde{\gt t}^*)^W$ such that
$\tilde{\varphi}_3=\tilde{\alpha}^4a+\tilde{\alpha}^2b+c$ for some
$a,b,c\in\cS({\gt t}^*)^{W_0}.$ Suppose $a=0$. Since
$\nu(\tilde{\varphi}_3)\ge 1$, we then have $b\ne 0$. Since $b$ is
a $W_0$-invariant of degree $4$, it is a polynomial in $\psi_1$
and $\psi_3$. Then $b_{|{\gt t}^\sigma}\ne 0$. Consequently,
$\tilde{\varphi_3}_{|\tilde{\gt
t}^\sigma}=\lambda\varphi_2+\mu\varphi_1^3$ where either
$\lambda\ne 0$ or $\mu\ne 0$. Part~(6) now yields
$\nu(\tilde{\varphi}_3)\ge 2$ forcing $a\ne 0$, a contradiction.
Thus, $\nu(\tilde{\varphi}_3)=2$, and it can be assumed without
loss that $a=\psi_1$.

Next we observe that $\dim \cS^8(\tilde{\gt t}^*)^{W}=3$. Because
$\nu(\tilde{\varphi}_1^4)=4$ and
$\nu(\tilde{\varphi}_3\tilde{\varphi}_1)=3$ by the above, the set
$\cS^8(\tilde{\gt t}^*)^{W}\setminus\{\mathbb
K\tilde{\varphi}_1^4\oplus \mathbb
K\tilde{\varphi}_2\tilde{\varphi}_1\}$ contains an element of the
form $\tilde{\alpha}^4a'+\tilde{\alpha}^2b'+c'$ with
$a',b',c'\in\cS({\gt t}^*)^{W_0}$, say $\tilde{\varphi}_4$. The
element $\tilde{\varphi}_4$ is a basic invariant of
$\cS(\tilde{\gt t}^*)^{W}$. As $\nu(\tilde{\varphi}_1^6)=6$ and
$\nu(\tilde{\varphi}_3\tilde{\varphi}_1^3)=5$, we can find a basic
invariant $\tilde{\varphi}_6\in \cS^{12}(\tilde{\gt t}^*)^{W}$ for
which $\nu(\tilde{\varphi}_6)\le 4$.

Suppose for a contradiction that $a'=0$. In view of our remarks in
(4.6) we then have $b'\ne 0$ and $\nu(\tilde{\varphi}_4)=1$.
Consequently, $$\sum_{i=2}^6\,\nu(\tilde{\varphi}_i)\le
1+2+1+2+4=10.$$ Since in the present case $\dim \g(1)=20$,
Proposition~\ref{weyl} shows that we have equalities everywhere and
the elements $\tilde{\varphi}_i^{(\nu(\tilde{\varphi}_i))}$ with
$2\le i\le 6$ are algebraically independent in $\cS(\gt t^*)^{W_0}$.
But then $\nu(\tilde{\varphi}_5)=2$ and $\nu(\tilde{\varphi}_6)=4$,
which implies that
$\tilde{\varphi}_5^{(\nu(\tilde{\varphi}_5))}\in\mathbb
K\psi_1\psi_2\oplus \mathbb K\psi_4$ and
$\tilde{\varphi}_6^{(\nu(\tilde{\varphi}_6))}\in\mathbb
K\psi_1^2\oplus \mathbb K\psi_3$. As $\nu(\tilde{\varphi}_4)=1$, we
have $\tilde{\varphi}_4^{(\nu(\tilde{\varphi}_4))}=
\mu_1\psi_5+\mu_2\psi_1\psi_3+\mu_3\psi_2^2+\mu_4\psi_1^3$ for some
$\mu_i\in\mathbb K$. Because
$\tilde{\varphi}_2^{(\nu(\tilde{\varphi}_2))},\ldots,
\tilde{\varphi}_6^{(\nu(\tilde{\varphi}_6))}$
are algebraically independent, the above shows that $\mu_1\ne 0$. In
conjunction with our earlier remarks this yields that for
$\tilde{\varphi_4}_{|\tilde{\gt t}^\sigma}\in\mathbb K[\tilde{\gt
t}^\sigma]^{W^\sigma}$ we have
$\nu\big(\tilde{\varphi_4}_{|\tilde{\gt t}^\sigma}\big)=1$. On the
other hand, $\tilde{\varphi_4}_{|\tilde{\gt t}^\sigma}$ is a linear
combination of $\varphi_3$, $\varphi_1\varphi_2$ and $\varphi_1^4$.
Since $\nu(\varphi_3)=2$, $\nu(\varphi_1\varphi_2)=3$ and
$\nu(\varphi_1^4)=4$, this is impossible. Therefore, $a'\ne 0$ and
$\nu(\tilde{\varphi}_4)=2$.

Since $a'$ is a $W_0$-invariant of degree $4$, we have
$a'=\lambda'\psi_3+\mu'\psi_1^2.$ Hence $a'_{|{\gt
t}^\sigma}=\lambda'\bar{\psi}_3+\mu'\bar{\psi}_1^2\ne 0$. It follows
that $\nu\big(\tilde{\varphi_4}_{|\tilde{\gt t}^\sigma}\big)=2$.
Then the above implies that $\tilde{\varphi_4}_{|\tilde{\gt
t}^\sigma}=\eta\varphi_3$ for some $\eta\in\mathbb K^{^\times}$.
Since $\varphi_3^{(\nu(\varphi_3))}=\,\eta^{-1}a'_{|{\gt t}^\sigma}$
is algebraically independent of $\varphi_2^{(\nu(\varphi_2))}$ by
part~(6), we now derive that $\lambda'\ne 0$. Since
$\tilde{\varphi}_6^{(\nu(\tilde{\varphi}_6))}\in\mathbb
K\psi_3\oplus\mathbb K\psi_1^2$, it follows that we can adjust
$\tilde{\varphi}_6$ by a suitable linear combination of
$\tilde{\varphi}_1^2\tilde{\varphi}_4$ and $\tilde{\varphi}_3^2$ to
achieve $\nu(\tilde{\varphi}_6)\le 3$. Then
$$\sum\limits_{i=2}^6\,\nu(\tilde{\varphi}_i)\le 1+2+2+2+3\,=\,\frac{1}{2}\dim\g(1).$$
Proposition~\ref{weyl} now shows that $\nu(\tilde{\varphi}_5)=2$,
$\nu(\tilde{\varphi}_6)=3$, and $\cS(\g)^{\g}$ admits a good
generating system for $\tilde{e}$. Hence
$\cS(\g_{\tilde{e}})^{\g_{\tilde{e}}}$ is a graded polynomial
algebra, and the degrees of basic invariants are $1$, $4$, $4$, $6$,
$7$, $9$.

\smallskip

\noindent (8) Finally, suppose $\g$ is of type ${\mathbf E}_7$. The
degrees of basic invariants in $\cS(\tilde{\gt t})^W$ are $2$, $6$,
$8$, $10$, $12$, $14$, $18$, and our arguments in part~(7) are not
easily adapted to the present situation. Fortunately, this will not
be necessary because a suitable for us system of basic invariants in
$\cS(\tilde{\gt t})^W$ is already recorded in the literature. It has
been constructed in \cite{KM} with the help of computer-aided
calculations.

We have to adopt the notation of \cite{KM}. So let
$\Delta'=\{v_0,v_1,\ldots, v_6\}$ be a basis of the root system
$\Phi$ with the simple roots numbered as follows:
\begin{equation}\label{E7}
{\beginpicture \setcoordinatesystem units <0.45cm,0.3cm>
 \setplotarea x from 0 to 16, y from -4 to 4
 \linethickness=0.03pt
 \put{$\circ$} at 2 0.02
 \put{$\circ$} at 4 0
 \put{$\circ$} at 6 0
 \put{$\circ$} at 8 0
 \put{$\circ$} at 10 0
 \put{$\circ$} at 12 0
 \put{$\circ$} at 6.03 -2.56
 \plot 2.18 .1 3.86 .1 /
 \plot 4.18 .1 5.86 .1 /
 \plot 6 -.2 6 -2.26 /
 \plot 6.18 .1 7.86 .1 /
 \plot 8.19 .1 9.86 .1 /
 \plot 10.19 .1  11.86 .1 /
  \put{$v_1$} at 2.2 1.5
   \put{$v_2$} at 4.2 1.5
    \put{$v_3$} at 6.2 1.5
     \put{$v_4$} at 8.2 1.5
      \put{$v_5$} at 10.2 1.5
       \put{$v_6$} at 12.2 1.5
        \put{$v_0$} at 6.97 -2.73
         \endpicture}
          \end{equation}
Since all roots in $\Phi$ are conjugate under $W$, we may (and will)
assume that $\tilde{\alpha}=v_1$.  Let
$\{v_0^*,v_1^*,\ldots,v_6^*\}$ be the basis of $\tilde{\gt t}$ such
that $v_i(v_j^*)=\delta_{ij}$ for all $0\le i,j\le 6$. As
$(v_1|v_1)=2$, it follows from (\ref{E7}) that
$\tilde{h}=2v_1^*-v_2^*$, whilst our choice of $\tilde{\alpha}$
ensures that $v_i^*\in\Ker\tilde{\alpha}$ for $i\in\{0,2,\ldots,
6\}$. For a root system type ${\mathbf E}_7$, the {\it distinguished
functionals} $t_1,t_2,\ldots, t_7$ are defined in \cite{KM} by the
following formulae:
\begin{eqnarray*}
&&t_1=-\frac{2}{3}v_0^*+v_1^*,\quad \ \ \ \ \,\,
t_2=-\frac{2}{3}v_0^*-v_1^*+v_2^*,\quad
t_3=-\frac{2}{3}v_0^*-v_2^*+v_3^*,\\
&&t_4=\frac{1}{3}v_0^*-v_3^*+v_4^*, \quad
t_5=\frac{1}{3}v_0^*-v_4^*+v_5^*,\quad \ \ \,
t_6=\frac{1}{3}v_0^*-v_5^*+v_6^*,\quad \ \
t_7=\frac{1}{3}v_0^*-v_6^*.
\end{eqnarray*}
We are particularly interested in the basic invariants $A_2, A_6,
A_8, A_{10}, A_{12}, A_{14}, A_{18}$ of $\cS(\tilde{\gt t})^W$
displayed in \cite[Appendix~2]{KM}. These are presented as
polynomials in the elementary symmetric functions $s_1,s_2,\ldots,
s_7$ of the distinguished functionals $t_1,t_2,\ldots, t_7$. The
coefficients of these polynomials are of no importance to us, but we
need to examine the monomials in $s_1, s_2,\ldots, s_7$ that occur
in the $A_k$'s.

Note that $\tilde{\alpha}(t_1)=v_1(t_1)=1$,
$\tilde{\alpha}(t_2)=-v_1(v_1^*)=-1$, and
$\tilde{\alpha}(t_i)=v_1(t_i)=0$ for $3\le i\le 7$. It follows
that $\nu(s_1)=\nu(s_1(t_1,\ldots, t_7))=0$ and
$\nu(s_i)=\nu(s_i(t_1,\ldots, t_7))=1$ for $2\le i\le 7$.
Therefore,
$$\nu(s_{1}^{j_1}s_{2}^{j_2}\cdots
s_{7}^{j_7})=j_2+\cdots+j_7\qquad\ (\forall\, j_k\in\mathbb Z_+,\
1\le k\le 7).$$ Taking this into account and using the explicit
formulae for $A_2, A_6, A_8, A_{10}, A_{12}, A_{14}, A_{18}$ in
\cite[Appendix~2]{KM} one finds out that $\nu(A_2)=1$, $\nu(A_6)\le
2$, $\nu(A_8)\le 2$, $\nu(A_{10})\le 2$, $\nu(A_{12})\le 3$,
$\nu(A_{14})\le 3$ and $\nu(A_{18})\le 4$. It follows that
$$\nu(A_6)+\nu(A_8)+\nu(A_{10})+\nu(A_{12})+
\nu(A_{14})+\nu(A_{18})\le 2+2+2+3+3+4=16.$$ Since in the present
case the derived subalgebra of $\g(0)$ has codimension $1$ in
$\g(0)$ and is isomorphic to $\mathfrak{so}_{12}$, we have
$\frac{1}{2}\dim\g(1)=(\dim
\g-\dim\mathfrak{so}_{12}-3)/4=(133-66-3)/4=16.$
Proposition~\ref{weyl} now shows that
$\nu(A_6)=\nu(A_8)=\nu(A_{10})=2$, $\nu(A_{12})=\nu(A_{14})=3$ and
$\nu(A_{18})=4$. This implies that
$\cS(\g_{\tilde{e}})^{\g_{\tilde{e}}}$ is a graded polynomial
algebra in seven variables, and the degrees of basic invariants are
$1$, $4$, $6$, $8$, $9$, $11$, $14$.

We summarise the results of this subsection:
\begin{cl}\label{notE8}
If $\g$ is not of type $\mathbf{E}_8$, then $\cS(\tilde{\gt t})^W$
contains a homogeneous generating system
$\varphi_1,\varphi_2,\ldots,\varphi_l$  such that $\deg\varphi_1=2$
and $\cS(\gt t)^{W_0}=\, \mathbb
K[\varphi_2^{(\nu_2)},\ldots,\varphi_l^{(\nu_l)}]$ where
$\nu_i=\nu(\varphi_i)$.
\end{cl}
\begin{proof}
We have shown that under the above assumption on $\g$ there exists a
homogeneous system of basic invariants
$\varphi_1,\varphi_2,\ldots,\varphi_l$ in $\cS(\tilde{\gt t})^W$
such that $\deg\varphi_1=2$ and the elements
$\varphi_2^{(\nu_2)},\ldots,\varphi_l^{(\nu_l)}$ are algebraically
independent in $\cS(\gt t)^{W_0}$. So the result follows by
comparing the Hilbert series of the graded polynomial algebra
$\cS(\gt t)^{W_0}$ and its graded subalgebra
$K[\varphi_2^{(\nu_2)},\ldots,\varphi_l^{(\nu_l)}]$.
\end{proof}
\begin{rmk}
If $\g$ is of type ${\bf E}_8$, then one can show by using {\sl ad
hoc} arguments that $\cS(\g_{\tilde{e}})^{\g_{\tilde{e}}}$ contains
an element of degree $4$ linearly independent of $\tilde{e}^4$.
Looking at the degrees of basic invariants in $\cS(\g)^\g$ and
taking into account (\ref{12}) one can observe that this element is
not of the form $^{\tilde{e}\!}F$ with $F\in\cS(\g)^\g$. It follows
that in type ${\bf E}_8$ the elements in $\co_{\min}$ do not admit
good generating systems in $\cS(\g)^\g$. Combining this with
Proposition~\ref{weyl} one obtains that for any homogeneous
generating system $\varphi_1,\varphi_2,\ldots,\varphi_8$ in
$\cS(\tilde{\gt t})^W$ with $\deg\varphi_1=2$ the elements
$\varphi_2^{(\nu(\varphi_2))},\ldots,\varphi_l^{(\nu(\varphi_l))}$
are algebraically dependent in $\cS(\gt t)^{W_0}$. This is in sharp
contrast with Corollary~\ref{notE8}.
\end{rmk}

\smallskip

\noindent 4.8. In this subsection we assume that $\gt g$ is of type
${\mathbf E}_8$, so that $l=\rk\g=8$. We adopt the notation
introduced in (3.9) and (3.10). In particular, $\gt n=\gt
g(1)\oplus\gt g(2)$. As before, we identify $\el^*$ with ${\rm
Ann}(\gt n)\subset\gt g_{\tilde e}^*$ and $\gt g(1)^*$ with ${\rm
Ann}(\el\oplus\gt g(2))$.

In the course of proving Theorem~\ref{min-nilp:sing} we established
that the principal open subset $Y=\g_{\tilde{e}}^*\setminus {\rm
Ann}(\tilde{e})$ of $\g_{\tilde{e}}^*$ decomposes as $Y\cong
\big((N/(N,N)\big)\times\big({\rm Ann}(\gt g(1))\cap Y\big).$ It
follows that restricting regular functions on $Y$ to ${\rm Ann}(\gt
g(1))\cap Y$ we get algebra isomorphisms
$$
\big({\mathcal S}(\gt g_{\tilde e})[1/{\tilde e}]\big)^N
\cong\,{\mathcal S}(\el)[\tilde e,1/\tilde e] \enskip \ \,\text{ and
}\ \enskip{\mathcal S}(\gt g_{\tilde e})^{\gt g_{\tilde
e}}[1/{\tilde e}]\,=\, \big({\mathcal S}(\gt g_{\tilde e})[1/{\tilde
e}]\big)^{\gt g_{\tilde e}}\,\cong\,
 ({\mathcal S}(\el)^L [\tilde e, 1/\tilde e] .
$$

The standard Poisson bracket of $\cS(\g_{\tilde{e}})$ (induced by
Lie product) gives $\cS(\g_{\tilde{e}})[1/\tilde{e}]$ a Poisson
algebra structure. As $\gt n$ is a Heisenberg Lie algebra, the
subspace $\cS^2(\gt g(1))/\tilde e$ is closed under the Poisson
bracket of $\cS(\g_{\tilde e})[1/\tilde{e}]$, i.e.,
$\cS^2(\g(1))/\tilde{e}$ is a Lie subalgebra of $\cS(\gt g_{\tilde
e})[1/\tilde{e}]$. This Lie algebra acts faithfully on $\gt g(1)$
and is isomorphic to $\gt{sp}(\gt g(1))$. Since the bilinear form
$\langle\,\cdot\,,\,\cdot\,\rangle$ is $(\ad \el)$-invariant, $\el$
acts on $\gt g(1)$ as a Lie subalgebra of $\gt{sp}(\gt g(1))$. From
this it follows that for every $x\in\el$ there exists a unique
$\omega(x)\in\cS^2(\g(1))$ for which  $x+\omega(x)/\tilde
e\in\big({\mathcal S}(\gt g_{\tilde e})[1/{\tilde e}]\big)^N$. Since
the restriction of $\omega(x)$ to ${\rm Ann}(\gt g(1))$ is zero,
$x+\omega(x)/\tilde e$ is the preimage of $x$ in $({\mathcal S}(\gt
g_{\tilde e})[1/\tilde e])^N$. It is straightforward to see that the
map $\omega\colon\,\el\to{\mathcal S}^2(\gt g(1))$ is linear.

Let $x_1,\ldots, x_m$ be a basis of $\el$. Given an $L$-invariant
$H=Q(x_1,\ldots, x_m)$ in $\cS(\el)$ we define
$$
\widehat H\,:=\,Q(x_1+\omega(x_1)/{\tilde
e},\ldots,x_m+\omega(x_m)/\tilde e).
$$ Clearly, $\widehat H\in\cS(\g_{\tilde{e}})^{\g_{\tilde{e}}}[1/\tilde{e}]$. Let $k=k(H)$ be the
smallest integer for which $\tilde e^{k}\widehat
H\in\cS(\g_{\tilde{e}})^{\g_{\tilde{e}}}$, and set
$\widetilde{H}:=\tilde e^{k}\widehat H$. Let $\omega(H)$ denote the
``constant term'' of $\widetilde{H}$ with respect to $\tilde e$, so
that $\omega(H)$ equals the restriction of $\widetilde{H}$ to ${\rm
Ann}(\tilde e)$. Note that $k\le\deg H$ and $\deg\omega(H)=\deg\widetilde{H}=
\deg H+k$.

Let $\{H_1,\ldots,H_{l-1}\}$ be a homogeneous generating set for
${\mathcal S}(\el)^L$. Then both
$\{\widehat{H_1},\ldots,\widehat{H_{l-1}}\}$ and
$\{\widetilde{H_1},\ldots,\widetilde{H_{l-1}}\}$ generate the
$\mathbb K[\tilde e,1/\tilde e]$-algebra $({\mathcal S}(\gt
g_{\tilde e})^{\gt g_{\tilde e}}[\tilde{e},1/\tilde e]$.

\begin{lm}\label{E8-1} The algebra
${\mathcal S}(\gt g_{\tilde e})^{\gt g_{\tilde e}}$ is free if and
only if ${\mathcal S}(\el)^L$ contains a homogeneous generating system
$H_1,\ldots,H_{l-1}$ such that the elements
$\omega(H_1),\ldots,\omega(H_{l-1})$ are algebraically independent.
\end{lm}
\begin{proof}
First suppose that ${\mathcal S}(\el)^L$ contains a required set of
generators $H_1,\ldots, H_{l-1}$, and let $H\in{\mathcal S}(\gt
g_{\tilde e})^{\g_{\tilde{e}}}$. Then $H$ is a polynomial in
$\widetilde{H_i}$, $\tilde e$ and $1/\tilde e$, hence can be
presented as a finite sum $H=\sum_{p\in\mathbb Z}\, \tilde e^p Q_p$,
where $Q_i$ are nontrivial polynomials in $\widetilde{H_i}$. Since
$\omega(H_i)$ are algebraically independent by our assumption, all
$Q_i$ are coprime to $\tilde e$. This implies that $H=\sum_{p\ge
0}\, \tilde e^p Q_p$, that is $H$ is a polynomial in
$\widetilde{H_i}$ and $\tilde{e}$.

Now suppose that ${\mathcal S}(\gt g_{\tilde e})^{\gt g_{\tilde e}}$
is a free algebra generated over $\mathbb K$ by $T_1,\ldots,T_l$.
Without loss of generality we may (and will) assume that all $T_i$
are homogeneous and $T_l=\tilde e$. As $({\mathcal S}(\gt g_{\tilde
e})[1/\tilde e])^{\gt g_{\tilde e}} \cong\,{\mathcal
S}(\el)^L[\tilde e,1/\tilde e]$, there exist $H_1,\ldots, H_{l-1}\in
\cS(\el)^L$ and $b_1,\ldots, b_{l-1}\in\mathbb Z$ such that
$T_i=e^{b_i}\widehat{H_i}$ for $1\le i\le l-1$. Moreover,
$H_1,\ldots, H_{l-1}$ generate ${\mathcal S}(\el)^L$. Because
the product $e^{b_i}\widehat{H_i}$ is irreducible and regular, it must be
that $b_i=k_i$. Hence $T_i=\widetilde{H_i}$ for all $i<l$.

Assume for a contradiction that
$P(\omega(H_1),\ldots,\omega(H_{l-1}))=0$ for a nonzero polynomial
$P\in\mathbb K[X_1,\ldots, X_{l-1}]$. Then
$H':=P(T_1,\ldots,T_{l-1})/\tilde e$ is a {\it regular} $\gt
g_{\tilde e}$-invariant. On the other hand, $H'$ is {\it uniquely}
expressed as a polynomial in $T_1,\ldots, T_{l-1}$ with coefficients
in $\mathbb K[\tilde e,1/\tilde e]$, and
$\cS(\g_{\tilde{e}})^{\g_{\tilde{e}}}=\,\mathbb
K[T_1,\ldots,T_{l-1},\tilde{e}]$ by our assumption. But then
$H'\not\in\cS(\g_{\tilde{e}})^{\g_{\tilde{e}}}$. By contradiction,
the result follows.
\end{proof}

It is well-known that in the present case $L$ has type ${\bf E}_7$
and the stationary subgroup $K=L\cap G_e$ is a simple algebraic
group of type ${\bf E}_6$. Recall from (3.9) that $e$ is a generic
point of the $L$-module $\g(1)$ and $K$ is a generic stabiliser in
$L$; see Definition~\ref{generic}. It is also known that $K$ is the
derived subgroup of the intersection of two opposite maximal
parabolics of $L$. More precisely, $K=(L^+\cap L^-,L^+\cap L^-)$,
where $L^+$ (resp., $L^-$) is the the normaliser in $L$ of the line
spanned by a highest (resp., lowest) weight vector of the $L$-module
$\g(1)$. These primitive vectors will be denoted by $e^+$ and $e^-$,
respectively. Note that $[e^+,e^-]$ is a nonzero multiple of
$\tilde{e}$ (equivalently, $\langle e^+,e^-\rangle\ne 0$). Choose a
maximal torus $\hat{\gt t}$ in the Levi subalgebra $\Lie (L^+\cap
L^-)$ of $\el$ and set $\gt t:=\hat{\gt t}\cap \ka$. It is easy to
see that $\gt t$ is a maximal torus in $\ka=\Lie K$.

It follows from the above description that $\gt g(1)^K={\mathbb
K}e^+\oplus\mathbb K e^-$. Hence it can be assumed without loss of
generality that $e=e^{+}+e^-$. Since the nondegenerate
skew-symmetric form $\langle\,\cdot\,,\cdot\,\rangle$ is
$L$-invariant, $\gt g(1)\cong\gt g(1)^*$ as $L$-modules. Set
$w_+:=\langle\,e^+,\,\cdot\,\rangle$,
$w_-:=\langle\,e^-,\,\cdot\,\rangle$, and $v:=\langle\,e,\,\cdot\,\rangle$. As
explained in the proof Theorem~\ref{min-nilp:sing}, the orbit $(\Ad
L)e$ has codimension $1$ in $\g(1)$ and $(\Ad\, L)({\mathbb
K}^{^\times}\!e)=(\Ad\, G(0))e$ is Zariski open in $\g(1)$. Hence
the tangent space $\el{\cdot}v$ (at $v$) to the orbit $L{\cdot}v$
has codimension $1$ in $\g(1)^*=\,\g(0){\cdot}v=\,\hat{\gt
t}{\cdot}v+{\el}{\cdot}v$. As $K$ is reductive and
$(\g(1)^*)^K=\,{\mathbb K}w_+\oplus{\mathbb K}w_-$ is $\hat{\gt
t}$-stable, we have that $\el{\cdot} v\,=\,\mathbb
Kh_0{\cdot}v\oplus V_0$, where $h_0\in\hat{\gt t}$ is orthogonal to
$\ka$ with respect to the Killing form and $V_0=\big\{\langle
x\,,\,\cdot\,\rangle\,|\,\,\,\langle x,e^+\rangle=\langle
x,e^-\rangle=0\big\}$.

As in (4.5), we regard the dual space $\hat{\gt t}^*$ as a subspace
of $\el^*\subset\gt g_{\tilde e}^*$. We identify $\gt t^*$ with the
subspace $\{\gamma\in\hat{\gt t}^*\,|\,\,\gamma(h_0)=0\}$ and view
$v\in\g(1)^*$ as a linear function on $\g_{\tilde{e}}$ vanishing on
$\el\oplus\mathbb K\tilde{e}.\,$ Set $\,W':=N_L(\gt{\hat
t})/Z_L(\hat{\gt t})\,$ and $\,W'_0:=N_K(\gt t)/Z_K(\gt t)\,$ (these
are reflection groups of type ${\bf E}_7$ and ${\bf E}_6$,
respectively).

\begin{lm}\label{E8-2} Let $H_1,\ldots,H_{l-1}$ be a homogeneous generating
set in ${\mathcal S}(\el)^L$. Then the elements
$\omega(H_1),\ldots,\omega(H_{l-1})$ are algebraically independent
if and only if their restrictions to $\gt t^*\oplus\mathbb K v$ are.
\end{lm}
\begin{proof}
Recall that $\widehat{H_i}\in\cS(\g_{\tilde{e}})^{\g_{\tilde{e}}}$
and $\omega(H_i)=\widetilde{H_i}_{|{\rm Ann}(\tilde e)}$ for $1\le
i\le l-1$. It follows that all $\omega(H_i)$ are invariant under the
coadjoint action of the semidirect product $\el\ltimes\gt g(1)$,
where $\gt g(1)$ is considered as a commutative Lie algebra.

By our earlier remarks, the $L$-saturation of $\mathbb K v$ is dense
in $\gt g(1)^*$. Also, for the same $v$, but regarded as an
element of $(\el\ltimes\gt g(1))^*$, we have
$\big(\ad^*\gt g(1)\big) v\cong(\gt l/\ka)^*$. Combining this two
facts we obtain natural embeddings
$$
\mathbb K[\omega(H_1),\ldots,\omega(H_{\ell-1})]\hookrightarrow\,
\mathbb K[\el^*\oplus\mathbb K v]^{\ka\,\ltimes\gt g(1)}
 \hookrightarrow\, \mathbb K[\ka^*\oplus\mathbb K v]^{\ka}
  \hookrightarrow\, \mathbb K[\gt t^*\oplus\mathbb K v].
$$
As the composition of these embeddings is also an embedding, the
result follows.
\end{proof}

Now we wish to express $\omega(H_i)$ in terms of polynomial invariants for
$W'$. Let $\alpha\in\gt g_{\tilde e}^*$ be such that $\alpha(\tilde
e)=1$ and $\alpha(\el\oplus\gt g(1))=0$, and set $$\gt s:=\gt
t^*\oplus\mathbb K v\oplus\mathbb K\alpha.$$ Then the restriction of
$\omega(H_i)$ to ${\gt t^*\oplus\mathbb Kv}$ is equal to the
``constant term'' (with respect to $\tilde e$) of
$\widehat{H_i}_{|\gt s}$. We thus need to describe the restrictions
of $\widehat{H_i}$ to $\gt s$. Let $\hat{\gt t}^{\perp}\subset\el$
be the orthogonal complement to $\hat{\gt t}=\gt t\oplus\mathbb K
h_0$ with respect to the Killing form, so that $\el=\gt
t\oplus\mathbb K h_0\oplus \hat{\gt t}^{\perp}$. Since $\hat{\gt
t}^\perp$ is spanned by root vectors of $\el$ with respect to
$\hat{\gt t}$ and $e=e^{+}+e^-$, it is straightforward to see that
$[[\hat{\gt t}^\perp,e],e]=0$.

\begin{lm}\label{E8-3}
The following statements are true:

\smallskip

\begin{itemize}
\item[(a)\ ]
$(x+\omega(x)/{\tilde e})_{|\gt s}\,=\,0\ $ for all $x\in\hat{\gt
t}^{\perp}$;
\smallskip

\item[(b)\ ]
$(x+\omega(x)/{\tilde e})_{|\gt s}\,=\,x\ $ for all $x\in\gt t$;
\smallskip

\item[(c)\ ]
$(h_0+\omega(h_0)/{\tilde e})_{|\gt s}\,=\,a(e^{+}-e^-)^2/{\tilde
e}\ $ for some $a\in\mathbb K^{^\times}$.
\end{itemize}
\end{lm}
\begin{proof} Let $x\in \el$ and let $\beta=\gamma+\lambda v+\mu\alpha\in\gt s$,
where $\gamma\in\gt t^*$ and $\lambda,\mu\in\mathbb K$. We shall
calculate the value of $x+\omega(x)/{\tilde e}$ at $\beta$. Without
loss of generality we may assume that both $\lambda$ and $\mu$ are
nonzero. 
Since $x+\omega(x)/\tilde e$ is
$N$-invariant, we can replace $\beta$ by
$\big(\Ad^*(\exp\,\frac{\lambda}{\mu}\,\ad e)\big)\beta$. Because
$v=\langle e,\,\cdot\,\rangle=-(\ad^* e)\alpha$ and
$[e,[e,\hat{\gt t}^{\perp}]]=0$, we have that
$\big(\Ad^*(\exp\,\frac{\lambda}{\mu}\,\ad
e)\big)\beta=\gamma+\delta+\mu\alpha,$ where
 $\delta$ is a nonzero linear function
 on $\g_{\tilde{e}}$
 which vanishes on $\gt t\oplus\hat{\gt
 t}^{\perp}\oplus\g(1)\oplus\g(2)$ and has the property that
$$
\delta(h_0)=\,\frac{\lambda^2}{2\mu}\,v([h_0,e]).
$$
Thus, $x+\omega(x)/\tilde e$ is zero on $\gt s$ for all
$x\in\hat{\gt t}^{\perp}$, proving (a). If $x\in\gt t$, then
$(x+\omega(x)/\tilde e)(\beta)=x(\gamma)=x(\beta)$, hence (b).
Finally, $(h_0+\omega(h_0)/\tilde e)(\beta)$ is a nonzero multiple
of $\lambda^2/\mu$, showing that the restriction of
$h_0+\omega(h_0)/\tilde e$ to $\mathfrak s$ is a nonzero multiple of
$(e^{+}-e^{-})^2/\tilde e$. One should keep in mind here that
$\psi(e^{+}-e^{-})=0$ for all $\psi\in\gt t^*\oplus\mathbb K\alpha$
and $v(e^{+}-e^{-})\ne 0$.
\end{proof}

For $1\le i\le l-1$, set $\varphi_i:={H_i}_{|\hat{\gt t}^*}$.
Then $\varphi_i$ is
homogeneous element in $\cS(\hat{\gt t})^{W'}$. It can be presented
uniquely as
$$
\varphi_i\,=\,\sum\limits_{j=0}^{\mu}\, \varphi_{i}^{(j)}\,
h_{0}^{j} \qquad \qquad\ \ \Big(\varphi_{i}^{(j)}\in\cS(\gt
t)^{W'_0},\, \,\,\,\, \varphi_i^{(\mu)}\ne
0,\,\,\,\,\,\mu=\mu(\varphi_i)\Big).
$$
Recall that $h_0$ spans the orthogonal complement to $\gt t$ in
$\hat{\gt t}$ with respect to the Killing form.
\begin{cl} For $1\le i\le l-1$ set $\mu_i=\mu(\varphi_i)$. Then in the above notation we
have
$$
\omega(H_i)_{|\gt t^*\oplus\,\mathbb
Kv}\,=\,a^{\mu_i}\,\varphi_i^{(\mu_i)}(e^+-e^-)^{2\mu_i} \qquad\ \
\, (1\le i\le l-1).
$$
\end{cl}
\begin{proof} This follows from Lemmas~\ref{E8-2} and \ref{E8-3}.
\end{proof}

Summing up the material of this subsection we obtain the following
result:

\begin{thm}\label{E8-4}
The algebra ${\mathcal S}(\gt g_{\tilde e})^{\gt g_{\tilde e}}$ is
free if and only if there is a homogeneous generating system
$\varphi_1,\ldots,\varphi_{7}$ in $\cS(\hat{\gt t})^{W'}$ such that
the elements
$\varphi_1^{(\mu_1)}h_0^{\mu_1},\ldots,\,\varphi_7^{(\mu_7)}h_0^{\mu_7}$
are algebraically independent.
\end{thm}

In type ${\bf E}_7$ it is difficult to calculate Weyl invariants by
hand, and the system of basic invariants used in the final part of
(4.7) is not very helpful in the present situation. Since this paper
is already quite long, we leave the ${\bf E}_8$ case open for the
time being.

\smallskip

\noindent 4.9. Assume now that $\g$ is not of type ${\bf A}_n$ or
${\bf E}_8$. Let $\tilde{e}$ be as before and put ${\gt
p}:=\n_{\g}(\mathbb K\tilde{e}).$ Recall that
$\p=\,\g(0)\oplus\g(1)\oplus\g(2)$ is a parabolic subalgebra of
$\g$. We are now going to apply our results on
$\cS(\g_{\tilde{e}})^{\g_{\tilde{e}}}$ to prove that the semi-centre
of the universal enveloping algebra $U(\gt p)$ is a polynomial
algebra. This will confirm a conjecture of Joseph for the parabolic
subalgebra $\p$.
\begin{cl}\label{semi}
Under the above assumptions, the semi-centre $U(\gt p)^{[\gt p,\gt
p]}$ is a polynomial algebra in $l=\rk \g$ variables.
\end{cl}
\begin{proof}
Since $\g$ is not of type $\mathbf A$, we have
$[\p,\p]=\g_{\tilde{e}}$ and $\p=\mathbb K\tilde{h}\oplus
\g_{\tilde{e}}$. Let $v\in \cS(\p)^{[\p,\p]}$ and write
$v=\tilde{h}^kv_k+\tilde{h}^{k-1}v_{k-1}+\cdots +v_0$ with
$v_i\in\cS(\g_{\tilde{e}})$. Since $\tilde{e}\in\z(\g_{\tilde{e}})$
and $\tilde{e}{\cdot}\tilde{h}^i=\,-2i\tilde{h}^{i-1}\tilde{e}$ for
all $i>0$, we get
$0=\,\tilde{e}{\cdot}v=-\textstyle{\sum}_{i=1}^k\,2i\tilde{h}^{i-1}\tilde{e}v_i.$
This yields $\cS(\p)^{[\p,\p]}=\,\cS(\g_{\tilde{e}})^{[\p,\p]}=\,
\cS(\g_{\tilde{e}})^{\g_{\tilde{e}}}$. Arguing in a similar fashion
we obtain
$$U(\p)^{[\p,\p]}=\,U(\g_{\tilde{e}})^{\g_{\tilde{e}}}=\,Z(\g_{\tilde{e}}),$$
where $Z(\g_{\tilde{e}})$ stands for the centre of
$U(\g_{\tilde{e}})$. As $\cS(\g_{\tilde{e}})^{\g_{\tilde{e}}}$ is a
polynomial algebra in $l$ variables, there exist algebraically
independent homogeneous elements $v_1,\ldots, v_l\in\cS(\p)$ such
that $\cS(\p)^{[\p,\p]}=\,\mathbb K[v_1,\ldots, v_l]$.

Let $r_i=\deg v_i$, where $1\le i\le l$, and let $(U_k)_{k\ge 0}$
denote the standard filtration of $U(\p)$. Using the
symmetrisation map $\cS(\p)\stackrel{\sim}{\to}U(\p)$ it is easy
to observe that there exist $u_1,\ldots, u_l\in U(\p)^{[\p,\p]}$
such that $u_i\in U_{r_i}$ and ${\rm gr}_{r_i}(u_i)=v_i$ for all
$i$. Since the $u_i$'s are central in $U(\g_{\tilde{e}})$, the
standard filtered-graded techniques now shows that
$U(\p)^{[\p,\p]}=\mathbb K[u_1,\ldots, u_l]$ is a polynomial
algebra in $l$ variables.
\end{proof}

\section{The null-cones in type A}

\noindent 5.1. In this section we assume that
$\g=\mathfrak{gl}(\VV)$ where $\dim\VV\ge 2$. Our goal is to prove
that for every $e\in\cN(\g)$ the null-cone $\gt \cN(e)\subset
\g_e^*$ has the expected codimension, i.e., $\dim \cN(e)=\dim \gt
g_e-n$. According to Theorem~\ref{deg1}, the variety $\cN(e)$ is the
zero locus of ${^e\!}F_1, \ldots, {^e\!}F_n,$ where
$F_i=\kappa^{-1}(\Delta_i)$. Thanks to the Affine Dimension Theorem,
in order to compute $\dim \cN(e)$ it suffices to find an
$n$-dimensional subspace $W\subset \gt g_e^*$ such that $W\cap
\cN(e)=0$. This will be achieved in a somewhat roundabout way: first
we shall construct a larger subspace $V^*\subset \gt g_e^*$ for
which the restrictions ${{^e\!}F_i}_{|V^*}$ can be described more or
less explicitly and then show that $V^*$ contains an $n$-dimensional
subspace transversal to $\cN(e)$.

For $m\in\{1,\ldots,k\}$, we partition the set $\{1,\ldots, m\}$
into pairs $(j,m-j+1)$. If $m$ is odd, then there will be a
``singular pair'' in the middle consisting of the singleton
$\{(m+1)/2\}$. We denote by $V_m$ the subspace of $\gt g_e$ spanned
by all $\xi_i^{j,s}$ with $i+j=m+1$, and set $V:=\bigoplus_{m\ge 1}
V_m$. Using the basis $\{(\xi_i^{j,s})^*\}$ of $\gt g_e^*$ dual to
the basis $\{\xi_i^{j,s}\}$, we shall regard the dual spaces $V_i^*$
and $V^*$ as subspaces of $\gt g_e^*$. Since $\mathbb
K[V^*]\cong{\mathcal S}(V),$ the restrictions
$\hat\varphi_i:={{{^e\!}F}_i}_{|V^*}$ are elements of ${\mathcal
S}(V)$. For $\bar s:=(s_1,\ldots,s_k)$ with $s_i\in \mathbb Z_{\ge
0}$ we set $|\bar s|:=s_1{+}s_2{+}\ldots{+}s_k$.

\begin{lm}\label{null-cone0}
Suppose $0\le q\le d_k$. Then $\hat\varphi_{n-q}\in{\mathcal S}(V_k)$. More precisely,
$$
\hat\varphi_{n-q}= \sum\limits_{|\bar s|=q}a(\bar
s)\,\xi_1^{k,d_k-s_k}\xi_{2}^{k-1,d_{k-1}-s_{k-1}} \cdots\,
\xi_{k}^{1,d_1-s_1} \ \ \mbox{ for some  }\   a(\bar s)\in\mathbb
K^{^\times}.$$
\end{lm}
\begin{proof} (a) According to Lemma~\ref{description},
${^e\!}F_{n-q}$ is a sum of monomials
$\xi_1^{\sigma(1),t_1}\ldots\xi_k^{\sigma(k),t_k}$, where $\sigma$
is a permutation of $\{1,\ldots,k\}$ and $t_1,\ldots,t_k$ are
nonnegative integers. Such a monomial does not vanish on $V^*$ only
if $\sigma(k)=1$, $\,\sigma(k-1)\in\{1,2\}\,$ and $\sigma(j)\le
k+1-j$ for all $j\le k$. Since $\sigma$ is a permutation, we then
have $\sigma(k-1)=2,\,$ $\sigma(k-2)=3\,$ and, in general,
$\sigma(j)=k+1-j$.

From (\ref{ksi}) we see that
$\xi_j^{k-j+1,\,d_{k-j+1}-s_{k-j+1}}\,\xi_{k-j+1}^{j,\,d_j-s_j}$ has
weight $\,2({d_j+d_{k-j+1}-s_j-s_{k-j+1}})$ with respect to $\ad h$.
As a consequence, the $h$-weight of
$$
\xi_1^{k,d_k-s_k}\xi_{2}^{k-1,d_{k-1}-s_{k-1}} \cdots\,
\xi_{k}^{1,d_1-s_1}
$$
equals $2({n-k-|\bar s|})$. Since $\deg F_{n-q}=n-q$ and
$\deg{{^e\!}F_{n-q}}=k$, this implies that only monomials with
$|\bar s|=q$ can occur in ${^e\!}F_{n-q}$. Because $|\bar s|=q\le
d_k\le d_i$ and all $s_i$ are nonnegative, we have that $s_i\le d_j$
for all $i,j$. This means that every $\xi_{k-i+1}^{i,d_i-s_i}$ is a
nonzero element of $\gt g_e$.

\smallskip

\noindent (b) We now prove by induction on $k$ that every $a(\bar
s)$ is nonzero. If $k=1$, then $V=\gt g_e$,
$\hat\varphi_{n-q}={{^e\!}F}_{n-q}=a(q)\xi_1^{1,d_1-q}$; and clearly
$a(q)\ne 0$. If  $k=2$, then $\big(\text{ad}^*
\xi_1^{1,1}\big){\cdot}V^*\subset V^*$. From this it follows that
the Poisson bracket $\{\xi_1^{1,1},\hat\varphi_{n-q}\}$ is zero. On
the other hand,
$$
\{\xi_1^{1,1},\,\hat\varphi_{n-q}\}\, =\, \,
\sum\limits_{i=0}^{q-1}\,\big(
a(q-i,i)-a(q-i-1,i+1)\big)\xi_1^{2,d_2-i}\xi_2^{1,d_1-q+i+1}\,.
$$
As the monomials $\xi_1^{2,d_2-i}\xi_2^{1,d_1-q+i+1}$ with $0\le
i\le q-1$ are nonzero in $\cS(\g_e)$, all coefficients $a(\bar s)$
with $|\bar s|=q$ must be equal. If one of them is zero, then all
are zeros. Assume that this is the case. By Lemma~\ref{description},
we then have
$$
{^e\!}F_{n-q}\,=\,\,\sum\limits_{|\bar s|=q}b(\bar
s)\,\xi_1^{1,d_1-s_1}\xi_2^{2,d_2-s_2}\,, \enskip \text{ where }
\,\, b(\bar s)\in\mathbb K.
$$
Let $s_2$ be the largest integer with $b(s_1,s_2)\ne 0$. As $s_2\le
q\le d_2$, the element $\xi:=\xi_2^{1,d_1-d_2+s_2}$ is nonzero in
$\gt g_e$. As ${{^e\!}F}_{n-q}$ belongs to the Poisson centre of
$\cS(\g_e)$, we have $\{\xi,\,{{^e\!}F}_{n-q}\}=0$. On the other
hand,
$$
\{\xi,\,{{^e\!}F}_{n-q}\}=\,b(s_1,s_2)\,\xi_1^{1,s_1}\xi_2^{1,d_1}+\,
(\text{multiples of monomials of the form $\xi_2^{1,*}\xi_2^{2,*}$
}).
$$
Since the RHS is nonzero, we reach a contradiction, proving the
lemma in case $k=2$.

\smallskip

\noindent (c) Now suppose $k>2$, and set $\gt
g':=\gt{gl}\big(\VV[1]{\,\oplus}\VV[k]\big)$ and $\gt
g'':=\gt{gl}\big(\VV[2]{\,\oplus}\cdots{\oplus}\VV[k{-}1]\big)$. These
are Lie subalgebras of $\g$ (embedded diagonally), and $e=e'+e''$
where $e'$ and $e''$ are the restrictions of $e$ to the $e$-stable
subspaces $\VV[1]{\,\oplus}\VV[k]\,$ and
$\VV[2]{\,\oplus}\cdots{\oplus}\VV[k{-}1]$.

We adopt the notation introduced in the course of proving
Lemma~\ref{description} and express $F_{n-q}$ as a polynomial in the
variables $E_{ij}$. Let $T$ be a monomial of $F_{n-q}$ such that
$T_{|V^*}$ is a nonzero multiple of a monomial of degree $k$ in
$\xi_{k-j+1}^{j,d_j-s_j}$. Then $T=T'T''$, where $T'$ and $T''$ are
polynomials in the variables coming from $\gt g'$ and $\gt g''$,
respectively. Suppose the restriction of $T'$ to $V^*$ equals
$a'\xi_1^{k,d_k-s_k}\xi_k^{1,d_1-s}$, where $a'\in\mathbb
K^{^\times}$. Then $T'$ is a monomial of $F'_{p'}\in{\mathcal S}(\gt
g')^{\gt g'}$ for $p'=d_1+d_2+2-s_1-s_2$. Likewise, $T''$ is a
monomial of $F_{p''}\in{\mathcal S}(\gt g'')^{\gt g''}$ for
$p''=n-q-p'$. It follows that $a(\bar
s)=a(s_1,s_k)a(s_2,\ldots,s_{k-1})$ where the coefficients
$a(s_1,s_k)$ and $a(s_2,\ldots,a_{k-1})$ are related to the
nilpotent elements $e'\in\gt g'$ and $e''\in\gt g''$, respectively.

Note that $e'\in\gt g'$ has two Jordan blocks of sizes $d_1+1$ and
$d_2+1$, and $a(s_1,s_k)$ is the coefficient of
$\xi_1^{k,d_k-s_k}\xi_k^{1,d_1-s_1}$ in the expression for
$\hat\varphi_{p'}$. This coefficient is nonzero by part~(b). The
coefficient $a(s_2,\ldots,s_{k-1})$ arises in a similar way from the
nilpotent element $e''\in\g''$. Since
$\g''\cong\mathfrak{gl}_{n-d_1-d_k-2}$ we can apply the inductive
hypothesis to conclude that $a(s_2,\ldots,s_{k-1})\ne 0$. Therefore
every $a(\bar s)$ is nonzero, as wanted.
\end{proof}

\noindent 5.2. Our next goal is to describe the zero locus
$X=X^{(d_k)}$ of $\hat\varphi_n,
\hat\varphi_{n-1},\ldots,\hat\varphi_{n-d_k}$ in $V_k^*$. Denote by
$X_{\bar s}$ the subspace of $V_k^*$ consisting of all $\gamma\in
V_k^*$ such that $\,\xi_{k-i+1}^{i,d_i-t}(\gamma)=0$ for $0\le t<
s_i$. Let ${\bf e}_i$ be the $k$-tuple whose $i$-th component equals
$1$ and the other components are zero.

\begin{lm}\label{null-cone1} The variety $X$ is a union of subspaces. More precisely,
$X\,=\,\bigcup_{|\bar s|=d_k+1}\,X_{\bar s}$.
\end{lm}
\begin{proof} Let $X^{(q)}\subset V_k^*$ be
the zero locus of
$\hat\varphi_n,\hat\varphi_{n-1},\ldots,\hat\varphi_{n-q}$. We are
going to prove by induction on $q$ that $X^{(q)}$ is a union of
subspaces in $V_k^*$ and the irreducible components of $X^{(q)}$
correspond bijectively to the $k$-tuples $\bar s$ with $|\bar
s|=q+1$. When $q=0$, our set of functions is a singleton containing
$\hat\varphi_n=
\xi_{k}^{1,d_1}\xi_{k-1}^{2,d_2}\cdots\,\xi_{1}^{k,d_k}$. Therefore,
$X^{(0)}$ is the union of $k$ hyperplanes in $V_k^*$ defined by the
equations $\,\,\xi_{k-i+1}^{i,d_i}=0,$ where $1\le i\le k$.

Assume that $X^{(q-1)}$ is a union of subspaces of $V_k^*$
parametrised by the $k$-tuples of size $q$. Let $\bar{s}$ be a
$k$-tuple of size $q-1$ and let $X_{\bar s}$ be the irreducible
component of $X^{(q-1)}$ corresponding to $\bar{s}$. Now consider an
arbitrary monomial
$f:=\,\xi_1^{k,d_k-t_k}\xi_{2}^{k-1,d_{k-1}-t_{k-1}}\cdots\,
\xi_{k}^{1,d_1-t_1}$ with $\sum t_i=q$, i.e., a typical summand of
$\hat\varphi_{n-q}$. If $\bar t=(t_1,\ldots, t_k)\ne \bar s$, then
there exists an index $i$ such that $t_i<s_i$. But then
$\xi_{k-i+1}^{i,d_i-t_i}$, and hence $f$, vanishes on $X_{\bar s}$.
This shows that the restriction of $\hat\varphi_{n-q}$ to $X_{\bar
s}$ coincides, up to a nonzero multiple, with that of
$\xi_1^{k,d_k-s_k}\xi_{2}^{k-1,d_{k-1}-s_{k-1}}\cdots\,
\xi_{k}^{1,d_1-s_1}$. As a consequence, the zero locus of
$\hat\varphi_n,\hat\varphi_{n-1},\ldots,\hat\varphi_{n-q}$ in
$X_{\bar s}$ is the union of $k$ linear subspaces $X_{\bar{s}+{\bf
e}_i}$, where $1\le i\le k$. Then $X^{(q)}=\bigcup_{|\bar s|=q+1}
X_{\bar s},\,$ and the statement follows by induction on $q$.
\end{proof}

\noindent 5.3. By Lemma~\ref{null-cone1}, all irreducible components
of the variety $X^{(d_k)}\subset V_k^*$ have dimension equal to
$\dim V_k-(d_k+1)$. Hence there is a linear subspace $W_k\subset
V_k^*$ such that $\,\dim W_k=d_k+1$ and $W_k\cap X^{(d_k)}=0$.

\begin{prop}    \label{null-cone2}
There exists an $n$-dimensional linear subspace $W=\bigoplus_{m\ge
1} W_m$ in $V^*$ such that $W_m\subset V_m^*$ for all $m$ and
$W\cap\cN(e)=0$.
\end{prop}
\begin{proof} We argue by induction on $k$.
If $k=1$, then $\cN(e)=0$ and there is nothing to prove. So assume
that $k\ge 2$, and set $\gt g_k:=\gt{gl}\big(\VV[k]\big)$ and
$\bar{\g}:=\gt{gl}\big(\VV[1]{\,\oplus}\cdots{\oplus}\VV[k{-}1]\big)$.
These Lie algebras are embedded diagonally into $\g$, and we regard
the dual spaces $\bar{\g}^*$ and $\g_k^*$ as subspaces of $\g^*$.
Note that $e=e_k+\bar{e}$ where $e_k$ and $\bar{e}$ are the
restrictions of $e$ to  $\VV[k]$ and
$\VV[1]{\,\oplus}\cdots{\oplus}\VV[k{-}1]$, respectively. Clearly, $e_k$
is a regular nilpotent element in $\g_k\cong\gt{gl}_{d_k+1}$ and
$\bar{e}\in\bar{\g}\cong\mathfrak{gl}_{n-d_k-1}$ is a nilpotent
element with Jordan blocks of sizes $d_1+1,\ldots, d_{k-1}+1$. For
$1\le i\le n-d_k-1$, put $\bar{F}_i:={F_i}_{|\bar{\gt g}^*}.\,$
Restricting the principal minors $\Delta_i$ from $\g$ to $\bar{\g}$
it is easy to see that the homogeneous generating system
$\bar{F}_i,\ldots,\bar{F}_{n-d_k-1}$ of $\cS(\bar{\gt
g})^{\bar{\g}}$ is good for $\bar{e}\in\bar{\g}$.

Next we observe that $\bar{\gt g}_{\bar{e}}$ is a Lie subalgebra of
$\gt g_e$ spanned by all $\xi_i^{j,s}$ with $1\le i,j<k$. Hence we
may identify the dual space $(\bar{\gt g}_{\bar{e}})^*$ with the
linear span of $\big\{(\xi_i^{j,s})^*\,|\,\, 1\le i,j<k\big\}$ in
$\g_e^*$. For every $i\in\{1,\ldots,n-d_k-1\}$ the restriction of
$\,{^e\!}F_i\,$ to $(\bar{\gt g}_{\bar{e}})^*$ equals
$\,^{{\bar{e}}\!}\bar{F}_i$.

Note that $V_m^*\subset(\bar{\gt g}_{\bar{e}})^*$ for $m<k$ and
$V_k^*\cap(\bar{\gt g}_{\bar{e}})^*=0$. By our inductive hypothesis,
there exists a subspace $\overline{W}=\,\bigoplus_{m=1}^{k-1}\, W_i$
such that $\dim\overline{W}=n-d_k-1$ and
$\overline{W}\cap\cN(\bar{e})=0$. Choose a $(d_k+1)$-dimensional
subspace $W_k$ in $V_k^*$ with $W_k\cap X^{(d_k)}=0$. Such a
subspace exists by Lemma~\ref{null-cone1}. Now set
$W:=\overline{W}\oplus W_k$. Then $\dim W=n$.

We claim that $W\cap\cN(e)=0$. By Lemma~\ref{null-cone0}, for
$n-d_k\le i\le n$ the restriction
$\hat\varphi_i={{{^e\!}F}_i}_{|V^*}$ belongs to ${\mathcal S}(V_k)$.
Therefore, the zero locus of
$\hat\varphi_n,\ldots,\hat\varphi_{n-d_k}$ in $V^*$ coincides with
$\big(\bigoplus_{m=1}^{k-1}\, V_m^*\big)\times X^{(d_k)}$. Since
$W_k\cap X^{(d_k)}=0$, we obtain
$W\cap\cN(e)\subset\bigoplus_{m=1}^{k-1}\, V_m^* \subset (\bar{\gt
g}_{\bar{e}})^*.$ But then $\,W\cap\cN(e)\subset
\overline{W}\cap\cN(\bar{e})=0,\,$ and we are done.
\end{proof}

The following is the main result of this section:

\begin{thm}\label{null-cone3}
Let $e$ be an arbitrary nilpotent element in $\gt g=\gt{gl}_n$. Then
all irreducible components of the null-cone $\cN(e)$ have
codimension $n$ in $\g_e^*$ and hence
${{^e\!}F}_1,\ldots,{{^e\!}F}_n$ is a regular sequence in $\mathbb
\cS(\gt g_e)$.
\end{thm}

\noindent 5.4. Let $X\subset\mathbb A^d_{\mathbb K}$ be a Zariski
closed set and let $x=(x_1,\ldots,x_d)$ be a point of $X$. Let $I$
denote the defining ideal of $X$ in the coordinate algebra
${\mathcal A}=\mathbb K[X_1,\ldots, X_d]$ of $\mathbb A^d_{\mathbb
K}$. Each nonzero $f\in\mathcal A$ can be expressed as a polynomial
in $X_1-x_1,\ldots,X_d-x_d$, say $f=f_k+f_{k+1}+\cdots$, where $f_i$
is a homogeneous polynomial of degree $i$ in $X_1-x_1,\ldots,
X_d-x_d$ and $f_k\ne 0$. We set ${\rm in}_x(f):=f_k$ and denote by
${\rm in}_x(I)$ the linear span of all ${\rm in}_x(f)$ with $f\in
I\setminus\{0\}$. This is an ideal of $\mathcal A$, and the affine
scheme $TC_x(X):=\,{\rm Spec}\,{\mathcal A}/{\rm in}_x(I)$ is called
the {\it tangent cone} to $X$ at $x$. Note that $(I\cap {\mathfrak
m}_x^k)_{k\ge 0}$ is a descending filtration of $I$, and the scheme
$TC_x(X)$ is nothing but the prime spectrum of the graded algebra
${\rm gr}_{{\mathfrak m}_x}{\mathcal A}/{\rm gr}\,I$. It is
well-known that the projectivised tangent cone ${\mathbb
P}TC_x(X)\subset{\mathbb P}T_x(X)$ is isomorphic to the special
divisor of the blow-up of $X$ at $x$; see \cite[Ex.~IV-24]{EH} for
example. Consequently, for $X$ irreducible, all irreducible
components of $TC_x(X)$ have dimension equal to $\dim X$.
\begin{cl}\label{cone}
Let $\cN$ be the nilpotent cone of $\g=\mathfrak{gl}_n$ and
$F_i=\kappa^{-1}(\Delta_i)$ where $1\le i\le n$. Let $e\in\mathcal
N$ and $r=\dim\g_e$. Then $TC_e({\mathcal N})\cong\,{\mathbb
A}^{n^2-r}_{\mathbb K}\!\times\, {\rm
Spec}\,\cS(\g_e)/({^{e\!}F_1},\ldots,{^{e\!}F_n})$ as affine
schemes.
\end{cl}
\begin{proof}
Since the map $x\mapsto (x,\,\cdot\,)$ takes $e$ to $\chi$ and $\cN$
isomorphically onto the zero locus of the ideal $J=(F_1,\ldots,
F_n)\subset\cS(\g)$, the scheme $TC_e(\cN)$ is isomorphic to ${\rm
Spec}\,\cS(\g)/{\rm in}_\chi(J)$. As $\chi(f)=1$, we have $\gt
g=\mathbb Kf\oplus e^{\perp}$ where $e^{\perp}$ is the orthogonal
complement to $\mathbb Ke$ in $\g$. For $1\le i\le n$ write
$F_i=f^{k(i)} p_{0,i}+f^{k(i)-1}p_{1,i}+\cdots+p_{k(i),i},$ where
$p_{j,i}\in\cS(e^{\perp})$ and $p_{0,i}\ne 0$. According to
Corollary~\ref{eF}, we have $p_{0,i}={^{e\!}F}_i$. Since $e^\perp$
and $f-\chi(f)$ lie in the maximal ideal of $\chi$ in
$\mathbb K[\g^*]=\cS(\g)$, it follows that ${\rm in}_\chi(F_i)={^{e\!}F}_i$
for all $1\le i\le n$.

By Theorem~\ref{null-cone3}, ${^{e\!}F}_1,\ldots,{^{e\!}F}_n$ is a
regular sequence in $\cS(\g_e)$. Therefore, it is also a regular
sequence in $\cS(\g)$. Since $J=(F_1,\ldots, F_n)$, it follows that
the ideal ${\rm in}_\chi(J)$ is generated by
${^{e\!}F}_1,\ldots,{^{e\!}F}_n$;\, see \cite[Prop.~2.1]{VV}. As a
consequence, $$TC_e(\cN)\cong\,{\rm
Spec}\,\cS(\g)/({^{e\!}F}_1,\ldots,{^{e\!}F}_n)\cong\,{\rm
Ann}(\g_e)\times{\rm
Spec}\,\cS(\g_e)/({^{e\!}F_1},\ldots,{^{e\!}F_n})$$ as affine
schemes. Since $\dim {\rm Ann}(\g_e)=n^2-r$, the result follows.
\end{proof}

\begin{conj} \label{red} If $\g=\mathfrak{gl}_n$, then for any
$e\in \cN$ the scheme $TC_e(\cN)$ is reduced.
\end{conj}

\begin{rmk}\quad\\
1.  It can be shown that in the subregular ${\mathbf G}_2$ case
the variety $TC_e(\cN(\g))_{\rm red}$ is isomorphic to an affine
space, but the scheme $TC_e(\cN(\g))$ is {\it not} reduced. Thus,
one cannot expect Conjecture~\ref{red} to be true for any simple
Lie algebra.

\smallskip

\noindent 2. It follows from Corollary~\ref{cone} that for
$\g=\mathfrak{gl}_n$ the affine variety $TC_e({\cN(\g)})_{\rm
red}$ is isomorphic to ${\mathbb A}^{m}_{\mathbb K}\times\,\cN(e)$
where $m=\dim\g-\dim\g_e$. It is possible that this isomorphism
continues to hold for any reductive Lie algebra $\g$. If this is
the case, then the variety $\cN(e)$ is always equidimensional.

\smallskip
\noindent 3. Although the variety $\cN(e)$ is irreducible in some
cases, in general it has many irreducible components. Due to
Theorem~\ref{summa}(iii), in order to prove Conjecture~\ref{red}
it would be sufficient to show that every irreducible component of
$\cN(e)$ intersects with $(\g_e^*)_{\rm reg}$. Describing the
irreducible components of $\cN(e)$ for $\g=\mathfrak{gl}_n$
appears to be an interesting combinatorial problem.
\end{rmk}

\section{Miscellany}  \label{sect:misc}

\noindent 6.1. In this section, Conjecture~\ref{conj:main} will be
verified in some special cases. The idea is that, for some
$e\in\cN(\g)$, we can prove that the algebra $\mathbb
K[\g_e]^{\g_e}$ is graded polynomial. If, in addition, it is known
that $\g_e\simeq \g_e^*$ as $\g_e$-modules, then we conclude that
Conjecture~\ref{conj:main} holds for such $e$.

We briefly recall the structure of the centralizer $\g_e$ of a
nilpotent element $e\in\g$ as described by the Dynkin--Kostant
theory; see e.g. \cite[Ch.\,4]{CM}. Let $\{e,h,f\}$ be an
$\tri$-triple and $\g=\bigoplus_{i\in\BZ}\g(i)$ the corresponding
$\BZ$-grading. Then $\g_e=\bigoplus_{i\ge 0}\g_e(i)$ and $\g_e(0)$
is a maximal reductive subalgebra of $\g_e$. Moreover,
$\g_e(0)=\z_\g(e,f)=\z_\g(e,h,f)$.
The element $e$ is called {\it even} if all the eigenvalues of $\ad
h$ are even, i.e., if $\g(i)=0$ for $i$ odd. By a classical result
of Dynkin, $e$ is even if and only if $\g(1)=0$; see
\cite[Thm.\,8.3]{ebd}. In this case the weighted Dynkin diagram of
$e$ contains only labels 0 and 2.

In the following theorem, we use some concepts and results on (1)
semi-direct products of Lie algebras and (2) contractions of Lie
algebras. All the necessary definitions can be found in
\cite[Sect.\,4]{dp} and \cite[Ch.\,7]{t41}, respectively.

\begin{thm}  \label{alg-princ}
Suppose that a principal nilpotent element in $\g_e(0)$ is also
principal in $\g(0)$ and $e$ is even. Then $\mathbb K[\g_e]^{\g_e}$
is a polynomial algebra and the degrees of basic invariants (=\,free
homogeneous generators) are the same as those for $\mathbb
K[\g(0)]^{\g(0)}$.
\end{thm}\begin{proof}
Associated to the triple $(e,h,f)$ and the corresponding
$\BZ$-grading, we have three Lie algebras: $\g(0)$, $\g_e$, and
$\q:=\g_e(0)\ltimes \big(\bigoplus_{i\ge 2}\,\g_e(i)\big)$. Here the
sign $\ltimes$ refers to the semi-direct product of Lie algebras and
the space $\bigoplus_{i\ge 2}\,\g_e(i)$ in $\q$ is regarded as
commutative Lie algebra. Clearly, $\dim\q=\dim\g_e$. The equality
$\dim\g(0)=\dim\g_e$ is equivalent to the fact that $e$ is even.
Thus, all three Lie algebras have the same dimension. Here we obtain
the chain of Lie algebra contractions:
\[
   \g(0) \leadsto \g_e \leadsto \q \ .
\]
The first contraction can be described as follows. Consider the
curve $e(t):=e+tf\in \g$, $t\in\mathbb K$. For $t\ne 0$, the element
$e(t)$ is $G$-conjugate to $h$. Therefore, $\g_{e(t)}$ is isomorphic
to $\g_h=\g(0)$. Hence $\lim_{t\to 0}\g_{e(t)}=\g_e$ yields a
contraction of $\g(0)$ to $\g_e$. Using the terminology of
\cite[Sect.\,9]{dp}, one can say that the passage $\g_e\leadsto \q$
is an {\it isotropy contraction\/} of $\g_e$. By
\cite[Theorem\,6.2]{dp}, the algebra of invariants of the adjoint
representation of $\q$ is polynomial. Moreover, if a regular
nilpotent element of $\g_e(0)$ is also regular in $\g(0)$, then by
\cite[Theorem\,9.5]{dp} the invariant algebras $\mathbb
K[\g(0)]^{\g(0)}$ and $\mathbb K[\q]^{\q}$ have the same Krull
dimension and the same degrees of basic invariants. It is easily
seen that the algebra of invariants of the adjoint representation
can only become larger under contractions. Since $\mathbb
K[\g(0)]^{\g(0)}$ and $\mathbb K[\q]^{\q}$ appear to be "the same",
the intermediate algebra $\mathbb K[\g_e]^{\g_e}$ must also be
polynomial with the same degrees of basic invariants.
\end{proof}
\noindent 6.2. By a result of Elashvili--Panyushev (Appendix to
\cite{vitya}), the assumptions on $e$ in Theorem~\ref{alg-princ}
precisely mean that $e$ is a member of a {\it rectangular principal
nilpotent pair}. The general theory of principal nilpotent pairs (to
be abbreviated as {\it {pn}-pairs} from now) was developed by Victor
Ginzburg \cite{vitya}. Because the general notion is not needed
here, we only recall the definition of a rectangular {pn}-pair.
\begin{df}
A pair of nilpotent elements $\mathbf e=(e_1,e_2)$ is called a {\it
rectangular pn-pair\/} if $\dim (\g_{e_1}\cap\g_{e_2})=\rk\g$ and
there are pairwise commuting $\tri$-triples $(e_1,h_1,f_1)$ and
$(e_2,h_2,f_2)$.
\end{df}
We say that a nilpotent orbit $G{\cdot}e$ is {\it very nice\/} if
$e$ is a member of a rectangular pn-pair and $\g_e\simeq\g_e^*$ as
$\g_e$-modules.
\begin{cl} Suppose $G{\cdot}e$ is very nice. Then
Conjecture~\ref{conj:main} holds for $\g_e$ and $\mathcal
S(\g)^{\g}$ admits a good generating system for $e$.
\end{cl}

 A  classification of rectangular {pn}-pairs is obtained by
Elashvili and Panyushev in \cite[Appendix]{vitya}. From that
classification one derives a description of very nice orbits. It is
worth mentioning that for a {pn}-pair $(e_1,e_2)$ the condition that
$G{\cdot}e_1$ is very nice does not in general guarantee that so is
$G{\cdot}e_2$.

Although there are not too many very nice nilpotent orbits
(especially in the exceptional Lie algebras), this approach does
provide new examples supporting Conjecture~\ref{conj:main}. The
examples for $\sln$ and $\spn$ are not new; see Section~4.

\smallskip

\noindent 6.3. Below we list the very nice nilpotent orbits in
exceptional Lie algebras. For each such orbit we give the
Dynkin-Bala-Carter label, the weighted Dynkin diagram, and the
degrees of basic invariants for $\mathbb \cS(\g_e)^{\g_e}$.

\bigskip

\begin{tabular}{ccrl}
\framebox{$\GR{E}{6}$} & \phantom{class} $\GR{D}{4}$ &
\begin{E6}{0}{0}{2}{0}{0}{2}\end{E6} &
\qquad 1,\,1,\,2,\,2,\,3,\,3 \\
\framebox{$\GR{E}{7}$} & \phantom{class} $\GR{E}{6}$ &
\begin{E7}{0}{2}{0}{2}{2}{2}{0}\end{E7} &
\qquad 1,\,1,\,1,\,1,\,2,\,2,\,2 \\
 & \phantom{class} $\GR{A}{2}{+}\GR{A}{4}$ & \begin{E7}{0}{0}{0}{2}{0}{0}{0}\end{E7} &
\qquad 1,\,2,\,2,\,2,\,3,\,3,\,4 \\
 & \phantom{class} $\GR{A}{6}$ & \begin{E7}{0}{2}{0}{2}{0}{0}{0}\end{E7} &
\qquad 1,\,1,\,2,\,2,\,2,\,2,\,3 \\
\end{tabular}

\bigskip

Let us give some details on the unique orbit for $\GR{E}{6}$. Here
$\dim\g_e=18$ and $\g_e$ is the direct sum of the 2-dimensional
centre and the Takiff Lie algebra $\es$ modelled on
$\mathfrak{sl}_3$. Namely, $\es$ is just the semi-direct product
$\mathfrak{sl}_3\ltimes\mathfrak{sl}_3$.

\smallskip

\noindent 6.4. The very nice nilpotent orbits in classical Lie
algebras are described below.

\smallskip

\noindent $1^\circ$.  $\g=\sln$. Here $e$ is a member of a
rectangular pn-pair if and only if the corresponding partition of
$n$ is a rectangle (i.e., all the parts are equal). That is, we may
assume that $n=rs$ and the partition of $e$ is $(r,\ldots, r)$, with
$s$ parts. We also write $e\sim (\underbrace{r,\dots,r}_{s})$ for
this. It is harmless but technically easier to work with $\g=\gln$
in place of $\sln$. Then $\g_e$ is a generalised Takiff Lie algebra
modelled on $\mathfrak{gl}_s$. More precisely, $\g_e$ is the factor
algebra of $\mathfrak{gl}_s\otimes \mathbb K[t]$ by its ideal
$\mathfrak{gl}_s\otimes t^r\mathbb K[t]$, where $t$ is an
indeterminate. It is easily seen that $\g_e\simeq \g_e^*$. (See
\cite{rt} and
 \cite[Sect.\,11]{dp} for more results on generalised Takiff Lie algebras.)
The second member of the rectangular pn-pair is given by the
conjugate partition $(s,\dots,s)$, with $r$ parts. This situation is
symmetric and both nilpotent orbits are very nice.

\smallskip

\noindent $2^\circ$.  $\g=\spn$. Here $e$ is a member of a
rectangular pn-pair if and only if the corresponding partition of
$2n$ is a rectangle whose sides have different parity. That is, we
may assume that $2n=rs$, where $r$ is even and $s$ is odd. The
situation here is not symmetric. Only the orbit corresponding to the
partition $(s,\dots,s)$ with $r$ parts is very nice.

\smallskip

\noindent $3^\circ$.  $\g=\son$. Here we have to distinguish the
series {\bf B} and {\bf D}.
\\
$\bullet$ \ \ If $n$ is odd, then the only suitable partitions are
the rectangles whose both sides are odd. That is, $n=rs$, where $r$
and $s$ are odd. Then $e\sim(\underbrace{s,\dots,s}_{r})$. Here both
members of the rectangular pn-pair give rise to very nice orbits.
\\
$\bullet$ \ \ For $n$ even, there are more possibilities for
rectangular pn-pairs.

(1) If a partition of $n$ is rectangle with both even sides, then
neither of the respective orbits is very nice.

(2) If $n=rs+1$, where $r,s$ are odd, the there is a rectangular
pn-pair $(e_1,e_2)$ with $e_1\sim (\underbrace{s,\dots,s}_{r},1)$
and $e_2\sim (\underbrace{r,\dots,r}_{s},1)$. Here both members of
the rectangular pn-pair give rise to very nice orbits.

(3) If $n=r+s$, where $r,s$ are odd, then there is a rectangular
pn-pair $(e_1,e_2)$ with $e_1\sim (s,\underbrace{1,\dots,1}_{r})$
and $e_2\sim (r,\underbrace{1,\dots,1}_{s})$. Here neither of the
respective orbits is very nice.
\section{Appendix}

Here we give an alternative (elementary) proof of
Proposition~\ref{invariants}, which is inspired by an unpublished
result of J.-Y.~Charbonnel (private communication).

Let $e^{\perp}\subset\gt g$ be the orthogonal complement of $\mathbb
Ke$. Since $(e,f)=1$, we have
 $\gt g=\mathbb Kf\oplus e^{\perp}$. Take a nonzero homogeneous $F\in\cS(\gt g)^G$ and
 express it as
$$
F=f^k p_0+f^{k-1}p_1+\cdots+p_k,
$$
where $p_i\in\cS(e^{\perp})$ and $p_0\ne 0$.

\begin{lm}\label{invariants2}
For any nonzero homogeneous $F\in\cS(\gt g)^G$ we have that
$p_0\in\cS(\gt g_e)^{G_e}$.
\end{lm}
\begin{proof}
If $g\in G_e$, then  $(\Ad\, g)e^{\perp}\subset e^{\perp}$ and $(\Ad
\,g)f\in f+e^{\perp}$. Therefore,
$$
F=g{\cdot}F=(g{\cdot}p_0)f^k+f^{k-1}p_1'+\cdots+p_k'$$ for some
$p_i'\in\cS(e^{\perp})$. Since $g{\cdot}p_0\in\cS(e^\perp)$, this
shows that $p_0$ is $G_e$-invariant.

Recall that $\gt g=\gt g_e\oplus\Im \ad f$. Choose a basis
$y_1,\ldots,y_{t}$ of $e^{\perp}\cap\Im \ad f$. If $p_0$ is not an
element of $\cS(\gt g_e)$, then renumbering the $y_i$'s if necessary
we may assume that $$p_0=y_1^sq_0+y_1^{s-1}q_1+\cdots+ q_s,$$ where
$s\ge 1$, $q_i\in{\mathcal S}(\gt g_e)[y_2,\ldots,y_{t}]$, and $q_0\ne 0$.
 Since the Killing form of $\g$
induces a nondegenerate pairing between $\Im\ad e$ and $\Im \ad f$,
there is a $z\in\gt g$ such that
 $([e,z],y_1)\ne 0$, $([e,z],f)=0$ and $([e,z],y_i)=0$ for all $i\ne 1$.
Note that $([z,y_1],e)=(y_1,[e,z])\ne 0$ and $([z,y_i],e)=0$ for all
$i\ne 1$. Also, $([z,x],e)=-(z,[e,x])=0$ for all $x\in\gt g_e$.

Rescaling $z$ if need be, we may assume that $([z,y_1],e)=1$. Then
$[z,y_1]\in f+e^\perp$ and $[z,f]\in e^\perp$, implying
$$
\{z,F\}=\,(sy_1^{s-1}q_0+(s-1)y_1^{s-2}q_1+\cdots+q_1)f^{k+1}+
\text{ (terms with smaller powers of $f$)}
$$
(here $\{\,\cdot,\,\cdot\,\}$ stands for the Poisson bracket of $
S(\g)$ induced by the Lie product in $\g$). This, however,
contradicts the equality $\{z,F\}=0$.
\end{proof}
\begin{cl}\label{eF}
For any homogeneous $F\in\cS(\g)^G$ we have that $\,^{e\!}F=p_0\,$
and $\,^{e\!}F\in\cS(\g_e)^{G_e}$.
\end{cl}
\begin{proof}
In view of  Lemma~\ref{invariants2} we have
$^{e}\!F=\kappa_e^{-1}(\kappa(p_0))=\kappa^{-1}(\kappa(p_0))=p_0$,
as stated.
\end{proof}
\noindent {\bf Acknowledgements.} Part of this work was done at the
Max Planck Institut f{\"u}r Mathematik (Bonn) and Manchester
Institute for Mathematical Sciences. We would like to thank both
institutions for hospitality and financial support. The third author
is a Humboldt fellow. She is grateful to the Humboldt Foundation for
financial support and to the Univesit\"at zu K\"oln, and especially
Peter Littelmann, for wonderful working conditions. We are also
thankful to the anonymous referee for very careful reading and
helpful comments.

\end{document}